\documentclass{article}
\usepackage[top = 1.5cm, left = 1.5cm, right = 1.5cm]{geometry}
\usepackage{fullpage}
\usepackage{microtype}
\usepackage{graphicx}
\usepackage{subfigure}
\usepackage{booktabs} 

\usepackage[utf8]{inputenc} 
\usepackage[T1]{fontenc}    

\usepackage{hyperref}
\usepackage{amsfonts}
\usepackage{amssymb,color,amsthm}
\usepackage{bbm}
\usepackage{mathrsfs}
\usepackage[none]{hyphenat}
\usepackage{subfigure}
\usepackage[tbtags]{amsmath}
\usepackage[authoryear,square]{natbib}
\usepackage{algorithm}
\usepackage{algpseudocode}

\def\remark{\addtocounter{remark}{1}\def\@currentlabel{\theremark}%
\emph{Remark~\theremark}. } \makeatother

\newcounter{remark}

\def\bd{\vec d}
\def\be{\mathbf e}
\def\bx{\mathbf x}
\def\bh{\mathbf h}
\def\tbh{\widetilde{\mathbf h}}
\def\wbx{\widetilde{\mathbf x}}
\def\wx{\widetilde{x}}

\def\bw{\mathbf w}

\def\by{\mathbf y}

\def\bu{\mathbf u}
\def\bv{\mathbf v}

\def\hbv{\widehat{\mathbf v}}

\def\bz{\mathbf z}
\def\tDelta{\widetilde{\Delta}}
\def\wde{\widehat{\delta}}

\def\bA{\mathbf A}

\def\bD{\mathbf D}

\def\bF{\mathbf F}
\def\bG{\mathbf G}
\def\bH{\boldsymbol{\mathcal{H}}}
\def\tbH{\widetilde{\boldsymbol{\mathcal{H}}}}
\def\bI{\mathbf I}

\def\bM{\mathbf M}

\def\bQ{\mathbf Q}

\def\bT{\mathbf T}
\def\bU{\mathbf U}

\def\bX{\mathbf X}
\def\bY{\mathbf Y}

\def\bZ{\mathbf Z}

\newcommand{\T}{\mathsf{T}}

\def\secref#1{Section~\ref{#1}}
\def\leref#1{Lemma~\ref{#1}}
\def\conref#1{Condition~\ref{#1}}
\def\thref#1{Theorem~\ref{#1}}

\def\coref#1{Corollary~\ref{#1}}
\def\figref#1{Figure~\ref{#1}}
\def\figtab#1{Table~\ref{#1}}
\def\algref#1{Algorithm~\ref{#1}}
\def\asref#1{Assumption~\ref{#1}}
\def\bydef{\triangleq}
\newtheorem{lemma}{Lemma}
\newtheorem{condition}{Condition}
\newtheorem{theorem}{Theorem}
\newtheorem{corollary}{Corollary}
\newtheorem{definition}{Definition}
\newtheorem{assumption}{Assumption}

\title{On the Sublinear Convergence of Randomly \\ Perturbed Alternating Gradient Descent \\  to Second Order Stationary Solutions}

\author{Songtao Lu \thanks{Department of Electrical and Computer Engineering, University of Minnesota -- Twin Cities}\;\thanks{Department of Electrical and Computer Engineering, Iowa State University}
\\
\texttt{\small lus@umn.edu}
\and
Mingyi Hong \footnotemark[1]
\\
\texttt{\small mhong@umn.edu}
\\
\and
Zhengdao Wang \footnotemark[2]
\\
\texttt{\small zhengdao@iastate.edu}
}
\date{}

\begin{document}
\maketitle

\begin{abstract}
The alternating gradient descent (AGD) is a simple but popular algorithm which has been applied to problems in optimization, machine learning, data ming, and signal processing, etc. The algorithm updates two blocks of variables in an alternating manner, in which a gradient step is taken on one block, while keeping the remaining block fixed. When the objective function is nonconvex, it is well-known the AGD converges to the first-order stationary solution with a global sublinear rate.

In this paper, we show that a variant of AGD-type algorithms will not be trapped by ``bad'' stationary solutions such as saddle points and local maximum points. In particular, we consider a smooth unconstrained optimization problem, and propose a perturbed AGD (PA-GD) which converges (with high probability) to the set of second-order stationary solutions (SS2) with a global sublinear rate. To the best of our knowledge, this is the first alternating type algorithm which takes $\mathcal{O}(\text{polylog}(d)/\epsilon^{7/3})$ iterations to achieve SS2 with high probability [where polylog$(d)$ is polynomial of the logarithm of dimension $d$ of the problem].
\end{abstract}

\section{Introduction}
\label{submission}

In this paper, we consider a smooth and unconstrained nonconvex optimization problem
\begin{equation}
\min_{\bx\in\mathbb{R}^{d\times 1}} f(\bx)\label{eq.op1}
\end{equation}
where $f:\mathbb{R}^{d}\rightarrow\mathbb{R}$ is twice differentiable.

There are many ways of solving problem \eqref{eq.op1}, such as gradient descent (GD), accelerated gradient descent (AGD), etc. When the problem dimension is large, it is natural to split the variables into multiple blocks and solve the subproblems with smaller size individually. The block coordinate descent (BCD) algorithm, and many of its variants such as block coordinate gradient descent (BCGD) and alternating gradient descent (AGD) \cite{bertsekas99,li17b},  are among the most powerful tools for solving large scale convex/nonconvex optimization problems \cite{nesterov2012efficiency,beck2013convergence,razaviyayn2013unified,hong13complexity}. The BCD-type algorithms
partition the optimization variables into multiple small blocks,
and optimize each block one by one following certain block selection rule, such as cyclic rule
\cite{tseng2001convergence}, Gauss-Southwell rule
\cite{tseng2009block}, etc.


In recent years, there are many applications of BCD-type algorithms in the areas of machine learning and data mining, such as matrix factorization \cite{zhao2015nonconvex,pmlrlu17,song17}, tensor decomposition, matrix completion/decomposition
\cite{xu2013block,jain2013low}, and training deep neural networks (DNNs) \cite{zima17}. Under relatively mild conditions, the convergence of BCD-type algorithms to first-order stationary solutions (SS1) have been broadly investigated for nonconvex and
non-differentiable optimization \cite{tseng2001convergence,Grippo00}. In particular, it is known that under mild conditions, these algorithms also achieve global sublinear rates \cite{meisam14nips}. 
However, despite its popularity and significant recent progress in understanding its behavior, it remains unclear whether BCD-type algorithms can converge to the set of second-order stationary solutions (SS2) with a provable global rate,  even for the simplest problem with two blocks of variables. 

\subsection{Motivation}

Algorithms that can escape from strict saddle points -- those stationary points that have negative eigenvalues --  have wide applications. Many recent works have analyzed the saddle points in machine learning problems \cite{kawaguchi2016deep}. Such as learning in shallow networks, the stationary points are either global minimum points or strict saddle points. In two-layer porcupine neural networks (PNNs), it has been shown that most local optima of PNN optimizations are also global optimizers \cite{soha17}.
Previous work in \cite{rong15fu} has shown that the saddle points in tensor decomposition are
indeed strict saddle points. Also, it has been shown that any saddle
points are strict in dictionary learning and phase retrieval problems theoretically and numerically in \cite{ju15qing,jusun17,nips2017wgsc,wagi17}. More recently, \cite{rong17} proposed a unified analysis of saddle points for a
board class of low rank matrix factorization problems, and they proved that
these saddle points are strict.

\subsection{Related Work}

Many recent works have been focused on the performance analysis and/or design of
algorithms with convergence guarantees to local minimum points/SS2 for nonconvex
optimization problems. These include the trust region method
\cite{conn2000trust}, cubic regularized Newton's method
\cite{nesterov2006cubic,carmon2016gradient}, and a mixed approach of the first-order and seconde-order methods \cite{sama17}, etc. However, these algorithms typically require second-order information, therefore they incur high computational complexity when problem dimension becomes large.

There has been a line of work on stochastic gradient descent algorithms, where properly scaled Gaussian noise is added to the
iterates of the gradient at each time [also known as stochastic gradient
Langevin dynamics, (SGLD)]. Some theoretical works have pointed out that SGLD
not only converges to the local minimum points asymptotically but also may
escape from local minima \cite{zhang2017hitting,raginsky2017non}.
Unfortunately, these algorithms require a large number of iterations with
$\mathcal{O}(1/\epsilon^{4})$ steps to achieve the optimal point. 
There are fruitful results that show some carefully designed algorithms can escape from strict saddle point efficiently, such as negative-curvature-originated-from noise (Neon) \cite{xu2017first},  Neon2 \cite{allen2017neon2}, Neon$^+$\cite{xuya2017} and gradient descent with one-step escaping (GOSE) \cite{yu2017saving}. The Neon-type of algorithms utilizes the stochastic first-order updates to find the negative curvature direction, and GOSE just needs one negative curvature descent step with calculation of eigenvectors when the iterates of the algorithm are near the saddle point for saving the computational burden.

On the other hand, there is also a line of work analyzing the deterministic GD type method. 
With random initializations, it
has been shown that GD only converges to SS2 for unconstrained smooth problems
\cite{jlee16jordan}. More recently, block coordinate descent, block mirror descent and proximal block coordinate descent have been proven to almost always converge to SS2 with random initializations \cite{jason17}, but there is no convergence rate reported. Unfortunately, a follow-up study indicated that GD requires exponential time to escape from saddle points for certain pathological problems \cite{sichi17}.
Adding some noise occasionally to the iterates of the algorithm is another way of finding the negative curvature. A perturbed version of GD has been proposed with convergence guarantees to SS2 \cite{jin2017jordan}, which shows a faster provable
convergence rate than the ordinary gradient descent algorithm with random initializations. Furthermore, the accelerated version of PGD (PAGD) is also proposed in \cite{jin2017accelerated}, which shows the fastest convergence rate among all Hessian free algorithms.


\begin{table*}[t]
\caption{Convergence rates of algorithms to SS2 with the first order information, where $p\ge4$, and $\widetilde{\mathcal{O}}$ hides factor  ploylog($d$).}
\label{tab.com}
\vskip 0.15in
\begin{center}
\begin{small}
\begin{sc}
\begin{tabular}{lll}
\toprule
Algorithm & Iterations  & $(\epsilon,\gamma)$-SS2 \\
\midrule
SGD \cite{rong15fu}    & $\mathcal{O}(d^p/\epsilon^{4}) $ & $(\epsilon,\epsilon^{1/4})$\\
SGLD \cite{zhang2017hitting} & $\mathcal{O}(d^p/\epsilon^{4})$ & $(\epsilon,\epsilon^{1/2})$\\
Neon+SGD \cite{xu2017first} & $\widetilde{\mathcal{O}}(1/\epsilon^{4})$ & $(\epsilon,\epsilon^{1/2})$ \\
Neon+Natasha \cite{xu2017first} & $\widetilde{\mathcal{O}}(1/\epsilon^{13/4})$ & $(\epsilon,\epsilon^{1/4})$ \\
Neon2+SGD \cite{allen2017neon2} & $\widetilde{\mathcal{O}}(1/\epsilon^4)$ & $(\epsilon,\epsilon^{1/2})$\\
$\textrm{Neon}^+$ \cite{xuya2017}  &   $\widetilde{\mathcal{O}}(1/\epsilon^{7/4})$ & $(\epsilon,\epsilon^{1/2})$ \\
PGD \cite{jin2017jordan}& $\widetilde{\mathcal{O}}(1/\epsilon^{2})$  & $(\epsilon,\epsilon^{1/2})$  \\
PAGD \cite{jin2017accelerated}& $\widetilde{\mathcal{O}}(1/\epsilon^{7/4})$ & $(\epsilon,\epsilon^{1/2})$\\
 PA-GD/PA-PP (This work)    & $\widetilde{\mathcal{O}}(1/\epsilon^{7/3})$  & $(\epsilon,\epsilon^{1/3})$\\
\bottomrule
\end{tabular}
\end{sc}
\end{small}
\end{center}
\vskip -0.1in
\end{table*}

\subsection{Scope of This Paper}

In this work, we consider a smooth unconstrained optimization problem, and develop a perturbed AGD algorithm (PA-GD) which converges (with high probability) to the set of SS2 with a global sublinear rate. Our work is inspired by the works \cite{jin2017jordan,rong15fu}, which developed novel perturbed GDs that escapes from strict
saddle points. Similarly as  in \cite{jin2017jordan}, we also divide the entire iterates of GD into three types of points: those whose gradients are
large, those that are local minimum, and those that are strict saddle
points. At a given point, when the size of the gradient is large
enough, we just implement the ordinary AGD. When the gradient norm is small,
which may be either strict saddle or local minimum, a perturbation will be
added on the iterates to help to escape from the saddle points.

From the above section, we know that  many works have been developed to make use of negative curvature information around the saddle points. Unfortunately, these
techniques cannot be directly applied to the BCD/AGD- type of algorithms. The {\it key challenge} here is that at each iteration only part of the variables are updated, therefore
we have access only to partial second order information at the points of interest. 
For example, consider a quadratic objective function shown in \figref{fig:1}. While fixing one block, the problem is strongly convex with respect to the other block, but the entire problem is nonconvex. Even if the iterates converge for each block to the minimum
points within the block, the stationary point could still be a saddle point
for the overall objective function. Therefore, the analysis of how AGD type of algorithms exploit the negative curvature is one of the main tasks in this paper.

To the best of our knowledge, there is no work on modifying AGD algorithms to escape from
strict saddle points with any convergence rate. The main contributions of this work are as follows.

\subsection{Contributions of This Work}

In this paper, we design and analyze a perturbed AGD algorithm for solving an
unconstrained nonconvex problem, namely perturbed AGD. Through the
perturbation of AGD, the algorithm is guaranteed to converge to a set of SS2 of a nonconvex problem with high probability. {By utilizing the matrix perturbation theory, convergence rate
of the proposed algorithm is also established, which shows that the
algorithm takes $\mathcal{O}(\text{polylog}(d)/\epsilon^{7/3})$ iterations to
achieve an ($\epsilon, \epsilon^{1/3}$)-SS2 with high probability}. Also, considering the fact that there is a strong relation between GD and proximal point algorithm, we also study a perturbed alternating proximal point (PA-PP) algorithm with some random perturbation. By leveraging the techniques proposed in this paper, we show that PA-PP, which may not need to calculate the gradient at each step, converges  as fast as PA-GD in the order of $\epsilon$ . The comparison of the algorithms which only use the first order information for escaping from strict saddle points is summarized as shown in \figtab{tab.com}.

The main contributions of the paper are highlighted below:
\vspace{-0.2cm}
\begin{enumerate}
\item  To the best of our knowledge, it is the first time that the convergence
analysis shows that some variants of AGD (using first-order information) can converge to SS2 for nonconvex optimization problems.

\vspace{-0.2cm}
\item The convergence rate of the perturbed AGD algorithm is analyzed, where the choice of the step size is only dependent on certain maximum Lipschitz constant over blocks rather than all variables. This is one of the major difference between GD and AGD.

\vspace{-0.2cm}
\item  By further extending the analysis in this paper, we also show that PA-PP can also escape from the strict points efficiently with the speed of $\mathcal{O}(\text{polylog}(d)/\epsilon^{7/3})$ .
\end{enumerate}
\vspace{-0.2cm}
\section{Preliminaries} \label{gen_inst}

\subsection{Notation}
{\bf Notation.} Bold upper case letters without subscripts (e.g., $\bX,\bY$) denote matrices
and bold lower case letters without subscripts (e.g., $\bx,\by$) represent
vectors. Notation $\bx_k$ denotes the $k$th block of vector
$\bx\in\mathbb{R}^{d\times 1}$. We use $\nabla_k f(\bx_{-k},\bx_k)$ to
denote the partial gradient with respect to its $k$th block variable while
the remaining one is fixed. Notation $\mathbb{B}_{\bx}(r)$ denotes a
$d$-dimensional ball centered at $\bx$ with radius $r$, and $\lambda_{\min}(\bX)$, $\lambda_{\max}(\bX)$ denote the smallest and largest eigenvalues of matrix $\bX$ respectively.

\subsection{Definitions}

The objective function has the following properties.
\begin{definition}
A differentiable function $f(\cdot)$ is L-smooth with gradient Lipschitz
constant $L$ (uniformly Lipschitz continuous), if
\begin{equation}
\notag
\|\nabla f(\bx)-\nabla f(\by)\|\le L\|\bx-\by\|,\quad\forall \bx,\by.
\end{equation}
The function is called block-wise smooth with gradient Lipschitz
constants $\{L_k\}$, if
\begin{equation}
\notag
\|\nabla_k f(\bx_{-k},\bx_k)-\nabla_k f(\bx_{-k},\bx'_k)\|\le L_k\|\bx_k-\bx'_k\|,\;\forall \bx,\bx'
\end{equation}
or with gradient Lipschitz
constants $\{\widetilde{L}_k\}$, if
\begin{equation}
\notag
\|\nabla_k f(\bx_{-k},\bx_k)-\nabla_k f(\bx'_{-k},\bx_k)\|\le \widetilde{L}_k\|\bx_{-k}-\bx'_{-k}\|,\;\forall\bx,\bx'.
\end{equation}
Further, let $L_{\max}\bydef\max\{ L_k, \widetilde{L}_k,\forall k\}\le L$.
\end{definition}
\begin{definition}
For a differentiable function $f(\cdot)$, if $\|\nabla f(\bx)\|=0$, then $\bx$
is a first-order stationary point. If $\|\nabla f(\bx)\|\le\epsilon$, then
$\bx$ is an $\epsilon$-first-order stationary point.
\end{definition}

\begin{definition}
For a differentiable function $f(\cdot)$, if $\bx$ is a SS1, and there exists $\epsilon>0$ so that for any $\by$ in the
$\epsilon$-neighborhood of $\bx$, we have $f(\bx)\le f(\by)$, then $\bx$ is a
local minimum. A saddle point $\bx$ is a SS1 that is
not a local minimum. If $\lambda_{\min}(\nabla^2f(\bx))<0$, $\bx$ is a strict
(non-degenerate) saddle point.
\end{definition}

\begin{definition}
A twice-differentiable function $f(\cdot)$ is $\rho$-Hessian Lipschitz if
\begin{equation}\label{eq.hlip}
\|\nabla^2f(\bx)-\nabla^2f(\by)\|\le\rho\|\bx-\by\|,\quad\forall \bx,\by.
\end{equation}
\end{definition}
\begin{definition}
For a $\rho$-Hessian Lipschitz function $f(\cdot)$, $\bx$ is a second-order
stationary point if $\|\nabla f(\bx)\|=0$ and
$\lambda_{\min}(\nabla^2f(\bx))\ge0$. If the following holds
\begin{equation}
\|\nabla f(\bx)\|\le\epsilon,\quad\textrm{and}\quad\lambda_{\min}(\nabla^2 f(\bx))\ge-\gamma\label{eq.csp}
\end{equation}
where $\epsilon,\gamma>0$, then $\bx$ is a $(\epsilon,\gamma)$-SS2.
\end{definition}
\begin{assumption}\label{as1}
Function $f(\cdot)$ is $L$-smooth, block-wise smooth with gradient
Lipschitz constants $\{L_k,\widetilde{L}_k\}, k=1,2$, and $\rho$-Hessian Lipschitz.
\end{assumption}

\begin{algorithm}[!t]
   \caption{Perturbed Alternating Gradient Descent (PA-GD) $(\bx^{(0)},L_{\max},L,\rho,\epsilon,\delta,\Delta f)$}
   \label{alg:p1}
\begin{algorithmic}
\State {\bfseries Input:} $\mathcal{P}_1=(1+\frac{L}{L_{\max}})$, $\mathcal{P}_2=(1+\frac{L\log(2d)}{L_{\max}})$, $\chi=6\max\{\log(\frac{\mathcal{P}^6_1\mathcal{P}^2_2dL^{5/3}_{\max}\Delta_f}{c^5\rho^{1/3}\epsilon^{7/3}\delta},4\}$,
$\eta=\frac{c}{L_{\max}}$,
$r=\frac{c^3}{\chi^{3}}\frac{\rho\epsilon}{L_{\max}\mathcal{P}^3_1\mathcal{P}_2}$,
$g_{\textsf{th}}=\frac{c^2\epsilon}{(\chi\mathcal{P}_1)^3\mathcal{P}_2}$,
$f_{\textsf{th}}=\frac{c^5\epsilon^{2}}{L_{\max}(\chi\mathcal{P}_1)^6\mathcal{P}^2_2}$,
$t_{\textsf{th}}=\frac{L_{\max}\chi\mathcal{P}_1}{c^2(L_{\max}\rho\epsilon)^{\frac{1}{3}}}$
   \For{$t=0,1,\ldots$}
   \If{$\sum^2_{k=1}\|\nabla_k f(\bh^{(t)}_{-k},\bx^{(t)}_k)\|^2\le
g^2_{\textsf{th}}$ and $t-t_{\textsf{p}}>t_{\textsf{th}}$}
   \State $\widetilde{\bx}^{(t)}\leftarrow\bx^{(t)}$ and $t_{\text{p}}\leftarrow t$
   \State $\bx^{(t)}=\widetilde{\bx}^{(t)}+\xi^{(t)}$, $\xi^{(t)}$
uniformly taken from $\mathbb{B}_0(r)$
   \EndIf
   \If{$t-t_{\textsf{p}}=t_{\textsf{th}}$ and
$f(\bx^{(t)})-f(\widetilde{\bx}^{(t_{\textsf{p}})})>-f_{\textsf{th}}$}
   \State {\bf return} $\widetilde{\bx}^{t_{\textsf{p}}}$
   \EndIf
   \For {$k=1,2$}
\State $\bx^{(t+1)}_k=\bx^{(t)}_{k}-\eta\nabla_k f(\bh^{(t)}_{-k},\bx^{(t)}_k)$
\EndFor
   \EndFor
\end{algorithmic}
\end{algorithm}

\section{Perturbed Alternating Gradient Descent}

\subsection{Algorithm Description}

AGD is a classical algorithm that optimizes the variables of an optimization problem in an alternating  manner \cite{bertsekas99}, meaning that when one block of variables is updated, the remaining block is fixed to be the same as its previous solution. Mathematically, the iterates of AGD are updated by the following rule
\begin{equation}\label{eq.gra}
\bx^{(t+1)}_k=\bx^{(t)}_{k}-\eta\nabla_k f(\bh^{(t)}_{-k},\bx^{(t)}_k),\quad k=1,2
\end{equation}
where superscript $(t)$ denotes the iteration counter;  $\bh^{(t)}_{-1}\bydef\bx^{(t)}_{2}$ and $\bh^{(t)}_{-2}\bydef\bx^{(t+1)}_1$; $\eta>0$ is the step size. AGD can be considered as a special case of block coordinate gradient descent \cite{nesterov2012efficiency,beck2013convergence}.

Our proposed algorithm is based on AGD, but modified in a way similar to  the recent work \citep{jin2017jordan}, which adds some noise in PGD. 
The details of the implementation of PA-GD are shown in
\algref{alg:p1}, where $c$ is a constant so that $\eta=c/L_{\max}$, $\Delta_f$ denotes the difference of the objective value at the initial point and global optimal solution, $\epsilon$ represents the predefined target error.


\begin{figure}[ht]
\begin{center}
\centerline{\includegraphics[width=0.46\columnwidth]{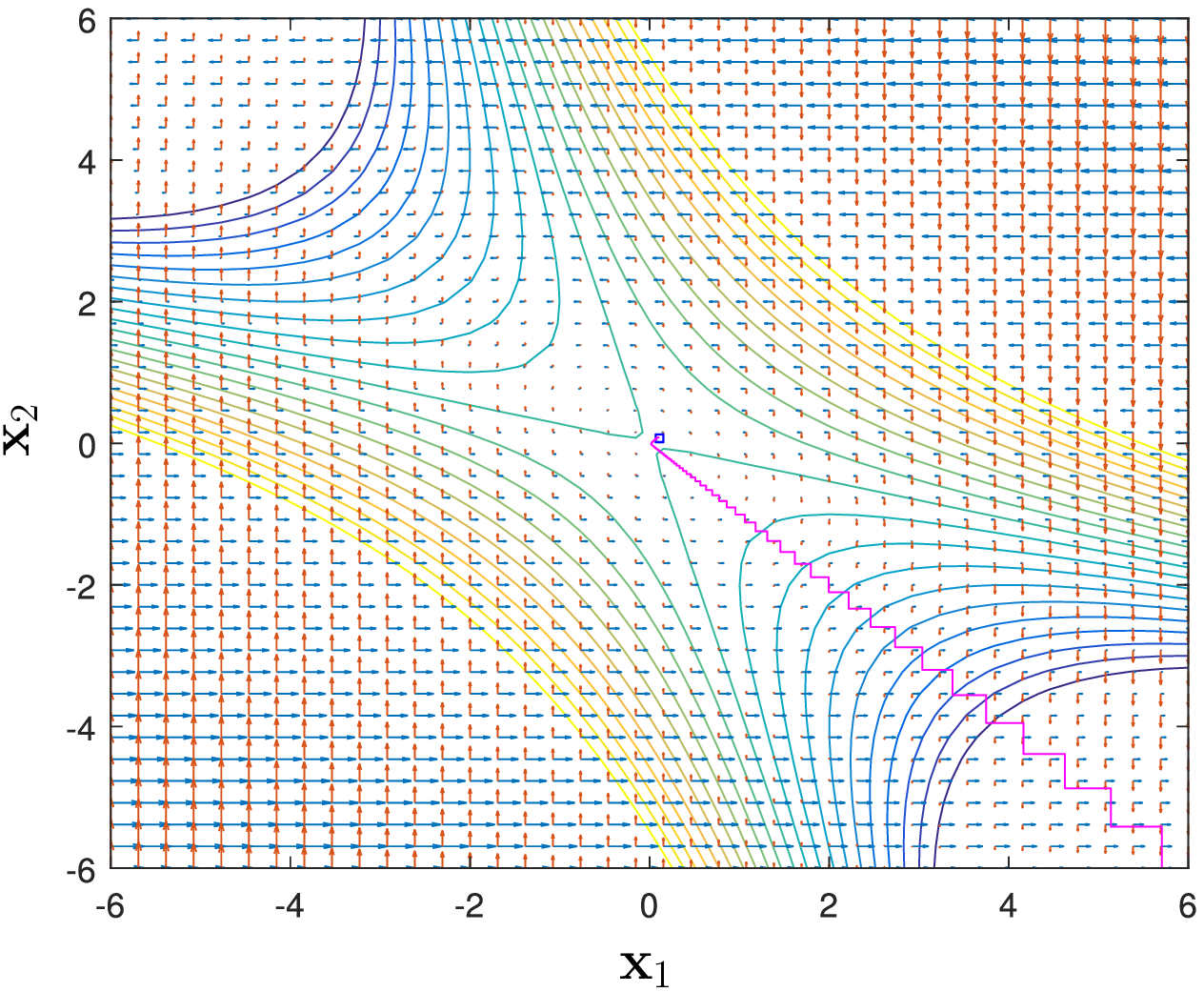}}
\caption{Contour of the objective values and the trajectory (pink color) of PA-GD started near strict saddle point $[0,0]$. The objective function is $f(\bx)=\bx^{\T}\bA\bx,\bx=[\bx_1;\bx_2]\in\mathbb{R}^{2\times 1}$ where $\bA\bydef[1\quad 2; 2 \quad 1]\in\mathbb{R}^{2\times 2}$, and the length of the arrows indicate the strength of $-\nabla f(\bx)$ projected onto directions $\bx_1,\bx_2$. }
\label{fig:1}
\end{center}
\end{figure}


In each update of variables, we implement one step of the block gradient descent, and then
proceed to the next block. Once the algorithm has sufficient decrease of the
objective value, it implies that the algorithm converges to some good
solution. Otherwise, some perturbation may be needed to help the iterates
escape from the saddle points. If after the perturbation the objective value
does not decrease sufficiently after a number of further iterations, the
algorithm terminates and returns the iterate before the last
perturbation.

To illustrate the practical behavior of the algorithm, we provide an example that shows the trajectory of AGD after a small perturbation at a stationary point. In \figref{fig:1}, it is clear that $\bx=[0;0]$ is a SS1 and also a strict saddle point since the eigenvalues of $\bA$ are $-1$ and $3$ respectively. When $\bx_1$ is fixed, function $f(\bx)$ is convex with respect to $\bx_2$ and vice versa, however, the objective function is nonconvex. It can be observed that PA-GD can escape from the strict saddle point efficiently. 

\subsection{Convergence Rate Analysis}

Despite the fact that PA-GD exploits a different way of updating variables,  we will show that it can still escape from strict saddle points with high probability
with suitable perturbation. The main theorem is presented as follows.

\begin{theorem}\label{th.mainth}
Under \asref{as1}, there exists a constant $c_{\max}$ such that: for any
$\delta\in(0,1]$, $\epsilon\le\frac{L^2_{\max}}{\rho}$, $\Delta_f\bydef f(\bh^{(0)}_{-1},\bx^{(0)}_1)-f^*$, and
constant $c\le c_{\max}$, with probability $1-\delta$, the iterates generated
by PA-GD converge to an $\epsilon$-SS2 $\bx$ satisfying
\begin{equation}\notag
\|\nabla f(\bx)\|\le\epsilon,\quad\textrm{and}\quad\lambda_{\min}(\nabla^2 f(\bx))\ge-(L_{\max}\rho\epsilon)^{1/3}
\end{equation}
in the following number of iterations:
\begin{align}\label{eq:rate:1}
\mathcal{O}\left(\frac{L^{5/3}_{\max}\mathcal{P}^7_1\mathcal{P}^2_2\Delta f}{\rho^{1/3}\epsilon^{7/3}}\log^{7}\left(\frac{\mathcal{P}^6_1\mathcal{P}^2_2dL^{5/3}_{\max}\Delta_f}{c^5\rho^{1/3}\epsilon^{7/3}\delta}\right)\right)
\end{align}
where $f^*$ denotes the global minimum value of the objective function, and $\mathcal{P}_1=(1+L/L_{\max})$ and $\mathcal{P}_2=(1+L\log(2d)/L_{\max})$.
\end{theorem}

%


\remark When $\eta=c_{\max}/L$ is used, the convergence rate of PA-GD is
\begin{equation}
\mathcal{O}\left(\frac{L^{5/3}_{\max}\log^2(2d)\Delta f}{\rho^{1/3}\epsilon^{7/3}}\log^{7}\left(\frac{\mathcal{P}^6_1\mathcal{P}^2_2dL^{5/3}_{\max}\Delta_f}{c^5\rho^{1/3}\epsilon^{7/3}\delta}\right)\right).
\end{equation}
It shows that if a smaller step size is used, the convergence rate of PA-GD is faster  (with smaller constants) {since the linear dependency of $\mathcal{P}^7_1$ and $\mathcal{P}^2_2$ in \eqref{eq:rate:1} both disappear. This property is consistent with the known result when BCD is used in convex optimization problems, i.e., when a smaller step size is used, the rate could become better; e.g., see {\citep[Theorem 2.1]{sun2015nips}}.}

\section{Perturbed  Alternating Proximal Point}

In many applications, AGD may not be efficient in the sense that the convergence rate of
gradient in each block may be very slow. For example, consider matrix factorization problem $\min_{\bX,\bY}\|\bZ-\bX\bY\|^2_F$
where $\bZ\in\mathbb{R}^{m\times d}$ is the given data, $d\gg m$, and
$\bX\in\mathbb{R}^{m\times r},\bY\in\mathbb{R}^{r\times d}$ are two block
variables. For this problem, the alternating least squares algorithm (which exactly minimizes each block) would be a faster algorithm compared with the AGD which only uses gradient steps.

In this section, we consider the classical proximal point algorithm \cite{parikh2014proximal} in which each block of variables is exactly minimized with respect to certain quadratic surrogate. To be specific, we can replace \eqref{eq.gra} in \algref{alg:p1} by
\begin{equation}\label{eq.bsum}
\bx^{(t+1)}_k=\arg\min_{\bx_k}f(\bh^{(t)}_{-k},\bx_k)+\frac{\nu}{2}\|\bx_k-\bx^{(t)}_k\|^2,\; k=1,2
\end{equation}
where $\nu>0$ is penalty parameter. The iteration can be explicitly written as 
\begin{equation}
\bx^{(t+1)}_k=\bx^{(t)}_k-\frac{1}{\nu}\nabla_k f(\bh^{(t)}_{-k},\bx^{(t+1)}_k),\quad\ k=1,2,\label{eq.upp}
\end{equation}
which has the similar form as the PA-GD algorithm, but with the step size
being $\eta\bydef1/\nu$, and with gradient evaluated at the new iterate. The resulting algorithm, detailed in the table above, is referred to as the {perturbed alternating proximal point (PA-PP)}. It is worth noting that when the subproblem is convex,
such as $\min_{\bX,\bY}\|\bZ-\bX\bY\|^2_F$, $\nu$ only needs to be a small number to make the
corresponding subproblem strongly convex. This property is useful in practice.



Next, we can also give the convergence rate of PA-PP.
\begin{algorithm}[tb]
\caption{Perturbed Alternating Proximal Point (PA-PP) $(\bx^{(0)},L_{\max},L,\rho,\epsilon,\delta,\Delta f)$}
   \label{alg:p2}
\begin{algorithmic}
\State {\bfseries Input:} $\mathcal{P}=(1+\frac{L\log(2d)}{L_{\max}})$, $\chi=6\max\{\log(\frac{\mathcal{P}^2dL^{5/3}_{\max}\Delta_f}{c^5\rho^{1/3}\epsilon^{7/3}\delta},4\}$,
$\nu=\frac{L_{\max}}{c}$,
$r=\frac{c^3}{\chi^{3}}\frac{\rho\epsilon}{L_{\max}\mathcal{P}}$,
$g_{\textsf{th}}=\frac{c^2\epsilon}{\chi^3\mathcal{P}}$,
$f_{\textsf{th}}=\frac{c^5\epsilon^{2}}{L_{\max}\chi^6\mathcal{P}^2}$,
$t_{\textsf{th}}=\frac{L_{\max}\chi}{c^2(L_{\max}\rho\epsilon)^{\frac{1}{3}}}$
   \For{$t=0,1,\ldots$}
   \If{$\|\bx^{(t+1)}-\bx^{(t)}\|\le
g_{\textsf{th}}/\nu$ and $t-t_{\textsf{p}}>t_{\textsf{th}}$}
   \State $\widetilde{\bx}^{(t)}\leftarrow\bx^{(t)}$ and $t_{\text{p}}\leftarrow t$
   \State $\bx^{(t)}=\widetilde{\bx}^{(t)}+\xi^{(t)}$, $\xi^{(t)}$
uniformly taken from $\mathbb{B}_0(r)$
   \EndIf
   \If{$t-t_{\textsf{p}}=t_{\textsf{th}}$ and
$f(\bx^{(t)})-f(\widetilde{\bx}^{(t_{\textsf{p}})})>-f_{\textsf{th}}$}
   \State {\bf return} $\widetilde{\bx}^{t_{\textsf{p}}}$
   \EndIf
   \For {$k=1,2$}
\State $\bx^{(t+1)}_k=\arg\min_{\bx_k}f(\bh^{(t)}_{-k},\bx_k)+\frac{\nu}{2}\|\bx_k-\bx^{(t)}_k\|^2$
\EndFor
   \EndFor
\end{algorithmic}
\end{algorithm}

\begin{corollary}\label{co.rate}
Under \asref{as1}, there exists a constant $c_{\max}$ such that: for any
$\delta\in(0,1]$, $\epsilon\le\frac{L^2_{\max}}{\rho}$, $\Delta_f\bydef f(\bh^{(0)}_{-1},\bx^{(0)}_1)-f^*$, and
constant $c\le c_{\max}$, with probability $1-\delta$, the iterates generated
by PA-PP converges to an $\epsilon$-SS2 $\bx$ satisfying
\begin{equation}\notag
\|\nabla f(\bx)\|\le\epsilon,\quad\textrm{and}\quad\lambda_{\min}(\nabla^2 f(\bx))\ge-(L_{\max}\rho\epsilon)^{1/3}
\end{equation}
in the following number of iterations:
\begin{equation}\notag
\mathcal{O}\left(\frac{L^{5/3}_{\max}\mathcal{P}^2\Delta f}{\rho^{1/3}\epsilon^{7/3}}\log^{7}\left(\frac{\mathcal{P}^2dL^{5/3}_{\max}\Delta_f}{c^5\rho^{1/3}\epsilon^{7/3}\delta}\right)\right)
\end{equation}
where $f^*$ denotes the global minimum value of the objective function, and $\mathcal{P}=(1+L\log(2d)/L_{\max})$.
\end{corollary}
Comparing with \thref{th.mainth}, we can find that term $\mathcal{P}^7_1, \mathcal{P}_1>2$ is removed so the convergence rate of PA-PP is slightly faster than PA-GD.

\section{Convergence Analysis}\label{sec.converge}

In this section, we will present the main proof steps
of convergence analysis of PA-GD.
\subsection{The Main Difficulty of the Proof}

\paragraph{Gradient Descent:}
GD searches the descent direction of the objective function in the entire space $\mathbb{R}^d$. Without loss of generality, we assume $\bx^{(0)}=0$.  According to the mean value theorem, the GD update can be expressed as
\begin{equation}
\bx^{(t+1)}=\bx^{(t)}-\eta\nabla f(\bx^{(t)})=\bx^{(t)}-\eta\nabla f(0)-\eta\left(\int^1_0\nabla^2f(\theta\bx^{(t)})d\theta\right)\bx^{(t)}.\label{eq.regd}
\end{equation}
It can be observed that the update rule of GD contains the information of the Hessian matrix at point $\bx^{(t)}$, i.e., $\nabla^2f(\theta\bx^{(t)})$. To be more specific, letting $\bH\triangleq\nabla^2f(\bx^*)$ where $\bx^*$ denotes an $\epsilon$-SS2 satisfying \eqref{eq.csp}, we can rewrite \eqref{eq.regd} as
\begin{equation}
\bx^{(t+1)}=(\bI-\eta\bH)\bx^{(t)}-\eta\Delta^{(t)}\bx^{(t)}-\eta\nabla f(0)
\end{equation}
where $\Delta^{(t)}\bydef\int^1_0(\nabla^2f(\theta\bx^{(t)})-\bH)d\theta$.

Based on the $\rho$-Hessian Lipschitz property, we can quantify $\|\Delta^{(t)}\|$ that is upper bounded by the difference of iterates. By exploiting the negative curvature of the Hessian matrix at saddle point $\bx^*$, we can project the iterate onto the direction $\bd$ where the eigenvalue of $\bI-\eta\bH$ is greater than 1. This leads to the fact that the norm of the iterates projected along direction $\bd$ will be increasing exponentially as the algorithm proceeds around point $\bx^*$, implying the sequence generated by GD is escaping from the saddle point. The details of characterizing the convergence rate have been analyzed previously in \cite{jin2017jordan}.

\paragraph{Alternating Gradient Descent:}
However, the AGD algorithm only updates partial variables of vector $\bx$, which belong to a subspace of the feasible set. Similarly, from the mean value theorem we can express the AGD rule of updating variables with assuming $\bx^{(0)}=0$ as follows:
\begin{align}
\notag
\bx^{(t+1)}&=\bx^{(t)}-\eta\left[\begin{array}{c}\nabla_1f(\bx^{(t)}_1,\bx^{(t)}_2) \\ \nabla_2f(\bx^{(t+1)}_1,\bx^{(t)}_2)\end{array}\right]
\\
&=\bx^{(t)}-\eta\nabla f(0)-\eta\int^1_0\bH^{(t)}_ld\theta\bx^{(t+1)}-\eta\int^1_0\bH^{(t)}_ud\theta\bx^{(t)}\label{eq.secagd}
\end{align}
where
\begin{equation}\notag
\bH^{(t)}_l\bydef\left[\begin{array}{cc} \boldsymbol{0} & \boldsymbol{0} \\
\\
\nabla^2_{21} f(\theta\bx^{(t+1)}_1, \theta\bx^{(t)}_2) & \boldsymbol{0}  \\
\end{array}\right]
\quad\textrm{and}\quad
\bH^{(t)}_u\bydef\left[\begin{array}{cc}\nabla^2_{11} f(\theta\bx^{(t)}_1, \theta\bx^{(t)}_2) & \nabla^2_{12} f(\theta\bx^{(t)}_1, \theta\bx^{(t)}_2) \\
\\
\boldsymbol{0} & \nabla^2_{22} f(\theta\bx^{(t+1)}_1, \theta\bx^{(t)}_2)  \\
\end{array}\right].
\end{equation}
From the above expression, it can be seen clearly that the update rule of AGD does not include a full Hessian matrix at any point but only partial ones. Furthermore, the right hand side of \eqref{eq.secagd} not only contains the second order information of the previous point, i.e., $[\bx^{(t)}_1,\bx^{(t)}_2]$ but also the one of the most recently updated point, i.e., $[\bx^{(t+1)}_1,\bx^{(t)}_2]$. These represent the main challenges in understanding the behavior of the sequence generated by the AGD algorithm.
\subsection{The Main Idea of the Proof}

Although the second order information is divided into two parts, we can still characterize the recursion of the iterates around strict saddle points. We can also split $\bH$ as two parts, which are
\begin{align}
\bH_u=\left[\begin{array}{cc}\nabla^2_{11} f(\bx^*) &  \nabla^2_{12} f(\bx^*)
\\
\boldsymbol{0} &   \nabla^2_{22} f(\bx^*)\end{array}\right],\quad
\bH_l=\left[\begin{array}{cc}\boldsymbol{0} & \boldsymbol{0}
\\
\nabla^2_{21} f(\bx^*) & \boldsymbol{0} \end{array}\right],\label{eq.defhl}
\end{align}
and obviously we have $\bH=\bH_l+\bH_u$.

Then, recursion \eqref{eq.secagd} can be written as
\begin{equation}\label{eq.reofx}
\bx^{(t+1)}+\eta\bH_l\bx^{(t+1)}=\bx^{(t)}- \eta\bH_u\bx^{(t)}-\eta\Delta^{(t)}_u \bx^{(t)} -\eta\Delta^{(t)}_l \bx^{(t+1)}
\end{equation}
where $\Delta^{(t)}_u\bydef\int^1_0(\bH^{(t)}_u(\theta)-\bH_u)d\theta$, $\Delta^{(t)}_l\bydef\int^1_0(\bH^{(t)}_l(\theta)-\bH_l)d\theta$. However, it is still unclear from \eqref{eq.reofx} how the iteration evolves around the strict saddle point.

To highlight ideas, let us define
\begin{equation}
\bM\bydef\bI+\eta\bH_l,\quad\bT\bydef\bI-\eta\bH_u.\label{eq.defofmt}
\end{equation}
It can be observed that $\bM$ is a lower triangular matrix where the diagonal entries are all 1s; therefore it is invertible. After taking the inverse of matrix $\bM$ on both sides of \eqref{eq.reofx}, we can obtain
\begin{equation}\notag
\bx^{(t+1)}=\bM^{-1}\bT\bx^{(t)}- \eta\bM^{-1}\Delta^{(t)}_u \bx^{(t)} -\eta\bM^{-1}\Delta^{(t)}_l \bx^{(t+1)}.
\end{equation}

Our goal of analyzing the recursion of $\bx^{(t)}$ becomes to find the maximum eigenvalue of $\bM^{-1}\bT$. With the help of the matrix perturbation theory, we can quantify the difference between the eigenvalues of matrix $\bH$ that contains the negative curvature and matrix $\bM^{-1}\bT$ that we are interested in analyzing. To be more precise, we give the following lemma.
\begin{lemma}\label{le.eiginver}
Under \asref{as1}, let $\bH\bydef\nabla^2f(\bx)$ denote the Hessian matrix at an $\epsilon$-SS2 $\bx$ where $\lambda_{\min}(\bH)\le-\gamma$ and $\gamma>0$. We have
\begin{align}
\lambda_{\max}(\bM^{-1}\bT)>1+\frac{\eta\gamma}{1+L/L_{\max}}
\end{align}
where $\bM,\bT$ are defined in \eqref{eq.defhl} and \eqref{eq.defofmt}.
\end{lemma}

\leref{le.eiginver} illustrates that there exits a subspace spanned by the eigenvector of $\bM^{-1}\bT$ whose eigenvalue is greater than 1, indicating that the sequence generated by AGD can still potentially escape from the strict saddle point by leveraging such negative curvature information. Next, we can give a sketch of the proof of \thref{th.mainth}.

\subsection{The Sketch of the Proof}

The structure of the proof for quantifying the sufficient decrease of the objective function after the perturbation is borrowed from the proof of PGD \cite{jin2017jordan}, but PA-GD updates the variables block by block, so we have to provide the new proofs to show that PA-GD can still escape from saddle points with the perturbation technique.

First, if the size of the gradient is large enough, \algref{alg:p1} just
implements the ordinary AGD. We give the descent lemma of AGD as follows.
\begin{lemma}\label{le.descent}
Under \asref{as1}, for the AGD algorithm with step size $\eta<1/L_{\max}$, we have
\begin{equation}
f(\bx^{(t+1)})\le f(\bx^{(t)})-\sum^2_{k=1}\frac{\eta}{2}\|\nabla_k f(\bh^{(t)}_{-k},\bx^{(t)}_k)\|^2.\notag
\end{equation}
\end{lemma}

Second, if the iterates are near a strict saddle point, we can show that the AGD algorithm
after a perturbation can give a sufficient decrease with high probability in terms of the
objective value. To be more precise, the statement is given as follows.

\begin{lemma}\label{le.escape}
Under \asref{as1}, there exists a absolute constant $c_{\max}$. Let $c\le c_{\max}$,
$\chi\ge 1$, and $\eta$, $r$, $g_{\textsf{th}}$, $t_{\textsf{th}}$ calculated
as \algref{alg:p1} describes. Let $\wbx^{(t)}$ be a strict saddle point, which
satisfies
\begin{equation}
\|\nabla f(\wbx^{(t)})\|^2\le 4\sum^2_{k=1}\|\nabla_k f(\tbh^{(t)}_{-k},\wbx^{(t)}_k)\|^2\le 4g^2_{\textsf{th}}\label{eq.conds}
\end{equation}
and $\lambda_{\min}(\nabla^2 f(\wbx^{(t)}))\le -\gamma$, where $\tbh^{(t)}_{-1}\bydef\wbx^{(t)}_2$ and $\tbh^{(t)}_{-2}\bydef\bx^{(t+1)}_1$.

Let $\bx^{(t)}=\wbx^{(t)}+\xi^{(t)}$ where $\xi^{(t)}$ is generated randomly
which follows the uniform distribution over $\mathbb{B}_0(r)$, and let
$\bx^{(t+t_{\textsf{th}})}$ be the iterates of PA-GD. With at least probability
$1-\frac{dL_{\max}}{(L_{\max}\rho\epsilon)^{1/3}}e^{-\chi}$, we have $f(\bx^{(t+t_{\textsf{th}})})-f(\wbx^{(t)})\le-f_{\textsf{th}}$.
\end{lemma}

We remark that Lemma \ref{le.descent} is well-known and Lemma \ref{le.escape} is the core technique. In the following, we outline the main
idea used in proving the latter. The formal statements of these steps are shown in the appendix; see \leref{le.layer31}--\leref{le.layer21} therein.

We emphasize that the main contributions of this paper lies in the analysis of the first two steps, {where the special update rule of PA-GD is analyzed so that the negative curvature of $\bH$ around the saddle points can be utilized.}
\paragraph{Step 1} (\leref{le.layer31}) Consider a generic sequence $\bu^{(t)}$ generated by {PA-GD}. As long as the initial point of $\bu^{(t)}$ is close to saddle point $\wbx^{(t)}$, the distance between $\bu^{(t)}$ and $\wbx^{(t)}$ can be upper bounded by using the $\rho$-Hessian  Lipschitz continuity property.
\paragraph{Step 2} (\leref{le.layer32})  Leveraging the negative curvature around the strict saddle point, we know that there exits a direction, i.e., $\vec{\be}$, which is spanned by the eigenvector of $\bM^{-1}\bT$ whose corresponding eigenvalue is largest (greater than 1). Consider two sequences generated by {PA-GD}, $\bu^{(t)},\bw^{(t)}$ initialized around the saddle point. When the initial points of these two iterates are separated apart away from each other along direction $\vec{\be}$ with a small distance, meaning that $\bw^{(0)}=\bu^{(0)}+\upsilon r \vec{\be},\;\upsilon\in[\delta/(2\sqrt{d}),1]$ where $r$ denotes the radius of the perturbation ball defined in \algref{alg:p1}, we can show that if iterate $\bu^{(t)}$ is still near the saddle point after $T$ steps, the other sequence $\bw^{(t)}$ will give a sufficient decrease of the objective value with less than $T$ steps, implying that iterates $\bw^{(t)}$ can escape from the saddle point with less than $T$ steps.
\paragraph{Step 3}(\leref{le.layer21}) Consider $\bu^{(0)},\bw^{(0)}$ as the points after the perturbation from the saddle point. We can quantify the probability that the AGD sequence will give a sufficient decrease of the objective value within $T$ iterations after the perturbation \citep[Lemma 14,15]{jin2017jordan}.  

\subsection{Extension to PA-PP}

By leveraging the convergence analysis of PA-GD and relation between PA-GD and PA-PP shown in \eqref{eq.upp}, we can also write the recursion of the PA-PP iteration as
\begin{equation}\label{eq.breofx}
\bx^{(t+1)}+\eta\bH'_l\bx^{(t+1)}=\bx^{(t)}- \eta\bH'_u\bv^{(t)}
-\eta\Delta'^{(t)}_u \bx^{(t)} -\eta\Delta'^{(t)}_l \bx^{(t+1)}
\end{equation}
where $\eta=1/\nu$, $\Delta'^{(t)}_u\bydef\int^1_0(\bH'^{(t)}_u(\theta)-\bH'_u)d\theta$, $\Delta'^{(t)}_l\bydef\int^1_0(\bH'^{(t)}_l(\theta)-\bH'_l)d\theta$,
\begin{equation}
\bH'_u=\left[\begin{array}{cc} \boldsymbol{0} &  \nabla^2_{12} f(\wbx^{(t)})
\\
\boldsymbol{0} &   \boldsymbol{0}\end{array}\right],
\quad
\bH'_l=\left[\begin{array}{cc} \nabla^2_{11} f(\wbx^{(t)}) & \boldsymbol{0}
\\
\nabla^2_{21} f(\wbx^{(t)}) & \nabla^2_{22} f(\wbx^{(t)}) \end{array}\right],\label{eq.bdefhl}
\end{equation}
and
\begin{equation}\notag
\bH'^{(t)}_l\bydef\left[\begin{array}{cc}\nabla^2_{11} f(\theta\bx^{(t+1)}_1, \theta\bx^{(t)}_2) & \boldsymbol{0} \\
\\
\nabla^2_{21} f(\theta\bx^{(t+1)}_1, \theta\bx^{(t+1)}_2) & \nabla^2_{22} f(\theta\bx^{(t+1)}_1, \theta\bx^{(t+1)}_2)  \\
\end{array}\right],
\quad
\bH'^{(t)}_u\bydef\left[\begin{array}{cc}\boldsymbol{0} & \nabla^2_{12} f(\theta\bx^{(t+1)}_1, \theta\bx^{(t)}_2) \\
\\
\boldsymbol{0} & \boldsymbol{0}  \\
\end{array}\right].
\end{equation}

Let
\begin{equation}
\bM'\bydef\bI+\eta\bH'_l\quad\bT'\bydef\bI-\eta\bH'_u.\label{eq.bdefofmt}
\end{equation}
We know that $\bT'$ is an upper triangular matrix where the diagonal entries are all 1s, so it is invertible. Different from the case of PA-GD, we take the inverse of matrix $\bT'$ on both sides of \eqref{eq.breofx} and obtain
\begin{equation}\notag
\bT'^{-1}\bM'\bx^{(t+1)}=\bx^{(t)}- \eta\bT'^{-1}\Delta'^{(t)}_u \bx^{(t)}-\eta\bT'^{-1}\Delta'^{(t)}_l \bx^{(t+1)}.
\end{equation}
Then, we can give the following result that characterizes the recursion of $\bx^{(t)}$ generated by PA-PP.
\begin{corollary}\label{co.eiginver}
Under \asref{as1}, let $\bH\bydef\nabla^2f(\bx)$ denote the Hessian matrix at an $\epsilon$-SS2 $\bx$ where $\lambda_{\min}(\bH)\le-\gamma$ and $\gamma>0$. Let ${\lambda}^{+}_{\min}(\cdot)$ denote the {\it minimum positive eigenvalue of a matrix}. Then we have
\begin{align}
{\lambda}^{+}_{\min}(\bT'^{-1}\bM') \le 1-\eta\gamma/2
\end{align}
where $\bM',\bT'$ are defined in \eqref{eq.bdefhl} and \eqref{eq.bdefofmt}; $\eta\le1/L_{\max}$ and $\gamma\le L_{\max}$. 
\end{corollary}
We remark that \coref{co.eiginver} is useful since it can be leveraged to show that the norm of the iterates around saddle points can increase exponentially. Then, we can apply the similar analysis steps as the case of proving the convergence rate of PA-GD and obtain the results shown in \coref{co.rate}.
\section{Connection with Existing Works}

\remark In Theorem \ref{th.mainth} we characterized the convergence rate to an $(\epsilon,{\epsilon}^{1/3})$-SS2. We can also translate this bound to the one for achieving an $(\epsilon,\sqrt{\epsilon})$-SS2, and in this case PA-GD needs $\widetilde{\mathcal{O}}(1/\epsilon^{3.5})$ iterations. Compared with the existing recent works \cite{jin2017jordan}, the convergence rate of PA-GD/PA-PP is slower than GD. The main reason is the fact that different from GD-type algorithms, PA-GD and PA-PP cannot fully utilize the Hessian information because they never see a full iteration. 
Similar situation happens for SGD-type of algorithms which also cannot get the exact negative curvature around strict saddle points.

From  \figtab{tab.com}, it can be seen that the convergence rate of PA-GD/PA-PP is still faster than SGD \cite{rong15fu}, SGLD \cite{zhang2017hitting}, Neon+SGD \cite{xu2017first}, and Neon2+SGD \cite{allen2017neon2} to achieve an $(\epsilon,\sqrt{\epsilon})$-SS2, but slower than the rest. We emphasize that PA-GD and PA-PP represent the first BCD-type algorithms with the convergence rate guarantee to escape from the strict saddle points efficiently. At this point, it is unclear whether our rate is the best that is achievable, and the question of whether the resulting rate can be improved will be left to future work.


\section{Numerical Results}

%

\begin{figure*}[htb]
  \centering
  \subfigure[Objective function in 2D.]{
    \label{fig:3}
    \includegraphics[width=.43\columnwidth]{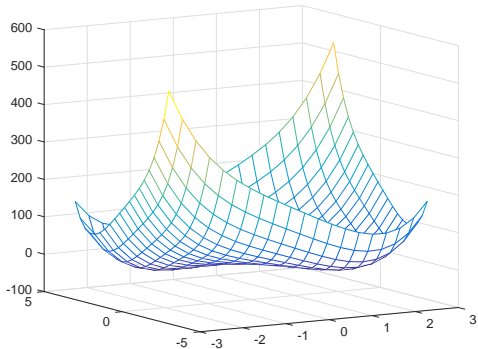}}
  \hspace{0.6in}
  \subfigure[{Objective value versus the number of iterations}]{
    \label{fig:2}
    \includegraphics[width=.43\columnwidth]{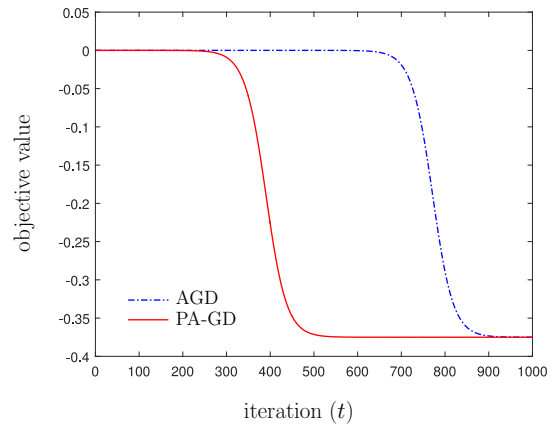}}
  \caption{Convergence comparison between AGD and PA-GD, where $\epsilon=10^{-4}$, $g_{\textsf{th}}=\epsilon/10$, $\eta=0.02$, $t_{\textsf{th}}=10/\epsilon^{1/3}$, $r=\epsilon/10$.}
  \label{fig:largesparsedata}
\end{figure*}

In this section, we present a simple example that shows the convergence behavior of PA-GD. Consider a nonconvex objective function, i.e.,
\begin{equation}
f(\bx)\bydef\bx^{\T}\bA\bx+\frac{1}{4}\|\bx\|^4_4.\label{eq.siobj}
\end{equation}
First, we have the following properties of function $f(\bx)$ such that $f(\bx)$ satisfies the assumptions of the analysis.
\begin{lemma}\label{le.simup}
For any $\tau\ge\lambda_{\max}(\bA)$ and $\bx\in\{\bx|\|\bx\|^2\le\tau\}$, $f(\bx)$ defined in \eqref{eq.siobj} is $5\tau$-smooth and $6\sqrt{\tau}$-Hessian Lipschitz.
\end{lemma}

Here, we can easily show the shape of objective function \eqref{eq.siobj} in the two dimensional (2D) case in \figref{fig:3}, where $\bA=[1\:2;2\;1]\in\mathbb{R}^{2\times 2}$. It can be observed clearly that there exits a strict saddle point at $[0,0]$ and two other local optimal points. We randomly initialize the algorithms around strict saddle point $[0,0]$. The convergence comparison between AGD and PA-GD is shown in \figref{fig:2}. It can be observed that PA-GD converges faster than AGD to a local optimal point.

\section{Conclusion}In this paper,
the perturbed variants of AGD and alternating proximal point (APP) algorithms are proposed, with the objective of finding the second order stationary solutions of nonconvex smooth problems. Leveraging the recently developed idea of random perturbation for the first-order methods, the
proposed algorithms add suitable perturbation to the AGD or APP iterates.
The main contribution of this work is a new analysis that takes into consideration the block structure of the updates for the perturbed AGD and APP algorithms. By exploiting the negative curvature, it is established that with high probability the algorithms can converge to an $(\epsilon,\epsilon^{1/3})$-SS2 with $\mathcal{O}(\text{polylog}(d)/\epsilon^{7/3})$ iterations.

\section{Acknowledgment}

The authors would like to thank Chi Jin for discussion on the perturbed gradient descent algorithm.

\appendix
\clearpage
\newpage


\bibliography{refs}
\bibliographystyle{plainnat}

\newpage
\onecolumn
\begin{center}
\bf{\Large Appendix}
\end{center}
\section{Preliminary}

We provide the proofs of some preliminary lemmas (\leref{le.layer00}--\leref{le.layer001grad}) used in the proof of \secref{sec.pofpagd}.

First, \leref{le.layer00} and \leref{le.layer001bl} give the property that quantify the size of
the difference of the second-order information of the objective values between
two points.
\begin{lemma}\label{le.layer00}
If function $f(\cdot)$ is $\rho$-Hessian Lipschitz, we have
\begin{equation}\label{eq.lemmah}
\left\|\int^1_0\nabla^2 f(\theta\bx)d\theta-\nabla^2
f(\by)\right\|\le\rho\left(\|\bx\|+\|\by\|\right),\quad\forall \bx,\by.
\end{equation}
\end{lemma}

\begin{lemma}\label{le.layer001bl}
Under \asref{as1}, we have block-wise Lipschitz continuity as follows:
\begin{equation}
\left\|\left[\begin{array}{cc}\nabla^2_{11}f(\bx) & \nabla^2_{12}f(\bx)
\\
\boldsymbol{0} & \nabla^2_{22}f(\by)\end{array}\right]-\left[\begin{array}{cc}\nabla^2_{11}f(\bz) & \nabla^2_{12}f(\bz)\\
\boldsymbol{0} & \nabla^2_{22}f(\bz)\end{array}\right]\right\|\le\rho\left(\|\bx-\bz\|+\|\by-\bz\|\right),\forall\bx,\by,\bz,\label{eq.deltatbd}
\end{equation}
and
\begin{equation}
\left\|\left[\begin{array}{cc} \boldsymbol{0} & \boldsymbol{0}
\\
\nabla^2_{21}f(\bx) & \boldsymbol{0}\end{array}\right]-\left[\begin{array}{cc} \boldsymbol{0} & \boldsymbol{0}\\
\nabla^2_{21}f(\by) & \boldsymbol{0}\end{array}\right]\right\|
\le\rho\|\bx-\by\|,\forall\bx,\by\label{eq.deltatbd2}.
\end{equation}
\end{lemma}
%


Then, we illustrate that the size of the partial gradient with one round update by the AGD algorithm has the following relation with the full size of the gradient.
\begin{lemma}\label{le.layer001grad}
If function $f(\cdot)$ is $L$-smooth with Lipschitz constant, then we have
\begin{equation}
\|\nabla f(\bx^{(t)})\|^2\le 4\sum^2_{k=1}\|\nabla_k f(\bh^{(t)}_{-k},\bx^{(t)}_{k})\|^2
\end{equation}
where sequence $\bx^{(t)}_k, k=1,2$ is generated by the AGD algorithm.
\end{lemma}

\subsection{Proof of \leref{le.layer00}}
\begin{proof}
If function $f(\cdot)$ is $\rho$-Hessian Lipschitz, then we have
\begin{align}
\notag
&\left\|\int^1_0(\nabla^2f(\theta\bx)-\nabla^2f(\by))d\theta\right\|
\le\int^1_0\left\|\nabla^2f(\theta\bx)-\nabla^2f(\by)\right\|d\theta
\\\notag
\mathop{\le}\limits^{(a)}&\rho\int^1_0\left\|\theta\bx-\by\right\|d\theta
\mathop{\le}\limits^{(b)}\rho\int^1_0\theta\|\bx\|d\theta+\rho\|\by\|
\le\rho\left(\|\bx\|+\|\by\|\right)
\end{align}
where $(a)$ is true because of Hessian Lipschitz, in $(b)$ we used
the triangle inequality.
\end{proof}

%
%
%

\subsection{Proof of \leref{le.layer001bl}}

There proof involves two parts:
\paragraph{Upper Triangular Matrix:}
Consider three different vectors $\bx$, $\by$ and $\bz$. We can have
\begin{align}
\notag
&\left\|\left[\begin{array}{cc}\nabla^2_{11}f(\bx) & \nabla^2_{12}f(\bx)
\\
0 & \nabla^2_{22}f(\by)\end{array}\right]-\left[\begin{array}{cc}\nabla^2_{11}f(\bz) & \nabla^2_{12}f(\bz)\\
0 & \nabla^2_{22}f(\bz)\end{array}\right]\right\|
\\\notag
\le&\left\|\bI_1\left(\left[\begin{array}{cc}\nabla^2_{11}f(\bx) & \nabla^2_{12}f(\bx)
\\
\nabla^2_{21}f(\bx) & \nabla^2_{22}f(\bx)\end{array}\right]-\left[\begin{array}{cc}\nabla^2_{11}f(\bz) & \nabla^2_{12}f(\bz)\\
\nabla^2_{21}f(\bz) & \nabla^2_{22}f(\bz)\end{array}\right]\right)\right\|
\\\notag
&+\left\|\bI_2\left(\left[\begin{array}{cc}\nabla^2_{11}f(\by) & \nabla^2_{12}f(\by)
\\
\nabla^2_{21}f(\by) & \nabla^2_{22}f(\by)\end{array}\right]-\left[\begin{array}{cc}\nabla^2_{11}f(\bz) & \nabla^2_{12}f(\bz)\\
\nabla^2_{21}f(\bz) & \nabla^2_{22}f(\bz)\end{array}\right]\right)\bI_2\right\|
\\\notag
\mathop{\le}\limits^{(a)}&\left\|\left[\begin{array}{cc}\nabla^2_{11}f(\bx) & \nabla^2_{12}f(\bx)
\\
\nabla^2_{21}f(\bx) & \nabla^2_{22}f(\bx)\end{array}\right]-\left[\begin{array}{cc}\nabla^2_{11}f(\bz) & \nabla^2_{12}f(\bz)\\
\nabla^2_{21}f(\bz) & \nabla^2_{22}f(\bz)\end{array}\right]\right\|
+\left\|\left[\begin{array}{cc}\nabla^2_{11}f(\by) & \nabla^2_{12}f(\by)
\\
\nabla^2_{21}f(\by) & \nabla^2_{22}f(\by)\end{array}\right]-\left[\begin{array}{cc}\nabla^2_{11}f(\bz) & \nabla^2_{12}f(\bz)\\
\nabla^2_{21}f(\bz) & \nabla^2_{22}f(\bz)\end{array}\right]\right\|
\\\notag
\le&\rho\left(\|\bx-\bz\|+\|\by-\bz\|\right)
\end{align}
where in $(a)$ we used
\begin{equation}
\bI_1=\left[\begin{array}{cc}\bI & 0 \\ 0 & 0\end{array}\right]\quad\quad\bI_2=\left[\begin{array}{cc}0 & 0\\0 &\bI
\end{array}\right]\label{eq.idef}
\end{equation}
and $\|\bI_1\|=\|\bI_2\|=1$.

\paragraph{Lower Triangular Matrix:}
\begin{align}
\notag
&\left\|\left[\begin{array}{cc}0 & 0
\\
\nabla^2_{21}f(\bx) & 0 \end{array}\right]-\left[\begin{array}{cc} 0 & 0\\
\nabla^2_{21}f(\by) & 0\end{array}\right]\right\|
\\\notag
&=\left\|\bI_2\left(\left[\begin{array}{cc}\nabla^2_{11}f(\bx) & \nabla^2_{12}f(\bx)
\\
\nabla^2_{21}f(\bx) & \nabla^2_{22}f(\bx)\end{array}\right]-\left[\begin{array}{cc}\nabla^2_{11}f(\by) & \nabla^2_{12}f(\by)\\
\nabla^2_{21}f(\by) & \nabla^2_{22}f(\by)\end{array}\right]\right)\bI_1\right\|
\\\notag
&\mathop{\le}\limits^{(a)}\rho\|\bx-\by\|
\end{align}
where $(a)$ is true because we know $\|\bI_1\|=\|\bI_2\|=1$.

\subsection{Proof of \leref{le.layer001grad}}
\begin{proof}
Recall the definition
\begin{equation}\notag
\bh^{(t)}_{-1}\bydef\bx^{(t)}_2\quad \textrm{and}\quad \bh^{(t)}_{-2}\bydef\bx^{(t+1)}_1.
\end{equation}

First, we have
\begin{equation}
\|\nabla_2 f(\bx^{(t)}_1,\bx^{(t)}_2)\|^2\le 2\|\nabla_2f(\bx^{(t+1)}_1,\bx^{(t)}_2)-\nabla_2 f(\bx^{(t)}_1,\bx^{(t)}_2)\|^2+2\|\nabla_2 f(\bx^{(t+1)}_1,\bx^{(t)}_2)\|^2.
\end{equation}


Using block-wise Lipschitz continuity, we have
\begin{align}
\notag
\|\nabla_2 f(\bx^{(t)}_1,\bx^{(t)}_2)\|^2&\le2L^2_{\max}\|\bx^{(t+1)}_1-\bx^{(t)}_1\|^2+2\|\nabla_2f(\bx^{(t+1)}_1,\bx^{(t)}_2)\|^2
\\\notag
&\mathop{=}\limits^{(a)} 2L^2_{\max}\|\eta\nabla_1 f(\bx^{(t)}_1,\bx^{(t)}_2)\|^2+2\|\nabla_2f(\bx^{(t+1)}_1,\bx^{(t)}_2)\|^2
\\
&\mathop{\le}\limits^{(b)}2\sum^2_{k=1}\|\nabla_{k}f(\bh^{(t)}_{-k},\bx^{(t)}_{k})\|^2
\end{align}
where $(a)$ is because we use the update rule of AGD, $(b)$ is true due to $\eta\le 1/L_{\max}$.

Summing $\|\nabla_1 f(\bx^{(t)}_1,\bx^{(t)}_2)\|^2$ on both sides of the above equation, we have
\begin{equation}
\|\nabla f(\bx^{(t)})\|^2\le\sum^2_{k=1}\|\nabla_kf(\bx^{(t)}_k)\|^2\le 4\sum^2_{k=1}\|\nabla_{k}f(\bh^{(t)}_{-k},\bx^{(t)}_{k})\|^2.
\end{equation}
\end{proof}

\section{Proofs of PA-GD}\label{sec.pofpagd}

As stated in the main body of the paper, we can use \leref{le.descent} and \leref{le.escape} to prove \thref{th.mainth}. \leref{le.descent} is basically well-known. The main task focuses on proving \leref{le.escape}, which consists of a sequence of
lemmas (\leref{le.layer31}--\leref{le.layer21}) that lead to \leref{le.escape}.

Before discussing the details of \leref{le.escape}, we need to introduce some constants defined as follows,

\begin{align}
\notag
\mathcal{F}\bydef& \eta^5 L^5_{\max}\frac{\gamma^{3}}{\kappa^3\rho^2}\log^{-6}\left(\frac{d\kappa}{\delta}\right)\mathcal{P}_1^{-6}\mathcal{P}^{-2}_2,
\\\notag
\mathcal{G}\bydef&\eta^2 L^2_{\max}\frac{\gamma^2}{\rho}\log^{-3}\left(\frac{d\kappa}{\delta}\right)\mathcal{P}^{-3}_1\mathcal{P}_2^{-1},
\\\notag
\mathcal{S}\bydef&\eta^2 L^2_{\max}\frac{\gamma}{\kappa\rho}\log^{-2}\left(\frac{d\kappa}{\delta}\right)\mathcal{P}_1^{-2}\mathcal{P}_2^{-1},
\\\notag
\mathcal{T}\bydef&\frac{\log\left(\frac{d\kappa}{\delta}\right)\mathcal{P}_1}{\eta\gamma}.
\end{align}

These quantities
refer to different units of the algorithm. Specifically, $\mathcal{F}$
accounts for the objective value, $\mathcal{G}$ for the size of the gradient,
$\mathcal{S}$ for the norm of the difference between iterates, and $\mathcal{T}$ for the number
of iterations. Also, we define a condition number in terms of $\gamma$ as
$\kappa\bydef\frac{L_{\max}}{\gamma}\ge 1$.

These quantities, $\mathcal{F}$, $\mathcal{G}$, $\mathcal{S}$ and $\mathcal{T}$ have some certain relations as follows, which are useful of simplifying the expressions in the proofs.
\begin{subequations}
\begin{align}
\sqrt{\mathcal{F}}=&\frac{\sqrt{\eta}\mathcal{G}}{\kappa},\label{eq.fandg}
\\
\frac{\eta\mathcal{G}\mathcal{T}}{\kappa}=&\mathcal{S},\label{eq.gandt}
\\
\rho\mathcal{S}^3=&\frac{\eta L_{\max}\mathcal{F}}{\mathcal{P}_2},\label{eq.sandf}
\\
\eta\rho\mathcal{S}\mathcal{T}=&\frac{\eta^2 L^2_{\max}}{\kappa\log(\frac{d\kappa}{\delta})\mathcal{P}_1\mathcal{P}_2}.\label{eq.sandt}
\end{align}
\end{subequations}
In the process of the proofs, we used conditions $\log(\frac{d\kappa}{\delta})\ge1$, $\mathcal{P}_1\ge2$ repeatedly to simply the expressions of the parameters. We also consider saddle point $\widetilde{\bx}^{(t)}$ that satisfies the following condition.
\begin{condition}\label{cond}
An $\epsilon$-second order stationary point $\widetilde{\bx}^{(t)}$ satisfies the following conditions:
\begin{equation}
\sum^2_{k=1}\|\nabla_k f(\tbh^{(t)}_{-k},\wbx^{(t)}_k)\|^2\le g^2_{\textsf{th}}\quad\textrm{and}\quad\lambda_{\min}(\nabla^2 f(\wbx^{(t)}))\le -\gamma\label{eq.gcond}
\end{equation}
where $g_{\textsf{th}}\bydef\frac{\mathcal{G}}{2\kappa}$.
\end{condition}
\conref{cond} implies that point $\widetilde{\bx}^{(t)}$ satisfies $\|\nabla f(\widetilde{\bx}^{(t)})\|\le \mathcal{G}/\kappa$ (see \leref{le.layer001grad}) and $\lambda_{\min}(\nabla^2 f(\wbx^{(t)}))\le -\gamma$.

\paragraph{Sufficient Decrease after Perturbation}


Consider $\wbx^{(t)}$ satisfy \conref{cond} and let $\bH\triangleq\nabla^2f(\widetilde{\bx}^{(t)})$. We consider a second order
approximation as the following
\begin{equation}
\widehat{f}_{\by}(\bx)\triangleq f(\by)+\nabla f(\by)^{\T}(\bx-\by)+\frac{1}{2}(\bx-\by)^{\T}\bH(\bx-\by).\label{eq.defoffhat}
\end{equation}

With these definitions of parameters, we will study how PA-GD can escape from
strict saddle points. The main part of the proof is to
show that when two sequences are apart from each other with a certain distance
along the $\vec{\be}$ direction at the starting points, where $\vec{\be}$ denotes the eigenvector of $\bM^{-1}\bT$ whose eigenvalue is maximum (greater than 1). Then, after a number of
iterations at least one of them can give a sufficient decrease of the objective
value. This property implies the iterates can easily escape from the saddle
points as long as there is a large enough perturbation between the initial
points of the two sequences along the $\vec{\be}$ direction. We will introduce the
following two lemmas formally which are the main contributions of this work.
\begin{lemma}\label{le.layer31}
Under \asref{as1}, consider $\widetilde{\bx}^{(t)}$ that satisfies \conref{cond} and a generic sequence $\bu^{(t)}$ generated by AGD.
For any constant $\widehat{c}\ge 2$, $\delta\in(0,\frac{d\kappa}{e}]$, when
initial point $\bu^{(0)}$ satisfies
\begin{equation}
\|\bu^{(0)}-\widetilde{\bx}^{(t)}\|\le 2r,
\end{equation}
then, with the definition of
\begin{equation}
r\bydef \frac{\eta L_{\max}\mathcal{S}}{\kappa\log(\frac{d\kappa}{\delta})\mathcal{P}_1},
\quad\textrm{and}\quad T\bydef\min\{\inf_t\{t|\widehat{f}_{\bu^{(0)}}(\bu^{(t)})-f(\bu^{(0)})\le-3\mathcal{F}\},\widehat{c}\mathcal{T}\},\label{eq.defoft}
\end{equation}
there exits constants $c^{(1)}_{\max}, \widehat{c}$ such that for any $\eta\le
c^{(1)}_{\max}/L_{\max}$, the iterates generated by PA-GD satisfy $\|\bu^{(t)}-\widetilde{\bx}^{(t)}\|\le5\widehat{c}\mathcal{S}, \forall t<T$.
\end{lemma}
\begin{lemma}\label{le.layer32}
Under \asref{as1}, consider $\widetilde{\bx}^{(t)}$ that satisfies \conref{cond}.
There exist constants $c^{(2)}_{\max}$, $\widehat{c}$ such that: for any
$\delta\in(0,\frac{d\kappa}{e}]$ and $\eta\le c^{(2)}_{\max}/L_{\max}$, with the
definition of
\begin{equation}\notag
T\bydef\min\left\{\inf_t\{t|\widehat{f}_{\bw_0}(\bw^{(t)})-f(\bw^{(0)})\le-3\mathcal{F}\},\widehat{c}\mathcal{T}\right\}
\end{equation}
where two iterates $\{\bu^{(t)}\}$ and $\{\bw^{(t)}\}$ that are generated
by PA-GD with initial points $\{\bu^{(0)}, \bw^{(0)}\}$ satisfying
\begin{equation}
\|\bu^{(0)}-\widetilde{\bx}^{(t)}\|\le r,\;\bw^{(0)}=\bu^{(0)}+\upsilon r \vec{\be},\;\upsilon\in[\delta/(2\sqrt{d}),1],\label{eq.inicond}
\end{equation}
where $\vec{\be}$ denotes the eigenvector of $\bM^{-1}\bT$ whose eigenvalue is maximum, then, if $\|\bu^{(t)}-\widetilde{\bx}^{(t)}\|\le5\widehat{c}\mathcal{S}, \forall t<T$, we will have $T<\widehat{c}\mathcal{T}$.
\end{lemma}

\leref{le.layer31} says that if the $\bu^{(t)}$-iterate generated by PA-GD
cannot provide a sufficient decrease of the objective value, then the iterates
are constrained within the area which is very close to the saddle point. With
this property, \leref{le.layer32} shows if there exists another PA-GD iterate
$\bw^{(t)}$, which is initialized with a certain distance along the $\vec{\be}$
direction from the $\bu$-iterate, then $\bw^{(t)}$ will provide a sufficient
decrease of the objective value. These two lemmas characterize the convergence
behavior of the PA-GD iterates.
%

\paragraph{Escaping from Saddle Points}

Then, we need to quantify the probability that after adding the perturbation the algorithm cannot escape from strict saddle points.
In previous work about escaping from
saddle points with GD, a characterization of the geometry around saddle points
has been given \citep[Lemma 15]{jin2017jordan}. Once we know that PA-GD also decreases the
objective value sufficiently in \leref{le.layer31} and \leref{le.layer32}, the
following lemma can be claimed straightforwardly. To be more specific, we can obtain the probability that iterates will be stuck at the strict points after $T$ iterations as follows.
\begin{align}
\notag
\mathbb{P}(\bw^{(0)}\in\mathcal{X}_{\textrm{stuck}})=&\int_{\mathbb{B}_{\wbx^{(t)}}(r)}\mathbb{P}(\bw^{(0)}\in\mathcal{X}_{\textrm{stuck}}|\bu^{(0)}\in\mathcal{X}_{\textrm{stuck}})\mathbb{P}(\bu^{(0)}\in\mathcal{X}_{\textrm{stuck}})d\bu^{(0)}
\\\notag
\le&\int_{\mathbb{B}_{\wbx^{(t)}}(r)}\mathbb{P}(\bw^{(0)}\in\mathcal{X}_{\textrm{stuck}}|\bu^{(0)}\in\mathcal{X}_{\textrm{stuck}})\mathbb{P}(\bu^{(0)})d\bu^{(0)}
\\\notag
\mathop{\le}\limits^{(a)}& \delta\int_{\mathbb{B}_{\wbx^{(t)}}(r)}\mathbb{P}(\bu^{(0)})d\bu^{(0)}=\delta
\end{align}
where $\mathcal{X}_{\textrm{stuck}}$ denotes the set where the algorithm starts such that the sequence cannot escape from the strict saddle point after $T$ iterations, $(a)$ is true because probability $\mathbb{P}(\bw^{(0)}\in\mathcal{X}_{\textrm{stuck}}|\bu^{(0)}\in\mathcal{X}_{\textrm{stuck}})$ can be upper bounded by $\delta$, which is proven in the following lemma.

\begin{lemma}\label{le.layer21}
Under \asref{as1}, there exists a universal constant $c_{\max}$, for any
$\delta\in(0,d\kappa/e]$: consider a saddle point $\wbx^{(t)}$ which satisfies
\conref{cond}, let $\bx^{(0)}=\wbx^{(t)}+\xi$ where $\xi$ is generated randomly
which follows the uniform distribution over a ball with radius
$r$, and let $\bx^{(t)}$ be the
iterates of PA-GD starting from $\bx^{(0)}$. Then, when step size $\eta\le
c_{\max}/L_{\max}$, with at least probability $1-\delta$, we have the
following for any $T\ge\mathcal{T}/c_{\max}$
\begin{equation}
f(\bx^{(T)})-f(\wbx^{(t)})\le-\mathcal{F}.
\end{equation}
\end{lemma}

Then, applying $\eta=\frac{c}{L_{\max}}$,$\gamma=(L_{\max}\rho\epsilon)^{1/3}$, and
$\delta=\frac{dL_{\max}}{(L_{\max}\rho\epsilon)^{1/3}}e^{-\chi}$ into \leref{le.layer21}, we can get
\leref{le.escape} immediately.

With these lemmas, we can give the proof of \thref{th.mainth} as the following.
\subsection{Proof of \thref{th.mainth}}

Next, we prove the main theorem.
\begin{proof}

Submitting $\eta=\frac{c}{L_{\max}}$,$\gamma=(L_{\max}\rho\epsilon)^{1/3}$, and
$\delta=\frac{dL_{\max}}{(L_{\max}\rho\epsilon)^{1/3}}e^{-\chi}$ into the definitions of $\mathcal{F},\mathcal{G}, \mathcal{T}$, we will have the following definitions.
\begin{align}
\notag
f_{\textsf{th}}\bydef&\mathcal{F}=\frac{c^5\epsilon^{2}}{L_{\max}(\chi\mathcal{P}_1)^6\mathcal{P}^2_2},
\\\notag
g_{\textsf{th}}\bydef&\frac{\mathcal{G}}{2\kappa}=\frac{c^2\epsilon}{2(\chi\mathcal{P}_1)^3\mathcal{P}_2},
\\\notag
t_{\textsf{th}}\bydef&\frac{\mathcal{T}}{c}=\frac{L_{\max}\chi\mathcal{P}_1}{c^2(L_{\max}\rho\epsilon)^{\frac{1}{3}}}.
\end{align}

After applying \leref{le.layer001grad}, we know that
\begin{equation}
\notag
\|\nabla f(\bx)\|\le \frac{c}{\chi^{3}\mathcal{P}^3_1\mathcal{P}_2}\epsilon
\end{equation}
where $c\le 1,\chi,\mathcal{P}_1,\mathcal{P}_2\ge 1$.

With a set of necessary lemmas and leveraging the proof of PGD \citep[Theorem 3]{jin2017jordan}, we have the following convergence analysis of PA-GD. Specifically, at any iteration, we need to
consider two cases (we use the first iteration as an example):
\begin{enumerate}
\item In this case the gradient is large such that $\sum^2_{k=1}\|\nabla_k
f(\bh^{(0)}_{-k},\bx^{(0)}_k)\|^2>g^2_{\textsf{th}}$: According to
\leref{le.descent}, we have
\begin{align}
\notag
f(\bx^{(1)})-f(\bx^{(0)})&\le-\sum^2_{k=1}\frac{\eta}{2}\|\nabla_k f(\bh^{(0)}_{-k},\bx^{(0)}_k)\|^2\le-\frac{\eta}{2}g^2_{\textsf{th}}
\\
&\mathop{=}\limits^{(a)}-\frac{c^5}{8(\chi\mathcal{P}_1)^6\mathcal{P}_2^2}\frac{\epsilon^2}{L_{\max}}\label{eq.gradd}
\end{align}
where in $(a)$  use the definition of $g^2_{\textsf{th}}$ and $\eta\le c/L_{\max}$.

\item The gradient is small in all block directions, namely $\sum^2_{k=1}\|\nabla_k
f(\bh^{(0)}_{-k},\bx^{(0)}_k)\|^2\le g^2_{\textsf{th}}$: in this case, we will
add the perturbation to the iterates, and implement AGD for the next
$t_{\textsf{th}}$ steps and then check the termination condition. If the
termination condition is not satisfied, we must have
\begin{equation}
f(\bx^{(t_{\textsf{th}})})-f(\bx^{(0)})\le-f_{\textsf{th}}=-\frac{c^5\epsilon^{2}}{L_{\max}(\chi\mathcal{P}_1)^6\mathcal{P}^2_2},
\end{equation}
which implies that the objective value in each step on average is decreased by
\begin{equation}
\frac{f(\bx^{(t_{\textsf{th}})})-f(\bx^{(0)})}{t_{\textsf{th}}}\le-\frac{c^{7}}{(\chi\mathcal{P}_1)^7\mathcal{P}^2_2}\frac{\epsilon^2}{L_{\max}}\frac{(L_{\max}\rho\epsilon)^{\frac{1}{3}}}{L_{\max}}.\label{eq.gradp}
\end{equation}
Since $\kappa=L_{\max}/(L_{\max}\rho\epsilon)^{1/3}\ge1$,  we know that the right-hand side (RHS) of \eqref{eq.gradp} is greater than RHS of \eqref{eq.gradd}.

With the results of these two cases, we can know that if there is a large size of the gradient, we can know the decrease of the objective function value by the result of case 1, and if not, we use the result of case 2. In summary, PA-GD can have a sufficient decrease of the objective function value by $\frac{c^{7}}{(\chi\mathcal{P}_1)^7\mathcal{P}^2_2}\frac{\epsilon^2}{L_{\max}}\frac{(L_{\max}\rho\epsilon)^{1/3}}{L_{\max}}$ per iteration on average. This means that \algref{alg:p1} must stop within a finite number of iterations, which is
\begin{equation}
\frac{f(\bh^{(0)}_{-1},\bx^{(0)}_1)-f^*}{\frac{c^{7}}{(\chi\mathcal{P}_1)^7\mathcal{P}^2_2}\frac{\epsilon^2}{L_{\max}}\frac{(L_{\max}\rho\epsilon)^{1/3}}{L_{\max}}}
=\frac{(\chi\mathcal{P}_1)^{7}\mathcal{P}^2_2}{c^{7}}\frac{L^2_{\max}\Delta f}{\epsilon^2(L_{\max}\rho\epsilon)^{1/3}}
=\mathcal{O}\left(\frac{\Delta f(\chi\mathcal{P}_1)^7\mathcal{P}^2_2L^{5/3}_{\max}}{\rho^{1/3}\epsilon^{7/3}}\right)
\end{equation}
where $\Delta f\bydef f(\bh^{(0)}_{-1},\bx^{(0)}_1)-f^*$.

According to \leref{le.escape}, we know that with probability $1-\frac{dL_{\max}}{(L_{\max}\rho\epsilon)^{1/3}}e^{-\chi}$ the algorithm can give a sufficient descent with the perturbation when $\sum^2_{k=1}\|\nabla_k
f(\bh^{(t)}_{-k},\bx^{(t)}_k)\|^2\le g^2_{\textsf{th}}$. Since the total number of perturbation we can add is at most
\begin{equation}
n=\frac{1}{t_{\textsf{th}}}\frac{(\chi\mathcal{P}_1)^{7}\mathcal{P}^2_2}{c^{7}}\frac{L^2_{\max}\Delta f}{\epsilon^2(L_{\max}\rho\epsilon)^{1/3}}
=\frac{(\mathcal{P}_1\chi)^{6}\mathcal{P}^2_2}{c^5}\frac{L_{\max}\Delta_f}{\epsilon^2}.
\end{equation}
Using the union bound, the probability of \leref{le.escape} being satisfied for all perturbations is
\begin{equation}\label{eq.bdoferror}
1-n\frac{dL_{\max}}{(L_{\max}\rho\epsilon)^{\frac{1}{3}}}e^{-\chi}=1-\frac{dL_{\max}}{(L_{\max}\rho\epsilon)^{\frac{1}{3}}}e^{-\chi}\frac{(\mathcal{P}_1\chi)^{6}\mathcal{P}^2_2}{c^5}\frac{L_{\max}\Delta_f}{\epsilon^2}
=1-\underbrace{\frac{dL_{\max}}{(L_{\max}\rho\epsilon)^{\frac{1}{3}}}\frac{\mathcal{P}^6_1\mathcal{P}^2_2}{c^5}\frac{\Delta_f}{\epsilon^2}}_{\bydef\mathcal{C}}\chi^{6}e^{-\chi}.
\end{equation}
With chosen $\chi=6\max\{\ln(\mathcal{C}/\delta),4\}$, we have $\chi^6e^{-\chi}\le e^{-\chi/6}$, which implies $\chi^6e^{-\chi}\mathcal{C}\le e^{-\chi/6}\mathcal{C}\le\delta$.
\end{enumerate}
The proof is complete.
\end{proof}

\subsection{Proof of \leref{le.eiginver}}

\begin{proof}
Recall the definitions:
\begin{equation}
\bH_u\bydef\left[\begin{array}{cc}\nabla^2_{11} f(\wbx^{(t)}) &  \nabla^2_{12} f(\wbx^{(t)})
\\
0 &   \nabla^2_{22} f(\wbx^{(t)})\end{array}\right]\quad\bH_l\bydef\left[\begin{array}{cc}0 & 0
\\
\nabla^2_{21} f(\wbx^{(t)}) & 0 \end{array}\right],
\end{equation}
where $\wbx^{(t)}$ is an $\epsilon$-second order stationary point, and
\begin{equation}
\bM\bydef\bI+\eta\bH_l,\quad\bT\bydef\bI-\eta\bH_u.
\end{equation}

Our goal of this lemma is to show that the maximum eigenvalue of $\bM^{-1}\bT$ is greater than 1 so that we can project iterates $\bv^{(t)}$ onto the two subspaces, where the first subspace is spanned by the eigenvector of $\bM^{-1}\bT$ whose eigenvalue is the largest (greater than 1) and the other one is spanned by the remaining eigenvectors.

Note that $\det(\bM)=1$, which implies that $\det(\bM^{-1}\bT-\lambda\bI)=\det(\bT-\lambda\bM)$, where $\lambda$ denotes the eigenvalue. We can analyze the determinant of $\bT-\lambda\bM$, i.e.,
\begin{align}\notag
\det[\bT-\lambda\bM]=&\det[\bI-\eta\bH_u-\lambda(\bI+\eta\bH_l)]
\\
=&\det\left[\underbrace{\begin{array}{cc}(1-\lambda)\bI-\eta\nabla^2_{11}f(\wbx^{(t)}) & -\eta\nabla^2_{12} f(\wbx^{(t)})
\\\notag
-\lambda\eta\nabla^2_{21}f(\wbx^{(t)})  & (1-\lambda)\bI-\eta\nabla^2_{22}f(\wbx^{(t)})
\end{array}}_{\bydef\bQ(\lambda)}\right].
\end{align}

Then, we use two steps to show $\lambda_{\max}(\bM^{-1}\bT)>1$: 1) we can show that all eigenvalues of $\bQ(\lambda)$ are real; 2) there exists a  $\lambda>1$ such that $\det(\bQ(\lambda))=0$.

Consider a $\delta>0$. We have
\begin{equation}
\bQ(1+\delta)=-\left(\underbrace{\eta\bH+\delta(\bI+\eta\bH_l)}_{\bydef\bF(\delta)}\right)\label{eq.defofq}
\end{equation}
where
\begin{align}
\notag
\bF(\delta)=&\delta\bI+\eta\left[\begin{array}{cc}
\nabla^2_{11}f(\wbx^{(t)})  &  \nabla^2_{12}f(\wbx^{(t)})
\\\notag
(1+\delta)\nabla^2_{21}f(\wbx^{(t)})  & \nabla^2_{22}f(\wbx^{(t)})
\end{array}\right]
\\\notag
=&\left[\begin{array}{cc}\bI &  \\ & \sqrt{1+\delta}\end{array}\right]\underbrace{\left[\begin{array}{cc}\delta\bI+\eta\nabla^2_{11}f(\wbx^{(t)})  & \eta\sqrt{1+\delta}\nabla^2_{12} f(\wbx^{(t)})
\\\notag
\eta\sqrt{1+\delta}\nabla^2_{21}f(\wbx^{(t)})  & \delta\bI+\eta\nabla^2_{22}f(\wbx^{(t)})
\end{array}\right]}_{\bG(\delta)}\left[\begin{array}{cc}\bI &  \\ & \frac{1}{\sqrt{1+\delta}} \end{array}\right],
\end{align}
meaning that  $\bF(\delta)$  is similar to $\bG(\delta)$. Consequently, we can conclude that $\bF(\delta)$ has the same eigenvalues of $\bG(\delta)$. Since we know that $\bH$ and $\bG(\delta)$ are diagonalizable (normal matrices), then we have the following result \citep{weyl1912asymptotische} (or \citep{holbrook1992spectral}) of quantifying the difference of the eigenvalues of the two normal matrices
\begin{equation}
\max_{1\le i\le d}|\lambda_i(\eta\bH)-\lambda_i(\bG(\delta))|\le\|\eta\bH-\bG(\delta)\|\label{eq.rehg}
\end{equation}
where $\lambda_i(\bH)$ and $\lambda_i(\bG(\delta))$ denote the $i$th eigenvalue of $\bH$ and $\bG(\delta)$, which are listed in a decreasing order.

With the help of \eqref{eq.rehg}, we can check
\begin{align}
\notag
&\|\eta\bH-\bG(\delta)\|
\\\notag
=&\left\|\delta\bI+\left[\begin{array}{cc}0 & (\sqrt{1+\delta}-1)\eta\nabla^2_{12}f(\wbx^{(t)})
\\\notag
(\sqrt{1+\delta}-1)\eta\nabla^2_{21}f(\wbx^{(t)}) & 0
\end{array}\right]\right\|
\\\notag
\le&\delta+(\sqrt{1+\delta}-1)\eta\|\bH\|+(\sqrt{1+\delta}-1)\eta\left\|\begin{array}{cc}\nabla^2_{11}f(\wbx^{(t)}) &0 \\ 0 &\nabla^2_{22}f(\wbx^{(t)})\end{array}\right\|
\\
\mathop{\le}\limits^{(a)} & \delta+(\sqrt{1+\delta}-1)(\frac{L}{L_{\max}}+1).\label{eq.diffeig}
\end{align}
where $(a)$ is true since we used $\eta\le c_{\max}/L_{\max}$ and the fact that $\|\bH\|\le L$ and $\|\bH_d\|\le L_{\max}$. Also, it can be observed that when $\delta=0$, matrix $\bG(\delta)$ is reduced to $\eta\bH$. Note that if $\eta=1/L$ is used, then we have $\|\eta\bH-\bG(\delta)\|\le\delta+2(\sqrt{1+\delta}-1)$.

We know that the minimum eigenvalue of $\eta\bH$ which is $-\eta\gamma$ and the maximum difference of the eigenvalues between  $\eta\bH$ and $\bG(\delta)$ is upper bounded by \eqref{eq.diffeig}. Then, we can choose a sufficient small $\delta$ such that $\bG(\delta)$ also has a negative eigenvalue, meaning that we need to find a $\delta$ such that
\begin{align}
\delta+(\sqrt{1+\delta}-1)(\frac{L}{L_{\max}}+1)<\eta\gamma. \label{eq.vadelta}
\end{align}

In other words, if we choose
\begin{equation}\notag
\delta^*=\frac{\eta\gamma}{1+\frac{L}{L_{\max}}}
\end{equation}
then we can conclude that $\bG(\delta^*)$ has a negative eigenvalue which is less than $-\eta\gamma+\delta^*=-\frac{\eta\gamma}{1+\frac{L_{\max}}{L}}$.

In the following, we will check that $\delta^*$ is a valid choice, meaning that equation \eqref{eq.vadelta} holds when $\delta^*=\frac{\eta\gamma}{1+\frac{L}{L_{\max}}}$.

\paragraph{First step}: since $L/L_{\max}\ge1$, we have $\eta\gamma/(1+L/L_{\max})\le\eta\gamma/2$.

\paragraph{Second step}: we only need to check
\begin{equation}\notag
(\sqrt{1+\delta}-1)(\frac{L}{L_{\max}}+1)<\frac{\eta\gamma}{2},
\end{equation}
meaning that it is sufficient to check
\begin{equation}
(\frac{L}{L_{\max}}+1)^2(1+\delta)\le\left(\frac{L}{L_{\max}}+1+\frac{\eta\gamma}{2}\right)^2\label{eq.deltaopt}.
\end{equation}

It can be easily check that the left-hand side (LHS) of \eqref{eq.deltaopt} with chosen $\delta^*$ is
\begin{equation}\notag
(\frac{L}{L_{\max}}+1)^2(1+\frac{\eta\gamma}{\frac{L}{L_{\max}}+1})\le(\frac{L}{L_{\max}}+1)^2+(\frac{L}{L_{\max}}+1)^2\eta\gamma<(\frac{L}{L_{\max}}+1)^2+(\frac{L}{L_{\max}}+1)^2\eta\gamma+\frac{\eta^2\gamma^2}{4},
\end{equation}
which is RHS of \eqref{eq.deltaopt}.

Therefore, we can conclude that $\bQ(1+\delta^*)$ has a negative eigenvalue.

When $\delta$ is large, it is easy to check $\bQ(1+\delta)$ has a positive eigenvalue, since term $\delta^2\bI$ dominates the spectrum of matrix $\bQ(1+\delta)$ in \eqref{eq.defofq}. Since the eigenvalue is continuous with respect to $\delta$, we can conclude there exists a largest $\delta$, i.e., $\widehat{\delta}$,  such that $\bQ(1+\widehat{\delta})$ has a zero eigenvalue, i.e., $\det(\bQ(1+\widehat{\delta}))=0$ where $1+\widehat{\delta}$ is at least
\begin{equation}
1+\delta^*=1+\frac{\eta\gamma}{L/L_{\max}+1}.\label{eq.gammap}
\end{equation}
Therefore, we can conclude that there exits a largest real eigenvalue of $\bM^{-1}\bT$ which is  $1+\wde>1+\delta^*>1$.
\end{proof}

\subsection{Proof of \leref{le.descent}}
\begin{proof}
Under \asref{as1}, we have (descent lemma)
\begin{align}
\notag
f(\bx^{(t+1)})\le& f(\bx^{(t)})+\sum^2_{k=1}\nabla_k f(\bh^{(t)}_{-k},\bx^{(t)}_k)^{\T}(\bx^{(t+1)}_{k}-\bx^{(t)}_k)+\sum^2_{k=1}\frac{L_k}{2}\|\bx^{(t+1)}_k-\bx^{(t)}_k\|^2
\\\notag
\mathop{\le}\limits^{(a)}&f(\bx^{(t)})-\sum^2_{k=1}\eta\|\nabla_k f(\bh^{(t)}_{-k},\bx^{(t)}_k)\|^2+\sum^2_{k=1}\frac{\eta^2L_k}{2}\|\nabla_k f(\bh^{(t)}_{-k},\bx^{(t)}_k)\|^2
\\
\mathop{\le}\limits^{(b)}&f(\bx^{(t)})-\sum^2_{k=1}\frac{\eta}{2}\|\nabla_k f(\bh^{(t)}_{-k},\bx^{(t)}_k)\|^2\label{eq.desbcd}
\end{align}
where (a) is true because of the update rule of gradient descent in each block
and \asref{as1}, in (b) we used $\eta\le1/L_{\max}$.
\end{proof}

\subsection{Proof of \leref{le.layer31}}
\begin{proof}



Without loss of generality, let $\bu^{(0)}$ be the origin, i.e.,
$\bu^{(0)}=0$. According to the AGD update rules, we have
\begin{align}
\bu^{(t+1)}=&\bu^{(t)}-\eta\left[\begin{array}{c}\nabla_1 f(\bu^{(t)}_1,\bu^{(t)}_2)\\ \nabla_2 f(\bu^{(t+1)}_1,\bu^{(t)}_2)\end{array}\right]\label{eq.bcdupdateo}.
\end{align}

Then, we use the mathematical induction to prove that
\begin{equation}
\|\bu^{(t)}\|\le 5\widehat{c}\mathcal{S}\label{eq.boundofu}, \forall t<T.
\end{equation}

When $t=0$, we have $\bu^{(0)}=0$, so \eqref{eq.boundofu} is true.

Suppose
\eqref{eq.boundofu} is true for the case where $\tau\le t$. We will show that
\eqref{eq.boundofu} is also true for the case where $\tau=t+1$.

First, we need to show the upper bound of $\|\bu^{(t+1)}-\bu^{(t)}\|$. According to the Taylor expansion and $\rho$-Hessian Lipschitz continuity, we have
\begin{equation}
\notag
f(\bu^{(t)})\le f(\bu^{(0)})+\nabla f(\bu^{(0)})^{\T}(\bu^{(t)}-\bu^{(0)})+\frac{1}{2}(\bu^{(0)}-\bu^{(t)})^{\T}\nabla^2 f(\bu^{(0)})(\bu^{(0)}-\bu^{(t)})+\frac{\rho}{6}\|\bu^{(t)}-\bu^{(0)}\|^3.
\end{equation}

Comparing with the definition of $\widehat{f}_{\bu^{(0)}}(\bu^{(t)})$, we have
\begin{align}
\notag
|f(\bu^{(t)})-\widehat{f}_{\bu^{(0)}}(\bu^{(t)})|\mathop{\le}\limits^{\eqref{eq.defoffhat}}&\frac{1}{2}(\bu^{(0)}-\bu^{(t)})^{\T}\left(\nabla^2 f(\bu^{(0)})-\bH\right)(\bu^{(0)}-\bu^{(t)})+\frac{\rho}{6}\|\bu^{(t)}-\bu^{(0)}\|^3
\\\notag
\mathop{\le}\limits^{(a)}&\frac{\rho}{2}\|\bu^{(0)}-\widetilde{\bx}^{(t)}\|\|\bu^{(t)}-\bu^{(0)}\|^2+\frac{\rho}{6}\|\bu^{(t)}-\bu^{(0)}\|^3
\end{align}
where in $(a)$ we also used $\rho$-Hessian Lipschitz continuity.

According to the definition of $T$, we know that $f(\bu^{(0)})-\widehat{f}_{\bu^{(0)}}(\bu^{(t)})\le 3\mathcal{F}$ for all $t<T$, which implies that
\begin{align}
\notag
f(\bu^{(0)})-f(\bu^{(t)})\le& |f(\bu^{(0)})-\widehat{f}_{\bu^{(0)}}(\bu^{(t)})|+|\widehat{f}_{\bu^{(0)}}(\bu^{(t)})-f(\bu^{(t)})|
\\\notag
\mathop{\le}\limits^{\eqref{eq.defoft}}& 3\mathcal{F}+\frac{\rho}{2}\|\widetilde{\bx}^{(t)}-\bu^{(0)}\|\|\bu^{(t)}-\bu^{(0)}\|^2+\frac{\rho}{6}\|\bu^{(t)}-\bu^{(0)}\|^3
\\
\le& 3\mathcal{F}+\frac{\rho}{2}\frac{\eta L_{\max}\mathcal{S}}{\kappa \log(\frac{d\kappa}{\delta})\mathcal{P}_1}(5\widehat{c}\mathcal{S})^2+\frac{\rho}{6}(5\widehat{c}\mathcal{S})^3\label{eq.bdrho}
\\\notag
\le&3\mathcal{F}+((5\widehat{c})^2/4+(5\widehat{c})^3/6)\rho\mathcal{S}^3
\\
\mathop{\le}\limits^{\eqref{eq.sandf}}&3\mathcal{F}+\eta L_{\max}(5\widehat{c})^3\mathcal{F}\mathcal{P}^{-1}_2 \label{eq.fdes}
\\
\le& 4\mathcal{F}\label{eq.diffbudes}
\end{align}
where in \eqref{eq.diffbudes} we used $c_{\max}=\mathcal{P}_2/(5\widehat{c})^3$ and $\eta\le c_{\max}/L_{\max}$.

From \eqref{eq.desbcd}, we also know that
\begin{equation}
f(\bu^{(t+1)})\le f(\bu^{(t)})-\frac{\eta}{2}\left(\|\nabla_1 f(\bu^{(t)}_1,\bu^{(t)}_2)\|^2+\|\nabla_2 f(\bu^{(t+1)}_1,\bu^{(t)}_2)\|^2\right),\quad \forall t<T.\label{eq.desbcdu}
\end{equation}

For simplification of expression, we define
\begin{equation}
\bz^{(t)}_{-1}\bydef\bu^{(t)}_2\quad \textrm{and}\quad \bz^{(t)}_{-2}\bydef\bu^{(t+1)}_1,\quad\forall t<T.
\end{equation}

Summing up \eqref{eq.desbcdu} for $\tau=0,\ldots,t$, we have
\begin{equation}
f(\bu^{(t)})\le f(\bu^{(0)})-\sum^{t-1}_{\tau=0}\sum^2_{k=1}\frac{\eta}{2}\|\nabla_k f(\bz^{(\tau)}_{-k},\bu^{(\tau)}_k)\|^2,\quad\forall t<T.\label{eq.desbd}
\end{equation}

Combining \eqref{eq.diffbudes} and \eqref{eq.desbd}, we know that
\begin{equation}
\sum^{t-1}_{\tau=0}\sum^2_{k=1}\frac{\eta}{2}\|\nabla_k f(\bz^{(\tau)}_{-k},\bu^{(\tau)}_k)\|^2\le4\mathcal{F},
\end{equation}
which implies
\begin{equation}
\max_{\tau}\sum^2_{k=1}\frac{\eta}{2}\|\nabla_k f(\bz^{(\tau)}_{-k},\bu^{(\tau)}_k)\|^2\le4\mathcal{F}, \tau\le t-1.\label{eq.bdofgu}
\end{equation}

According to \eqref{eq.bcdupdateo}, we know
\begin{align}
\notag
&\|\bu^{(t+1)}-\bu^{(t)}\|^2
\\\notag
=&\eta^2\sum^{2}_{k=1}\|\nabla_k f(\bz^{(t)}_{-k},\bu^{(t)}_k)\|^2
\\\notag
=&2\eta^2\sum^2_{k=1}\|\nabla_k f(\bz^{(t)}_{-k},\bu^{(t)}_k)-\nabla_k f(\bz^{(t-1)}_{-k},\bu^{(t-1)}_k)\|^2+2\eta^2\sum^2_{k=1}\|\nabla_k f(\bz^{(t-1)}_{-k},\bu^{(t-1)}_k)\|^2
\\\notag
=&2\eta^2\left(2\sum^2_{k=1}\|\nabla_k f(\bz^{(t)}_{-k},\bu^{(t)}_k)-\nabla_k
f(\bz^{(t-1)}_{-k},\bu^{(t)}_k)\|^2+2\sum^2_{k=1}\|\nabla_k f(\bz^{(t-1)}_{-k},\bu^{(t)}_k)-\nabla_k
f(\bz^{(t-1)}_{-k},\bu^{(t-1)}_k)\|^2\right)
\\\notag
&+2\eta^2\sum^2_{k=1}\|\nabla_k f(\bz^{(t-1)}_{-k},\bu^{(t-1)}_k)\|^2
\\\notag
\mathop{\le}\limits^{(a)} &8\eta^2L^2_{\max} \|\bu^{(t+1)}-\bu^{(t)}\|^2+4\eta^2L^2_{\max}\|\bu^{(t)}-\bu^{(t-1)}\|^2+16\eta\mathcal{F}.
\end{align}
where in $(a)$ we used Lipschitz continuity, i.e., $\sum^2_{k=1}\|\nabla_k f(\bz^{(t)}_{-k},\bu^{(t)}_k)-\nabla_k
f(\bz^{(t-1)}_{-k},\bu^{(t)}_k)\|^2\le L^2_{\max}\|\bu^{(t+1)}_1-\bu^{(t)}_1\|^2+L^2_{\max}\|\bu^{(t)}_2-\bu^{(t-1)}_2\|^2$, and $\sum^2_{k=1}\|\nabla_k f(\bz^{(t-1)}_{-k},\bu^{(t)}_k)-\nabla_k
f(\bz^{(t-1)}_{-k},\bu^{(t-1)}_k)\|^2\le L^2_{\max}\|\bu^{(t+1)}_1-\bu^{(t)}_1\|$.

Then, we have
\begin{align}
\notag
\|\bu^{(t+1)}-\bu^{(t)}\|^2\le& \underbrace{\frac{4\eta^2L^2_{\max}}{(1-8\eta^2L^2_{\max})}}_{\bydef\omega}\|\bu^{(t)}-\bu^{(t-1)}\|^2+\frac{16\eta\mathcal{F}}{(1-8\eta^2L^2_{\max})}
\\\notag
=&\omega^t\|\bu^{(1)}-\bu^{(0)}\|^2+\sum^{t-1}_{\tau=0}\omega^{\tau}\frac{16\eta\mathcal{F}}{(1-8\eta^2L^2_{\max})}
\\\notag
\mathop{\le}\limits^{(a)}&\frac{1-\omega^t}{1-\omega}\frac{16\eta\mathcal{F}}{(1-8\eta^2L^2_{\max})}\le\frac{1}{1-\omega}\frac{16\eta\mathcal{F}}{(1-8\eta^2L^2_{\max})}<1.14*16\eta\mathcal{F}<18.2\eta\mathcal{F}
\end{align}
where $(a)$ is true because we have $\|\bu^{(1)}-\bu^{(0)}\|^2\le16\eta\mathcal{F}$ since $t<T$ and \eqref{eq.bdofgu}, and we used $\eta\le c'_{\max}/L_{\max}$ where $c'_{\max}=1/10$ such that $\omega\approx0.0435<1$.

Then, we can obtain
\begin{equation}
\|\bu^{(t+1)}-\bu^{(t)}\|\le4.3\sqrt{\eta\mathcal{F}}\mathop{\le}\limits^{\eqref{eq.fandg}} \frac{4.3\eta \mathcal{G}}{\kappa}\label{eq.bdofdiffu}.
\end{equation}

Based on \eqref{eq.bdofdiffu}, we can get the upper bound of the sum of $\|\bu^{(t+1)}-\bu^{(t)}\|,\forall t<T$ as the following,
\begin{equation}
\sum^{t+1}_{\tau=1}\|\bu^{(\tau)}-\bu^{(\tau-1)}\|\le\sqrt{t\sum^{t+1}_{\tau=1}\|\bu^{(\tau)}-\bu^{(\tau-1)}\|^2}
\mathop{\le}\limits^{\eqref{eq.bdofdiffu}}T\cdot\frac{4.3\eta\mathcal{G}}{\kappa}\le \widehat{c}\mathcal{T}\frac{4.3\eta\mathcal{G}}{\kappa}\mathop{\le}\limits^{\eqref{eq.gandt}}4.3\widehat{c}\mathcal{S},
\end{equation}
which implies
\begin{equation}
\|\bu^{(t+1)}\|\mathop{\le}\limits^{(a)}\sum^{t+1}_{\tau=1}\|\bu^{(\tau)}-\bu^{(\tau-1)}\|+\|\bu^{(0)}\|\le4.3\widehat{c}\mathcal{S}
\end{equation}
where in $(a)$ we used the triangle inequality and $\bu^{(0)}=0$.

Due to the following fact
\begin{equation}
\|\bu^{(t+1)}-\wbx^{(t)}\|=\|\bu^{(t+1)}-\bu^{(0)}+\bu^{(0)}-\wbx^{(t)}\|\le\|\bu^{(t+1)}-\bu^{(0)}\|+\|\bu^{(0)}-\wbx^{(t)}\|\le 4.3\widehat{c}\mathcal{S}+\mathcal{S}/(2\kappa\log(\frac{d\kappa}{\delta})),
\end{equation}
we have $\|\bu^{(t+1)}-\wbx^{(t)}\|\le 5\widehat{c}\mathcal{S}$ since $\widehat{c}\ge2$. Therefore, we know that there exits $c^{(1)}_{\max}=\min\{c_{\max},c'_{\max}\}$ such that $\|\bu^{(t)}-\wbx^{(t)}\|\le 5\widehat{c}\mathcal{S},\forall t<T$ when $\eta\le c^{(1)}_{\max}/L_{\max}$, which completes the proof.
\end{proof}

\subsection{Proof of \leref{le.layer32}}
\begin{proof}
Let $\bu^{(0)}=0$ and define $\bv^{(t)}\bydef\bw^{(t)}-\bu^{(t)}$.
According
to the assumption of \leref{le.layer32}, we know that
$\bv^{(0)}=\upsilon[\eta L_{\max}\mathcal{S}/(\kappa\log(\frac{d\kappa}{\delta})\mathcal{P}_1)]\vec{\be}$ when
$\upsilon\in[\delta/(2\sqrt{d}),1]$.  First, we define an auxiliary function
\begin{equation}
\notag
h(\theta)\bydef\left[\begin{array}{c}\nabla_1 f(\bu^{(t)}_1+\theta\bv^{(t)}_1,\bu^{(t)}_2+\theta\bv^{(t)}_2)
\\
\nabla_2 f(\bu^{(t+1)}_1+\theta\bv^{(t+1)}_1,\bu^{(t)}_2+\theta\bv^{(t)}_2)\end{array}\right],
\end{equation}
then have
\begin{align}
\notag
& h(0)=\left[\begin{array}{c}\nabla_1 f(\bu^{(t)}_1,\bu^{(t)}_2)\\ \nabla_2 f(\bu^{(t+1)}_1,\bu^{(t)}_2)\end{array}\right],\quad
h(1)=\left[\begin{array}{c}\nabla_1 f(\bu^{(t)}_1+\bv^{(t)}_1,\bu^{(t)}_2+\bv^{(t)}_2)\\ \nabla_2 f(\bu^{(t+1)}_1+\bv^{(t+1)}_1,\bu^{(t)}_2+\bv^{(t)}_2)\end{array}\right],
\\\notag
&g(\theta)=\frac{dh(\theta)}{d\theta}=\underbrace{\left[\begin{array}{cc}\nabla^2_{11} f(\bu^{(t)}_1+\theta\bv^{(t)}_1,\bu^{(t)}_2+\theta\bv^{(t)}_2) & \nabla^2_{12} f(\bu^{(t)}_1+\theta\bv^{(t)}_1,\bu^{(t)}_2+\theta\bv^{(t)}_2) \\
\\\notag
0 & \nabla^2_{22} f(\bu^{(t+1)}_1+\theta\bv^{(t+1)}_1,\bu^{(t)}_2+\theta\bv^{(t)}_2)  \\
\end{array}\right]}_{\tbH^{(t)}_u(\theta)}\bv^{(t)}
\\
&\quad\quad\quad+\underbrace{\left[\begin{array}{cc}0 & 0 \\
\\\notag
\nabla^2_{21} f(\bu^{(t+1)}_1+\theta\bv^{(t+1)}_1,\bu^{(t)}_2+\theta\bv^{(t)}_2) & 0  \\
\end{array}\right]}_{\tbH^{(t)}_l(\theta)}\bv^{(t+1)},
\\\notag
& \left[\begin{array}{c}\nabla_1 f(\bw^{(t)}_1,\bw^{(t)}_2)\\  \nabla_2 f(\bw^{(t+1)}_1,\bw^{(t)}_2)\end{array}\right]=\int^1_0g(\theta)d\theta+\left[\begin{array}{c}\nabla_1 f(\bu^{(t)}_1,\bu^{(t)}_2)\\ \nabla_2 f(\bu^{(t+1)}_1,\bu^{(t)}_2)\end{array}\right].
\end{align}

Then, we consider sequence $\bw^{(t)}$, i.e.,
\begin{align}
&\bu^{(t+1)}+\bv^{(t+1)}=\bw^{(t+1)}=\bw^{(t)}-\eta\left[\begin{array}{c}\nabla_1 f(\bw^{(t)}_1,\bw^{(t)}_2)\\  \nabla_2 f(\bw^{(t+1)}_1,\bw^{(t)}_2)\end{array}\right]\label{eq.updateofw}
\\\notag
=&\bu^{(t)}+\bv^{(t)}-\eta\left[\begin{array}{c}\nabla_1 f(\bu^{(t)}_1+\bv^{(t)}_1,\bu^{(t)}_2+\bv^{(t)}_2)\\  \nabla_2 f(\bu^{(t+1)}_1+\bv^{(t+1)}_1,\bu^{(t)}_2+\bv^{(t)}_2)\end{array}\right]
\\
=&\bu^{(t)}+\bv^{(t)}-\eta\left[\begin{array}{c}\nabla_1 f(\bu^{(t)}_1,\bu^{(t)}_2)\\  \nabla_2 f(\bu^{(t+1)}_1,\bu^{(t)}_2)\end{array}\right]-\int^1_0g(\theta)d\theta
\\
\mathop{=}\limits^{(a)} &\bu^{(t)}+\bv^{(t)}-\eta\left[\begin{array}{c}\nabla_1 f(\bu^{(t)}_1,\bu^{(t)}_2)\\  \nabla_2 f(\bu^{(t+1)}_1,\bu^{(t)}_2)\end{array}\right]-\eta\tDelta^{(t)}_u \bv^{(t)} - \eta\bH_u\bv^{(t)}-\eta\tDelta^{(t)}_l \bv^{(t+1)} - \eta\bH_l\bv^{(t+1)}
\label{eq.dynofv}
\end{align}
where in $(a)$ we used the following definitions:
\begin{align}
\tDelta^{(t)}_u\bydef\int^1_0\tbH^{(t)}_u(\theta)d\theta-\bH_u,
\\
\tDelta^{(t)}_l\bydef\int^1_0\tbH^{(t)}_l(\theta)d\theta-\bH_l,
\end{align}
and
\begin{equation}
\bH_u\bydef\left[\begin{array}{cc}\nabla^2_{11} f(\wbx^{(t)}) &  \nabla^2_{12} f(\wbx^{(t)})
\\
0 &   \nabla^2_{22} f(\wbx^{(t)})\end{array}\right]\quad\bH_l\bydef\left[\begin{array}{cc}0 & 0
\\
\nabla^2_{21} f(\wbx^{(t)}) & 0 \end{array}\right].
\end{equation}
Obviously, $\bH=\bH_l+\bH_u$.

\paragraph{Dynamics of $\bv^{(t)}$:}

Since the first two terms at RHS of \eqref{eq.dynofv} combined with $\bu^{(t)}$ at LHS of \eqref{eq.dynofv} are exactly the same as \eqref{eq.bcdupdateo}. It can be observed that equation \eqref{eq.dynofv} gives the dynamic of
$\bv^{(t)}$, i.e.,
\begin{equation}\label{eq.rev}
\bv^{(t+1)}=\bv^{(t)}-\eta\tDelta^{(t)}_u \bv^{(t)} - \eta\bH_u\bv^{(t)}-\eta\tDelta^{(t)}_l \bv^{(t+1)} - \eta\bH_l\bv^{(t+1)}.
\end{equation}

Then, we can rewrite \eqref{eq.rev} in a matrix form as the following.
\begin{equation}
\underbrace{(\bI+\eta\bH_l)}_{\bydef\bM}\bv^{(t+1)} +\eta\tDelta^{(t)}_l \bv^{(t+1)}\mathop{=}\limits^{\eqref{eq.dynofv}}\underbrace{(\bI-\eta\bH_u) }_{\bydef\bT}\bv^{(t)}-\eta\tDelta^{(t)}_u \bv^{(t)}.\label{eq.iterav}
\end{equation}
It is worth noting that matrix $\bM$ is a lower triangular matrix where the diagonal entries are all 1s, so it is invertible.

Taking the inverse of $\bM$ on both sides of \eqref{eq.iterav}, we can obtain
\begin{equation}
\bv^{(t+1)}+\bM^{-1}\eta\tDelta^{(t)}_l\bv^{(t+1)}=\bM^{-1}\bT\widehat{\bv}^{(t)}-\bM^{-1}\eta\tDelta^{(t)}_u \bv^{(t)}.\label{eq.iterav2}
\end{equation}

Let $\mathbb{P}_{\texttt{left}}$ denote the projection operator that projects the vector onto the space spanned by the eigenvector of $\bM^{-1}\bT$ whose  eigenvalue is maximum. Taking the projection on both sides of \eqref{eq.iterav2}, we have
\begin{equation}
\mathbb{P}_{\texttt{left}}\widehat{\bv}^{(t+1)}+\mathbb{P}_{\texttt{left}}\bM^{-1}\eta\tDelta^{(t)}_l \bv^{(t+1)}=\mathbb{P}_{\texttt{left}}(\bM^{-1}\bT)\widehat{\bv}^{(t)}-\mathbb{P}_{\texttt{left}}\bM^{-1}\eta\tDelta^{(t)}_u\bv^{(t)}.\label{eq.keyre}
\end{equation}

From \leref{le.eiginver}, we know that the maximum eigenvalue of $\bM^{-1}\bT$ is greater than 1.

\paragraph{Relationship of the Norm of $\bv^{(t)}$ Projected in the Two Subspaces:}

Let $\phi^{(t)}$ denote the norm of $\bv^{(t)}$ projected onto the space spanned by the eigenvector of $\bM^{-1}\bT$ whose maximum eigenvalue is $1+\wde$
where $\wde\ge\eta\gamma/(1+L/L_{\max})$ due to \leref{le.eiginver}, and $\theta^{(t)}$ denote the norm of $\bv^{(t)}$ projected onto the remaining space. From \eqref{eq.keyre}, we can have
\begin{align}
\phi^{(t+1)}\mathop{\ge}\limits^{(a)}(1+\wde)\phi^{(t)}-\eta\|\bM^{-1}\|\|\tDelta^{(t)}_l\|\|\hbv^{(t+1)}\|-\eta\|\bM^{-1}\|\|\tDelta^{(t)}_u\|\|\bv^{(t)}\|,\label{eq.recurphi}
\\
\theta^{(t+1)}\le(1+\wde)\theta^{(t)}+\eta\|\bM^{-1}\|\|\tDelta^{(t)}_l\|\|\hbv^{(t+1)}\|+\eta\|\bM^{-1}\|\|\tDelta^{(t)}_u\|\|\bv^{(t)}\|.\label{eq.recurtheta}
\end{align}
where $(a)$ is true because we applied the triangle inequality since $\eta$ is sufficiently small. Also, since $\bM^{-1}=\bI-\eta\bH_l$, we have
\begin{align}
\notag
\|\bM^{-1}\|\le&1+\eta\|\bH_l\|
\\\notag
\mathop{=}\limits^{(a)}&1+\|\eta\bH\odot\bD-\eta\bH_d\|
\\\notag
\le&1+\eta\|\bH\odot\bD\|+\eta\|\bH_d\|
\\\notag
\mathop{\le}\limits^{(b)}&1+\eta (1+\frac{1}{\pi}+\frac{\log(d)}{\pi})\|\bH\|+\eta\|\bH_d\|
\\\notag
\mathop{\le}\limits^{(c)}&1+\eta\log(2d)\|\bH\|+\eta\|\bH_d\|
\\\notag
\mathop{\le}\limits^{(d)}&1+\eta L\log(2d)+\eta L_{\max}
\\
\le& 1+\frac{L}{L_{\max}}\log(2d)+1<2(1+\frac{L\log(2d)}{L_{\max}})\label{eq.minv}
\end{align}
where in $(a)$ $\odot$ denotes the Hadamard product and
\begin{equation}\notag
\bH_d\bydef\left[\begin{array}{cc}\nabla^2_{11}f(\wbx^{(t)}) & 0\\
0 & \nabla^2_{22}f(\wbx^{(t)}) \end{array}\right]\quad\bD=\left[\begin{array}{cccc}  1 & 0 & \cdots & 0
\\
1 & 1 & \cdots & 0
\\ \vdots & \vdots & \ddots & \vdots  \\ 1 & \cdots & 1 & 1\end{array}\right]\in\mathbb{R}^{d\times d}
\end{equation}
and inequality $(b)$ comes from the result on the spectral norm of the triangular truncation operator (please see [Theorem 1]\cite{angelos1992triangular}). In particular, by defining
\begin{equation}\notag
Y(\bD)\bydef\max\left\{\frac{\|\bH\odot\bD\|}{\|\bH\|},\bH\neq 0\right\},
\end{equation}
we have
\begin{equation}
\left|\frac{Y(\bD)}{\log(d)}-\frac{1}{\pi}\right|\le(1+\frac{1}{\pi})\frac{1}{\log(d)},
\end{equation}
$(c)$ is true for $d\ge3$, in $(d)$ we used the fact that $\|\bH\|\le L$ and $\|\bH_d\|\le L_{\max}$.

Since
$\|\bw^{(0)}-\wbx^{(t)}\|\le\|\bu^{(0)}-\wbx^{(t)}\|+\|\bv^{(0)}\|\le 2r$,
we can apply \leref{le.layer31}. Then, we know
$\|\bw^{(t)}-\wbx^{(t)}\|\le 5\widehat{c}\mathcal{S},\forall t< T$. According to the
assumptions of \leref{le.layer32}, we have $\|\bu^{(t)}-\wbx^{(t)}\|\le
5\widehat{c}\mathcal{S}$, and
\begin{equation}\label{eq.sizeofv}
\|\bv^{(t)}\|=\|\bw^{(t)}-\bu^{(t)}\|\le\|\bu^{(t)}-\wbx^{(t)}\|+\|\bw^{(t)}-\wbx^{(t)}\|\le 10\widehat{c}\mathcal{S}.
\end{equation}

From \eqref{eq.bdofdiffu}, we know that
\begin{equation}\notag
\|\bw^{(t+1)}-\bw^{(t)}\|\le\frac{4.3\eta\mathcal{G}}{\kappa}=
\frac{4.3\eta^3L^3_{\max}\frac{\gamma}{\rho}}{\kappa^2\log^3\frac{d\kappa}{\delta}\mathcal{P}^3_1\mathcal{P}_2}\le\mathcal{S},
\end{equation}
since $\mathcal{P}_1\ge2$ and we choose $\eta\le c_{\max}/L_{\max}$ and $c_{\max}=1/10$. Similarly, we also have $\|\bu^{(t+1)}-\bu^{(t)}\|\le\mathcal{S}$.

According to Lipsichiz continuity, we have the following bounds of $\|\bv^{(t+1)}\|$, $\|\tDelta^{(t)}_u\|$ and  $\|\tDelta^{(t)}_l\|$.
\begin{enumerate}
\item Relation between $\|\bv^{(t)}\|$ and $\|\bv^{(t+1)}\|$:
We also know that
\begin{align}
\notag
\|\bv^{(t+1)}\|^2=&\|\bw^{(t+1)}-\bu^{(t+1)}\|^2=\left\|\bw^{(t)}-\eta\left[\begin{array}{c}\nabla_1 f(\bw^{(t)}_1,\bw^{(t)}_2)\\  \nabla_2 f(\bw^{(t+1)}_1,\bw^{(t)}_2)\end{array}\right] -\left(\bu^{(t)}-\eta\left[\begin{array}{c}\nabla_1 f(\bu^{(t)}_1,\bu^{(t)}_2)\\ \nabla_2 f(\bu^{(t+1)}_1,\bu^{(t)}_2)\end{array}\right]\right)\right\|^2
\\\notag
\le&2\|\bv^{(t)}\|^2+4\eta^2\left\|\left[\begin{array}{c}\nabla_1 f(\bw^{(t)}_1,\bw^{(t)}_2)\\ \nabla_2 f(\bw^{(t+1)}_1,\bw^{(t)}_2)\end{array}\right]-\left[\begin{array}{c}\nabla_1 f(\bu^{(t)}_1,\bw^{(t)}_2)\\  \nabla_2 f(\bu^{(t+1)}_1,\bw^{(t)}_2)\end{array}\right]\right\|^2
\\\notag
\quad&+4\eta^2\left\|\left[\begin{array}{c}\nabla_1 f(\bu^{(t)}_1,\bw^{(t)}_2)\\  \nabla_2 f(\bu^{(t+1)}_1,\bw^{(t)}_2)\end{array}\right]-\left[\begin{array}{c}\nabla_1 f(\bu^{(t)}_1,\bu^{(t)}_2)\\ \nabla_2 f(\bu^{(t+1)}_1,\bu^{(t)}_2)\end{array}\right]\right\|^2
\\
\mathop{\le}\limits^{(a)}& 2\|\bv^{(t)}\|^2+4\eta^2L^2_{\max}(\|\bv^{(t+1)}_1\|^2+\|\bv^{(t)}_1\|^2)+8\eta^2L^2_{\max}\|\bv^{(t)}_2\|^2\label{eq.revv}
\end{align}
where $(a)$ is true due to Lipschitz continuity.

We can express \eqref{eq.revv} as
\begin{equation}\notag
(1-4\eta^2L^2_{\max})\|\bv^{(t+1)}\|\le(2+8\eta^2L^2_{\max})\|\bv^{(t)}\|^2
\end{equation}
which implies
\begin{equation}
\|\bv^{(t+1)}\|\le\sqrt{\frac{2+\frac{8}{100}}{1-\frac{4}{100}}}\|\bv^{(t)}\|<\sqrt{2.2}\|\bv^{(t)}\|<1.5\|\bv^{(t)}\|,\label{eq.relavv}
\end{equation}
where we choose $\eta\le c_{\max}/L_{\max}$ and $c_{\max}=1/10$. 

\item Bounds of $\|\tDelta^{(t)}_u\|$ and $\|\tDelta^{(t)}_l\|$:

According to $\rho$-Hessian Lipschitz continuity and \leref{le.layer001bl}, we have the size of $\tDelta^{(t)}_u$ as the following.
\begin{align}
\notag
\|(\tDelta^{(t)}_u)\|\le&\int^1_0\|\tbH^{(t)}_u(\theta)-\bH_u\|d\theta
\\
\mathop{\le}\limits^{\eqref{eq.deltatbd}}&\int^1_0\rho\left(\|\bu^{(t)}+\theta\bv^{(t)}-\wbx^{(t)}\|+\left\|\left[\begin{array}{c}\bu^{(t+1)}_1+\theta\bv^{(t+1)}_1 \\ \bu^{(t)}_2+\theta\bv^{(t)}_2\end{array}\right]-\wbx^{(t)}\right\|\right)d\theta
\\\notag
\mathop{\le}\limits^{(a)}&\int^1_0\rho\left(2\|\bu^{(t)}+\theta\bv^{(t)}-\wbx^{(t)}\|+\|\bu^{(t+1)}+\theta\bv^{(t+1)}-\wbx^{(t)}\|\right)d\theta
\\\notag
\le&\rho(\|\bu^{(t+1)}-\wbx^{(t)}\|+2\|\bu^{(t)}-\wbx^{(t)}\|)+\rho\int^1_0\theta(\|\bv^{(t+1)}\|+\|\bv^{(t)}\|)d\theta
\\\notag
\le&\rho\left(\|\bu^{(t+1)}-\bu^{(t)}\|+\|\bu^{(t)}-\wbx^{(t)}\|+2\|\bu^{(t)}-\wbx^{(t)}\|)+0.5\|\bv^{(t+1)}\|+0.5\|\bv^{(t)}\|\right)
\\\notag
\mathop{\le}\limits^{\eqref{eq.relavv}} &\rho\left(\|\bu^{(t+1)}-\bu^{(t)}\|+3\|\bu^{(t)}-\wbx^{(t)}\|+1.25\|\bv^{(t)}\|\right)
\\\notag
\le & \rho (1+27.5\widehat{c})\mathcal{S}
\end{align}
where $(a)$ is true because
\begin{align}
\notag
\left\|\left[\begin{array}{c}\bu^{(t+1)}_1+\theta\bv^{(t+1)}_1 \\ \bu^{(t)}_2+\theta\bv^{(t)}_2\end{array}\right]-\wbx^{(t)}\right\|\le&\left\|\bI_1\left(\bu^{(t+1)}+\theta\bv^{(t+1)}-\wbx^{(t)}\right)\right\|+\left\|\bI_2\left(\bu^{(t)}+\theta\bv^{(t)}-\wbx^{(t)}\right)\right\|
\\
\mathop{\le}\limits^{\eqref{eq.idef}}&\|\bu^{(t+1)}+\theta\bv^{(t+1)}-\wbx^{(t)}\|+\|\bu^{(t)}+\theta\bv^{(t)}-\wbx^{(t)}\|.
\end{align}
%

Applying \leref{le.layer001bl}, we can also get the upper bound of $\|\tDelta^{(t)}_l\|$, i.e.,
\begin{align}
\notag
\|(\tDelta^{(t)}_l)\|\le&\int^1_0\|\tbH^{(t)}_l(\theta)-\bH_l\|d\theta
\\
\mathop{\le}\limits^{\eqref{eq.deltatbd2}}&\int^1_0\rho\left\|\left[\begin{array}{c}\bu^{(t+1)}_1+\theta\bv^{(t+1)}_1 \\ \bu^{(t)}_2+\theta\bv^{(t)}_2\end{array}\right]-\wbx^{(t)}\right\|d\theta
\\\notag
\le&\int^1_0\rho(\|\bu^{(t)}+\theta\bv^{(t)}-\wbx^{(t)}\|+\|\bu^{(t+1)}+\theta\bv^{(t+1)}-\wbx^{(t)})\|d\theta
\\\notag
\le&\rho(\|\bu^{(t+1)}-\wbx^{(t)}\|+\|\bu^{(t)}-\wbx^{(t)}\|)+\rho\int^1_0\theta(\|\bv^{(t+1)}\|+\|\bv^{(t)}\|)d\theta
\\\notag
\mathop{\le}\limits^{\eqref{eq.relavv}} &\rho\left(\|\bu^{(t+1)}-\bu^{(t)}\|+2\|\bu^{(t)}-\wbx^{(t)}\|+1.25\|\bv^{(t)}\|\right)
\\\notag
\le&\rho(1+22.5\widehat{c})\mathcal{S}.
\end{align}
\end{enumerate}

With the upper bounds of $\|\bv^{(t+1)}\|$, $\|\tDelta^{(t)}_u\|$, $\|\tDelta^{(t)}_l\|$  and relation between $\|\bv^{(t+1)}\|$ and $\|\bv^{(t)}\|$, we can further simply \eqref{eq.recurphi} and \eqref{eq.recurtheta} as follows,
\begin{align}
\notag
\phi^{(t+1)}\mathop{\ge}\limits^{\eqref{eq.recurphi}}(1+\wde)\phi^{(t)}-\eta(1.5\|\tDelta^{(t)}_l\|+\|\tDelta^{(t)}_u\|)\|\bM^{-1}\|\|\bv^{(t)}\|
\\\notag
\theta^{(t+1)}\mathop{\le}\limits^{\eqref{eq.recurtheta}}(1+\wde)\theta^{(t)}+\eta(1.5\|\tDelta^{(t)}_l\|+\|\tDelta^{(t)}_u\|)\|\bM^{-1}\|\|\bv^{(t)}\|
\end{align}
and further we have
\begin{align}
\notag
\phi^{(t+1)}\ge(1+\wde)\phi^{(t)}-\eta(1.5\|\tDelta^{(t)}_l\|+\|\tDelta^{(t)}_u\|)\|\bM^{-1}\|\sqrt{(\phi^{(t)})^2+(\theta^{(t)})^2},
\\\notag
\theta^{(t+1)}\le(1+\wde)\theta^{(t)}+\eta(1.5\|\tDelta^{(t)}_l\|+\|\tDelta^{(t)}_u\|)\|\bM^{-1}\|\sqrt{(\phi^{(t)})^2+(\theta^{(t)})^2},
\end{align}
since $\|\bv^{(t)}\|=\sqrt{(\phi^{(t)})^2+(\theta^{(t)})^2}$.

Consequently, we can arrive at
\begin{align}
\phi^{(t+1)}\ge(1+\wde)\phi^{(t)}-\mu\sqrt{(\phi^{(t)})^2+(\theta^{(t)})^2},\label{eq.phii}
\\
\theta^{(t+1)}\le(1+\wde)\theta^{(t)}+\mu\sqrt{(\phi^{(t)})^2+(\theta^{(t)})^2},\label{eq.thei}
\end{align}
where $\mu$ is the upper bound of $\eta(1.5\|\tDelta^{(t)}_l\|+\|\tDelta^{(t)}_u\|)\|\bM^{-1}\|$ and can be obtained by
\begin{equation}
\mu\bydef\eta\rho\mathcal{S}\mathcal{P}_2(2.5+62\widehat{c})\label{eq.defofmu}.
\end{equation}

\paragraph{Quantifying the Norm of $\bv^{(t)}$ Projected at Different Subspaces:}
Then, we will use mathematical induction to prove
\begin{equation}\label{eq.phinorm}
\theta^{(t)}\le 4\mu t\phi^{(t)}.
\end{equation}
It is true when $t=0$ since $\|\theta^{(0)}\|\mathop{=}\limits^{\eqref{eq.inicond}}0$.


Assuming that equation $\eqref{eq.phinorm}$ is
true at the $t$th iteration, we need to prove
\begin{equation}
\theta^{(t+1)}\le 4\mu (t+1)\phi^{(t+1)}.\label{eq.goa}
\end{equation}

Applying \eqref{eq.phii} into RHS of \eqref{eq.goa}, we have
\begin{equation}\label{eq.rephi}
4\mu(t+1)\phi^{(t+1)}\ge4\mu(t+1)\left((1+\wde)\phi^{(t)}-\mu\sqrt{(\phi^{(t)})^2+(\theta^{(t)})^2}\right)
\end{equation}
and substituting \eqref{eq.thei} into LHS of \eqref{eq.goa}, we have
\begin{equation}\label{eq.retheta}
\theta^{(t+1)}\le (1+\wde)(4\mu t\phi^{(t)})+\mu\sqrt{(\phi^{(t)})^2+(\theta^{(t)})^2}.
\end{equation}

Then, our goal is to prove RHS of \eqref{eq.rephi} is greater
than RHS of \eqref{eq.retheta}. After some manipulations, it is sufficient to show
\begin{equation}
\left(1+4\mu(t+1)\right)\left(\sqrt{(\phi^{(t)})^2+(\theta^{(t)})^2}\right)\le4\phi^{(t)}.
\end{equation}

In the following, we will show that the above relation is true.
\paragraph{First step}:
We know that
\begin{equation}\label{eq.bdofmu}
4\mu(t+1)\le4\mu T\mathop{\le}\limits^{\eqref{eq.defofmu}}4\eta\rho\mathcal{S}\mathcal{P}_2(2.5+62\widehat{c})\widehat{c}\mathcal{T}\mathop{\le}\limits^{\eqref{eq.sandt}\eqref{eq.defofmu}}\frac{4\widehat{c}\eta^2 L^2_{\max}(2.5+62\widehat{c})}{\kappa\log(\frac{d\kappa}{\delta})\mathcal{P}_1}\mathop{\le}\limits^{(a)}1
\end{equation}
where $(a)$ is true because $\mathcal{P}_1\ge2$ and we choose $c'_{\max}=1/(2\widehat{c}(2.5+62\widehat{c}))$ and
$\eta\le c'_{\max}/L_{\max}$.

\paragraph{Second step}:
Also, we know that
\begin{equation}
4\phi^{(t)}\ge2\sqrt{2(\phi^{(t)})^2}\mathop{\ge}\limits^{\eqref{eq.phinorm},\eqref{eq.bdofmu}}(1+4\mu(t+1))\sqrt{(\phi^{(t)})^2+(\theta^{(t)})^2}.
\end{equation}

With the above two steps, we have $\theta^{(t+1)}\le 4\mu(t+1)\phi^{(t+1)}$, which completes the induction.

\paragraph{Recursion of $\phi^{(t)}$}:Using \eqref{eq.phinorm}, we have $\theta^{(t)}\mathop{\le}\limits^{\eqref{eq.phinorm}} 4\mu t\phi^{(t)}\mathop{\le}\limits^{\eqref{eq.bdofmu}} \phi^{(t)}$, which implies
\begin{align}
\notag
\phi^{(t+1)}\mathop{\ge}\limits^{\eqref{eq.phii}}&(1+\wde)\phi^{(t)}-\mu\sqrt{(\phi^{(t)})^2+(\theta^{(t)})^2}
\\\notag
\mathop{\ge}\limits^{(a)}&(1+\frac{\gamma\eta}{1+L/L_{\max}})\phi^{(t)}-\mu\sqrt{(\phi^{(t)})^2+(\theta^{(t)})^2}
\\
\mathop{\ge}\limits^{(b)}&(1+\frac{1}{1+L/L_{\max}}\frac{\gamma\eta}{2})\phi^{(t)}\label{eq.relation}
\end{align}
where in $(a)$ we used \leref{le.eiginver}, and $(b)$ is true because
\begin{align}
\notag
\mu=&\eta\rho\mathcal{S}\mathcal{P}_2(2.5+62\widehat{c})
\\\notag
\le&\frac{\gamma\eta}{1+L/L_{\max}}\frac{\eta^2 L^2_{\max}(2.5+62\widehat{c})}{\log^{2}(\frac{d\kappa}{\delta})\mathcal{P}_1}
\\\notag
\mathop{\le}\limits^{(a)}&\frac{1}{1+L/L_{\max}}\frac{\gamma\eta}{2\sqrt{2}}
\end{align}
where in $(a)$ we choose $c''_{\max}= 1/(2\sqrt{2}(2.5+62\widehat{c}))$ and $\eta\le c''_{\max}/L_{\max}$.

\paragraph{Quantifying Escaping Time:}

From  \eqref{eq.sizeofv}, we have
\begin{align}
\notag
10\mathcal{S}\widehat{c}\ge&\|\bv^{(t)}\|\ge\phi^{(t)}
\\\notag
\mathop{\ge}\limits^{\eqref{eq.relation}}&(1+\frac{\gamma\eta}{2(1+L/L_{\max})})^{t}\phi^{(0)}
\\\notag
\mathop{\ge}\limits^{(a)}&(1+\frac{\gamma\eta}{2(1+L/L_{\max})})^{t}\frac{\delta}{2\sqrt{d}}\frac{\eta L_{\max}\mathcal{S}}{\kappa}\log^{-1}(\frac{d\kappa}{\delta})\mathcal{P}^{-1}_1
\\
\mathop{\ge}\limits^{(b)}&(1+\frac{\gamma\eta}{2(1+L/L_{\max})})^{t}\frac{\delta}{2\sqrt{d}}\frac{c\mathcal{S}}{\kappa}\log^{-1}(\frac{d\kappa}{\delta})\mathcal{P}^{-1}_1 \quad\forall t<T \label{eq.upoft}
\end{align}
where in $(a)$ we use condition $\upsilon\in[\delta/(2\sqrt{d}),1]$, in $(b)$ we used $\eta =c/L_{\max}$.

Since \eqref{eq.upoft} is true for all $t<T$,  we can have
\begin{align}
\notag
T-1\le&\frac{\log(20\frac{\widehat{c}}{c}(\frac{\kappa\sqrt{d}}{\delta})\log(\frac{d\kappa}{\delta})\mathcal{P}_1)}{\log(1+\frac{\eta\gamma}{2(1+L/L_{\max})})}
\\\notag
\mathop{<}\limits^{(a)}&\frac{4(1+L/L_{\max})\log(20(\frac{\sqrt{d}\kappa}{\delta})\frac{\widehat{c}}{c}\log(\frac{d\kappa}{\delta})\mathcal{P}_1)}{\eta\gamma}
\\\notag
\mathop{<}\limits^{(b)}&\frac{4(1+L/L_{\max})\log(20(\frac{d\kappa }{\delta})^2\frac{\widehat{c}}{c}\mathcal{P}_1)}{\eta\gamma}
\\
\mathop{<}\limits^{(c)}&4(2+\log(20\frac{\widehat{c}}{c}))\mathcal{T}\label{eq.bdt}
\end{align}
where $(a)$ comes from inequality $\log(1+x)>x/2$ when $x<1$,  in $(b)$ we used relation $\log(x)<x, x>0$,   and $(c)$ is true because $\delta\in(0,\frac{d\kappa}{e}]$ and $\log(d\kappa/\delta)>1$ and $\mathcal{P}_1>1$ we have
\begin{equation}\notag
\log(\frac{d\kappa}{\delta}\mathcal{P}_1)\le\log(\frac{d\kappa}{\delta})+\log(1+\frac{L}{L_{\max}})\le\log(\frac{d\kappa}{\delta})+\frac{L}{L_{\max}}\le\log(\frac{d\kappa}{\delta})\mathcal{P}_1.
\end{equation}

From \eqref{eq.bdt}, we know that
\begin{equation}
T<4(2+\log(20\frac{\widehat{c}}{c}))\mathcal{T}+1\mathop{<}\limits^{(a)}4(2\frac{1}{4}+\log(20\frac{\widehat{c}}{c})\mathcal{T}
\end{equation}
where $(a)$ is true due to the fact that $\eta L_{\max}\ge1$, $\log(d\kappa/\delta)>1$ and $\mathcal{P}_1>1$ so we know $\mathcal{T}\ge1$.

When
\begin{equation}
4(2.25+\log(20\frac{\widehat{c}}{c}))\le\widehat{c},\label{eq.chat}
\end{equation}
we will have $T<\widehat{c}\mathcal{T}$ where $c^{(2)}_{\max}\bydef\min\{c_{\max},c'_{\max},c''_{\max}\}$.

Since $\widehat{c}\ge2$, we have $c_{\max}=\min\{c^{(1)}_{\max},c^{(2)}_{\max}\}\le 1/(5\widehat{c})^3$. Also, we know that $c\le c_{\max}$. Combining with \eqref{eq.chat}, we need
\begin{equation}
\frac{\widehat{c}}{2^{\frac{\widehat{c}}{4}-2.25-\log(20)}}\le c\le \frac{1}{(5\widehat{c})^3},
\end{equation}
meaning that
\begin{equation}
125(2^{2.25+\log(20)}\widehat{c}^4)\le 2^{\frac{\widehat{c}}{4}}.\label{eq.bdhatc}
\end{equation}
It can be observed that LHS of \eqref{eq.bdhatc} is a polynomial with respect to $\widehat{c}$ and RHS of \eqref{eq.bdhatc} is a exponential function in terms of $\widehat{c}$, implying there exists a universal $\widehat{c}$ such that \eqref{eq.bdhatc} holds. The proof is complete.
\end{proof}

\subsection{Proof of \leref{le.layer21}}
\begin{proof}

The proof of \leref{le.layer21} is similar as the one of proving convergence of PGD shown in \citep[Lemma 14,15]{jin2017jordan}. Considering the completeness of the whole proof in this paper, here we give the following proof of this lemma in details.

First, after the random perturbation, the objective function value in the worst case is increased at most by
\begin{align}
\notag
f(\bu^{(0)})-f(\wbx^{(t)})\le& \sum^2_{k=1}\nabla_kf(\tbh^{(t)}_{-k},\wbx^{(t)}_k)^{\T}\xi_k+\frac{L_k}{2}\|\xi_k\|^2
\\\notag
\le&\sum^2_{k=1}\|\nabla_kf(\tbh^{(t)}_{-k},\wbx^{(t)}_k)\|\|\xi_k\| +\frac{L_{\max}}{2}\|\xi\|^2
\\\notag
\mathop{\le}\limits^{(a)}&\|\xi\|\sqrt{\sum^2_{k=1}2\|\nabla_kf(\tbh^{(t)}_{-k},\wbx^{(t)}_k)\|^2}+\frac{L_{\max}}{2}\|\xi\|^2
\\
\mathop{\le}\limits^{(b)}& \frac{\mathcal{G}}{\kappa}\frac{\eta L_{\max}\mathcal{S}}{\kappa\log(\frac{d\kappa}{\delta})\mathcal{P}_1}+\frac{L_{\max}}{2}(\frac{\eta L_{\max}\mathcal{S}}{\kappa\log(\frac{d\kappa}{\delta})\mathcal{P}_1})^2\le\frac{3}{2}\mathcal{F}\label{eq.sufin}
\end{align}
where $\bu^{(0)}$ is a vector that follows uniform distribution within the ball $\mathbb{B}^{(d)}_{\wbx^{(t)}}(r)$, $\mathbb{B}^{(d)}_{\wbx^{(t)}}$ denotes the $d$-dimensional ball centered at $\wbx^{(t)}$ with radius $r$, $\xi_k$ represents the $k$th block of the vector which is the difference between random generated vector $\bu^{(0)}$ and $\wbx^{(t)}$, and $(a)$ is true because $\xi\bydef[\xi_1,\ldots,\xi_K]$, $\|\xi_k\|\le\|\xi\|,\forall k$, and in $(b)$ we used $\kappa>1$, $\log(d\kappa/\delta)>1$ and \conref{cond}.

Second, under \asref{as1}, let $\widetilde{\bx}^{(t)}$ satisfy conditions \conref{cond}, and two
PA-GD iterates $\{\bu^{(t)}\}$ $\{\bw^{(t)}\}$ satisfy the conditions as in
\leref{le.layer32}. Selecting $c_{\max}=\min\{c^{(1)}_{\max},c^{(2)}_{\max}\}$, so we have that $\eta\le c_{\max}/L_{\max}$ is small enough such that \leref{le.layer31} and \leref{le.layer32} can both hold.

Let $T^*\bydef\widehat{c}\mathcal{T}$ and $T'\bydef\inf_t\{t|\widehat{f}_{\bu^{(0)}}(\bu^{(t)})-f(\bu^{(0)})\le-3\mathcal{F}\}$. Then, we have the following two cases to analyze the decrease of the objective value after $T$ iterations with the random perturbation.
\begin{enumerate}
\item Case $T'\le T^*$:
\begin{align}
\notag
f(\bu^{(T')})-f(\bu^{(0)})\le&\nabla f(\bu^{(0)})^{\T}(\bu^{(T')}-\bu^{(0)})+\frac{1}{2}(\bu^{(T')}-\bu^{(0)})^{\T}\nabla^2f(\bu^{(0)})(\bu^{(T')}-\bu^{(0)})+\frac{\rho}{6}\|\bu^{(T')}-\bu^{(0)}\|^3
\\\notag
\le&\widehat{f}_{\bu^{(0)}}(\bu^{(t)})-f(\bu^{(0)})+\frac{\rho}{2}\|\bu^{(0)}-\wbx^{(t)}\|\|\bu^{(T')}-\bu^{(0)}\|^2+\frac{\rho}{6}\|\bu^{(T')}-\bu^{(0)}\|^3
\\
\mathop{\le}\limits^{\eqref{eq.bdrho}-\eqref{eq.fdes}}& -3\mathcal{F}+0.5\rho\mathcal{S}^3\mathop{\le}\limits^{\eqref{eq.sandf}}-2.5\mathcal{F}\label{eq.ftf}.
\end{align}
Based on \leref{le.descent}, we know that AGD is always decreasing the objective function. For any $T\ge\mathcal{T}/c_{\max}\ge\widehat{c}\mathcal{T}=T^*\ge T'$, we have
\begin{equation}
\notag
f(\bu^{(T)})-f(\bu^{(0)})\le f(\bu^{(T^*)})-f(\bu^{(0)})\le f(\bu^{(T')})-f(\bu^{(0)})\le -2.5\mathcal{F}
\end{equation}
where $c_{\max}=\min\{1,1/\widehat{c}\}$.

\item Case $T'> T^*$: Applying \leref{le.layer31}, we know that $\|\bu^{(t)}-\bu^{(0)}\|\le 5\widehat{c}\mathcal{S}$ for $t\le T^*$. Define $T''=\inf_t\{t|\widehat{f}_{\bw^{(0)}}(\bw^{(t)})-f(\bw^{(0)})\le -3\mathcal{F}\}$. Then, after applying \leref{le.layer32}, we know $T''\le T^*$. Similar as \eqref{eq.ftf}, for $T\ge 1/c_{\max}\mathcal{T}$, we also have $f(\bw^{(T)})-f(\bw^{(0)})\le f(\bw^{T^*})-f(\bw^{(0)})\le f(\bw^{T''})-f(\bw^{(0)})\le -2.5\mathcal{F}$.
\end{enumerate}

Combining the above two cases, we have
\begin{equation}
\min\{f(\bu^{(T)})-f(\bu^{(0)}),f(\bw^{(T)})-f(\bw^{(0)})\}\le -2.5\mathcal{F},\label{eq.suffdes}
\end{equation}
meaning that at least one of the sequences can give a sufficient decrease of the objective function if the initial points of the two sequences are separated apart with each other far enough along direction $\vec{\be}$.

Therefore, we can conclude that if $\bu^{(0)}\in\mathcal{X}_{\textsf{stuck}}$, then $(\bu^{(0)}\pm\upsilon r\vec{\be})\notin\mathcal{X}_{\textsf{stuck}}$ where $\upsilon\in[\frac{\delta}{2\sqrt{d}},1]$.

Finally, we give the upper
bound of the volume of $\mathcal{X}_{\textsf{stuck}}$,
\begin{align}
\notag
\textsf{Vol}(\mathcal{X}_{\textsf{stuck}})=&\int_{\mathbb{B}^{(d)}_{\wbx^{(t)}}}d\bu I_{\mathcal{X}_{\textsf{stuck}}}(\bu)=\int_{\mathbb{B}^{(d-1)}_{\wbx^{(t)}}}du_{-1}\int^{\wx^{(t)}_1+\sqrt{r^2-\|\wbx^{(t)}_{-1}-u_{-1}\|^2}}_{\wx^{(t)}_1-\sqrt{r^2-\|\wbx^{(t)}_{-1}-u_{-1}\|^2}}du_1 I_{\mathcal{X}_{\textsf{stuck}}}(\bu)
\\\notag
\le&\int_{\mathbb{B}^{(d-1)}_{\wbx^{(t)}}}du_{-1}\left(2\frac{\delta}{2\sqrt{d}r}\right)=\textsf{Vol}(\mathbb{B}^{(d-1)}_{\wbx^{(t)}}(r))\frac{r\delta}{\sqrt{d}}
\end{align}
where $I_{\textsf{stuck}}(\bu)$ is an indicator function showing that $\bu$ belongs to set $\mathcal{X}_{\textsf{stuck}}$, and $u_1$ represents the component of vector $\bu$ along $\vec{\be}$ direction, and $\bu_{-1}$ is the remaining $d-1$ dimensional vector.

Then, the ratio of $\textsf{Vol}(\mathcal{X}_{\textsf{stuck}})$ over the whole volume of the perturbation ball can be upper bounded by
\begin{equation}\notag
\frac{\textsf{Vol}(\mathcal{X}_{\textsf{stuck}})}{\textsf{Vol}(\mathbb{B}^{(d)}_{\wbx^{(t)}}(r))}\le\frac{\frac{r\delta}{\sqrt{d}}
\textsf{Vol}(\mathbb{B}^{(d-1)}_{\wbx^{(t)}}(r))}{\textsf{Vol}(\mathbb{B}^{(d)}_{\wbx^{(t)}}(r))}
=\frac{\delta}{\sqrt{d\pi}}\frac{\Gamma(\frac{d}{2}+1)}{\Gamma(\frac{d}{2}+1)}\le\frac{\delta}{\sqrt{d\pi}}\sqrt{\frac{d}{2}+\frac{1}{2}}\le\delta
\end{equation}
where $\Gamma(\cdot)$ denotes the Gamma function, and inequality is true due to the fact that $\Gamma(x+1)/\Gamma(x+1/2)<\sqrt{x+1/2}$ when $x\ge0$.

Combining \eqref{eq.sufin} and \eqref{eq.suffdes}, we can show that
\begin{equation}
f(\bx^{(T)})-f(\wbx^{(t)})=f(\bx^{(T)})-f(\bu^{(0)})+f(\bu^{(0)})-f(\wbx^{(t)})\le-2.5\mathcal{F}+1.5\mathcal{F}\le-\mathcal{F}
\end{equation}
with at least probability $1-\delta$.
\end{proof}

\clearpage
\newpage

\section{Proof of PA-PP}

First, we need to introduce some constants defined as follows,

\begin{align}
\notag
\mathcal{F}\bydef& \eta^5 L^5_{\max}\frac{\gamma^{3}}{\kappa^3\rho^2}\log^{-6}\left(\frac{d\kappa}{\delta}\right)\mathcal{P}^{-2},
\quad\quad
\mathcal{G}\bydef\eta^2 L^2_{\max}\frac{\gamma^2}{\rho}\log^{-3}\left(\frac{d\kappa}{\delta}\right)\mathcal{P}^{-1},
\\\notag
\mathcal{S}\bydef&\eta^2 L^2_{\max}\frac{\gamma}{\kappa\rho}\log^{-2}\left(\frac{d\kappa}{\delta}\right)\mathcal{P}^{-1},
\quad\quad
\mathcal{T}\bydef\frac{\log\left(\frac{d\kappa}{\delta}\right)}{\eta\gamma}
\end{align}
where $\eta=1/\nu$. In order to keep the completeness of the proof, the certain relations of these quantities are listed as follows, which are useful of simplifying the expressions in the proofs.
\begin{subequations}
\begin{align}
\sqrt{\mathcal{F}}=&\frac{\sqrt{\eta}\mathcal{G}}{\kappa},\label{eq.bfandg}
\\
\frac{\eta\mathcal{G}\mathcal{T}}{\kappa}=&\mathcal{S},\label{eq.bgandt}
\\
\rho\mathcal{S}^3=&\frac{\eta L_{\max}\mathcal{F}}{\mathcal{P}},\label{eq.bsandf}
\\
\eta\rho\mathcal{S}\mathcal{T}=&\frac{\eta^2 L^2_{\max}}{\kappa\log(\frac{d\kappa}{\delta})\mathcal{P}},\label{eq.bsandt}
\\
\eta\rho\mathcal{S}=&\eta L_{\max}\frac{\eta^2\gamma^2}{\log^2(\frac{d\kappa}{\delta})\mathcal{P}}.\label{eq.seg}
\end{align}
\end{subequations}
We also consider saddle point $\widetilde{\bx}^{(t)}$ that satisfies the following condition.
\begin{condition}\label{cond2}
An $\epsilon$-second order stationary point $\widetilde{\bx}^{(t)}$ satisfies the following conditions:
\begin{equation}
\|\bx^{(t+1)}-\bx^{(t)}\|\le g_{\textsf{th}}/\nu\quad\textrm{and}\quad\lambda_{\min}(\nabla^2 f(\wbx^{(t)}))\le -\gamma\label{eq.bgcond}
\end{equation}
where $g_{\textsf{th}}=\frac{\mathcal{G}}{2\kappa}$.
\end{condition}
Then, we have the following preliminary lemmas.

\begin{lemma}\label{le.blayer001grad}
If function $f(\cdot)$ is $L$-smooth with Lipschitz constant, the we have
\begin{equation}
\|\nabla f(\bx^{(t)})\|^2\le4\nu\|\bx^{(t+1)}-\bx^{(t)}\|^2
\end{equation}
where sequence $\bx^{(t)}_k, k=1,2$ is generated by the APP algorithm.
\end{lemma}

\begin{lemma}\label{le.blayer001bl}
Under \asref{as1}, we have block-wise Lipschitz continuity as the follows:
\begin{equation}
\left\|\left[\begin{array}{cc}\nabla^2_{11}f(\bx) & \boldsymbol{0}
\\
\nabla^2_{21}f(\by) & \nabla^2_{22}f(\by)\end{array}\right]-\left[\begin{array}{cc}\nabla^2_{11}f(\bz) & \boldsymbol{0}\\
\nabla^2_{21}f(\bz) & \nabla^2_{22}f(\bz)\end{array}\right]\right\|
\le\rho\left(\|\bx-\bz\|+\|\by-\bz\|\right),\forall\bx,\by,\bz\label{eq.bdeltatbd}
\end{equation}
and
\begin{equation}
\left\|\left[\begin{array}{cc} \boldsymbol{0} & \nabla^2_{21}f(\bx)
\\
\boldsymbol{0} & \boldsymbol{0}\end{array}\right]-\left[\begin{array}{cc} \boldsymbol{0} & \nabla^2_{12}f(\by) \\
\boldsymbol{0} & \boldsymbol{0}\end{array}\right]\right\|
\\
\le\rho\|\bx-\by\|,\forall\bx,\by.\label{eq.bdeltatbd2}
\end{equation}
\end{lemma}

Second, we can have the descent lemma as the following

\begin{lemma}\label{le.bdescent}
Under \asref{as1}, for the APP algorithm with penalizer $\nu\ge 3L_{\max}$, we have
\begin{equation}
f(\bx^{(t+1)})\le f(\bx^{(t)})-\frac{\nu}{2}\|\bx^{(t+1)}-\bx^{(t)}\|^2.\notag
\end{equation}
\end{lemma}

Third, we need to characterize the convergence behaviour of PA-PP when $\|\bx^{(t+1)}-\bx^{(t)}\|$ is small. In this case, we need three steps to arrive the final results.
\paragraph{Step 1}: Quantify upper bound of the distance between generic iterate $\bu^{(t)}$ and saddle point $\wbx^{(t)}$.
\begin{lemma}\label{le.blayer31}
Under \asref{as1}, consider saddle point $\widetilde{\bx}^{(t)}$ that satisfies \conref{cond2}.
For any constant $\widehat{c}\ge 2$, $\delta\in(0,\frac{d\kappa}{e}]$, when
initial point $\bu^{(0)}$ satisfies
\begin{equation}
\|\bu^{(0)}-\widetilde{\bx}^{(t)}\|\le 2r,
\end{equation}
then, with the definition of
\begin{equation}
r\bydef \frac{\frac{L_{\max}}{\nu}\mathcal{S}}{\kappa\log(\frac{d\kappa}{\delta})\mathcal{P}_1}\quad \textrm{and}\quad T\bydef\min\{\inf_t\{t|\widehat{f}_{\bu^{(0)}}(\bu^{(t)})-f(\bu^{(0)})\le-3\mathcal{F}\},\widehat{c}\mathcal{T}\},\label{eq.bdefoft}
\end{equation}
there exits constants $c^{(1)}_{\max}, \widehat{c}$ such that for any $\nu\ge
L_{\max}/c^{(1)}_{\max}$, the iterates generated by PA-PP satisfy $\|\bu^{(t)}-\widetilde{\bx}^{(t)}\|\le5\widehat{c}\mathcal{S}, \forall t<T$.
\end{lemma}

\paragraph{Step 2}: Quantify the escaping time of iterates near a strict saddle point.
\begin{lemma}\label{le.blayer32}
Under \asref{as1}, consider saddle point $\widetilde{\bx}^{(t)}$ that satisfies satisfies \conref{cond2}.
There exist constants $c^{(2)}_{\max}$, $\widehat{c}$ such that: for any
$\delta\in(0,\frac{d\kappa}{e}]$ and $\nu\ge L_{\max}/c^{(2)}_{\max}$, with the
definition of
\begin{equation}
T\bydef\min\left\{\inf_t\{t|\widehat{f}_{\bw_0}(\bw^{(t)})-f(\bw^{(0)})\le-3\mathcal{F}\},\widehat{c}\mathcal{T}\right\}
\end{equation}
where two iterates $\{\bu^{(t)}\}$ and $\{\bw^{(t)}\}$ that are generated
by PA-PP with initial points $\{\bu^{(0)}, \bw^{(0)}\}$ satisfying
\begin{equation}
\|\bu^{(0)}-\widetilde{\bx}^{(t)}\|\le r,\quad\bw^{(0)}=\bu^{(0)}+\upsilon r \vec{\be}',\quad\upsilon\in[\delta/(2\sqrt{d}),1],\label{eq.inicond2}
\end{equation}
where $\vec{\be}'$ denotes the eigenvector of $\bT'^{-1}\bM'$ whose corresponding positive eigenvalue is minimum, if $\|\bu^{(t)}-\widetilde{\bx}^{(t)}\|\le5\widehat{c}\mathcal{S}, \forall t<T$, we will have $T<\widehat{c}\mathcal{T}$.
\end{lemma}

\paragraph{Step 3}: Quantify sufficient decrease with random perturbation.
With \leref{le.blayer31} and \leref{le.blayer32}, we can apply \leref{le.layer21} directly and obtain the following lemma.
\begin{lemma}\label{le.blayer21}
Under \asref{as1}, there exists a universal constant $c_{\max}$, for any
$\delta\in(0,d\kappa/e]$: consider a saddle point $\wbx^{(t)}$ which satisfies
\eqref{eq.csp}, let $\bx^{(0)}=\wbx^{(t)}+\xi$ where $\xi$ is generated randomly
which follows the uniform distribution over a ball with radius
$r$, and let $\bx^{(t)}$ be the
iterates of PA-PP starting from $\bx^{(0)}$. Then, when step size $\nu\ge
L_{\max}/c_{\max}$, with at least probability $1-\delta$, we have the
following for any $T\ge\mathcal{T}/c_{\max}$
\begin{equation}
f(\bx^{(T)})-f(\wbx^{(t)})\le-\mathcal{F}.
\end{equation}
\end{lemma}

Substituting $\nu=\frac{L_{\max}}{c}$,$\gamma=(L_{\max}\rho\epsilon)^{1/3}$, and
$\delta=\frac{dL_{\max}}{(L_{\max}\rho\epsilon)^{1/3}}e^{-\chi}$ in to \leref{le.blayer21}, we can obtain the following lemma immediately.
\begin{lemma}\label{le.bescape}
Under \asref{as1}, there exists a absolute constant $c_{\max}$. Let $c\le c_{\max}$,
$\chi\ge 1$, and $\eta$, $r$, $g_{\textsf{th}}$, $t_{\textsf{th}}$ calculated
as \algref{alg:p2} describes. Let $\wbx^{(t)}$ be a strict saddle point, which
satisfies
\begin{equation}
\|\nabla f(\wbx^{(t)})\|^2\le 4\nu\|\bx^{(t+1)}-\bx^{(t)}\|^2\le 4g^2_{\textsf{th}}\label{eq.condd2}
\end{equation}
and
\begin{equation}
\lambda_{\min}(\nabla^2 f(\wbx^{(t)}))\le -\gamma\notag.
\end{equation}

Let $\bx^{(t)}=\wbx^{(t)}+\xi^{(t)}$ where $\xi^{(t)}$ is generated randomly
which follows the uniform distribution over $\mathbb{B}_{\wbx^{(t)}}(r)$, and let
$\bx^{(t+t_{\textsf{th}})}$ be the iterates of PA-PP. With at least probability
$1-\frac{dL_{\max}}{(L_{\max}\rho\epsilon)^{1/3}}e^{-\chi}$, we have
\begin{equation}
f(\bx^{(t+t_{\textsf{th}})})-f(\wbx^{(t)})\le-f_{\textsf{th}}.
\end{equation}
\end{lemma}

Finally, we can get the convergence rate of PA-PP as the following.
\subsection{Proof of \coref{co.rate}}

Next, we prove the main theorem.
\begin{proof}

Submitting $\nu=\frac{L_{\max}}{c}$,$\gamma=(L_{\max}\rho\epsilon)^{1/3}$, and
$\delta=\frac{dL_{\max}}{(L_{\max}\rho\epsilon)^{1/3}}e^{-\chi}$ into the definition of $\mathcal{F},\mathcal{G}, \mathcal{T}$, we will have the following definitions.
\begin{align}
\notag
f_{\textsf{th}}\bydef&\mathcal{F}=\frac{c^5\epsilon^{2}}{L_{\max}\chi^6\mathcal{P}^2},
\\
g_{\textsf{th}}\bydef&\frac{\mathcal{G}}{2\kappa}=\frac{c^2\epsilon}{2\chi\mathcal{P}},\label{eq.defgth}
\\\notag
t_{\textsf{th}}\bydef&\frac{\mathcal{T}}{c}=\frac{L_{\max}\chi}{c^2(L_{\max}\rho\epsilon)^{\frac{1}{3}}}.
\end{align}

After applying \leref{le.layer001grad}, we know that
\begin{equation}
\|\nabla f(\bx)\|\le \frac{c}{\chi^{3}\mathcal{P}}\epsilon
\end{equation}
where $c\le 1,\chi,\mathcal{P}\ge 1$.

Similarly, at any iteration, we need to
consider two cases (we use the first iteration as an example):
\begin{enumerate}
\item In this case the gradient is large such that $\|\bx^{(1)}-\bx^{(0)}\|>g_{\textsf{th}}/\nu$: According to
\leref{le.bdescent}, we have
\begin{align}
\notag
f(\bx^{(1)})-f(\bx^{(0)})&\le-\frac{\nu}{2}\|\bx^{(1)}-\bx^{(0)}\|^2\le-\frac{\nu}{2}g^2_{\textsf{th}}
\\
&\mathop{=}\limits^{(a)}-\frac{c^5}{8\chi^6\mathcal{P}^2}\frac{\epsilon^2}{L_{\max}}\label{eq.bgradd}
\end{align}
where in $(a)$  use the definition of $g^2_{\textsf{th}}$ and $\nu\ge L_{\max}/c$.

\item The gradient is small in all block directions, namely $\|\bx^{(t+1)}-\bx^{(t)}\|^2\le g_{\textsf{th}}/\nu$: in this case, we will
add the perturbation to the iterates, and implement APP for the next
$t_{\textsf{th}}$ steps and then check the termination condition. If the
termination condition is not satisfied, we must have
\begin{equation}
f(\bx^{(t_{\textsf{th}})})-f(\bx^{(0)})\le-f_{\textsf{th}}=-\frac{c^5\epsilon^{2}}{L_{\max}\chi^6\mathcal{P}^2},
\end{equation}
which implies that the objective value in each step on average is decreased by
\begin{equation}
\frac{f(\bx^{(t_{\textsf{th}})})-f(\bx^{(0)})}{t_{\textsf{th}}}\le-\frac{c^{7}}{\chi^7\mathcal{P}^2}
\frac{\epsilon^2}{L_{\max}}\frac{(L_{\max}\rho\epsilon)^{\frac{1}{3}}}{L_{\max}}.\label{eq.bgradp}
\end{equation}
Since $\kappa=L_{\max}/(L_{\max}\rho\epsilon)^{1/3}\ge1$ and $c\le 1/3$,  we know that RHS of \eqref{eq.bgradp} is greater than RHS of \eqref{eq.bgradd}.

With the results of these two cases, we can know that if there is a large size of the gradient, we can know the decrease of the objective function value by the result of case 1, and if not, we use the result of case 2. In summary, PA-PP can have a sufficient decrease of the objective function value by $\frac{c^{7}}{\chi^7\mathcal{P}^2}\frac{\epsilon^2}{L_{\max}}\frac{(L_{\max}\rho\epsilon)^{1/3}}{L_{\max}}$ per iteration on average. This means that \algref{alg:p1} must stop within a finite number of iterations, which is
\begin{equation}
\frac{f(\bh^{(0)}_{-1},\bx^{(0)}_1)-f^*}{\frac{c^{7}}{\chi^7\mathcal{P}^2}\frac{\epsilon^2}{L_{\max}}\frac{(L_{\max}\rho\epsilon)^{1/3}}{L_{\max}}}
=\frac{\chi^{7}\mathcal{P}^2}{c^{7}}\frac{L^2_{\max}\Delta f}{\epsilon^2(L_{\max}\rho\epsilon)^{1/3}}
=\mathcal{O}\left(\frac{\Delta f\chi^7\mathcal{P}^2L^{5/3}_{\max}}{\rho^{1/3}\epsilon^{7/3}}\right)
\end{equation}
where $\Delta f\bydef f(\bh^{(0)}_{-1},\bx^{(0)}_1)-f^*$.

According to \leref{le.escape}, we know that with probability $1-\frac{dL_{\max}}{(L_{\max}\rho\epsilon)^{1/3}}e^{-\chi}$ the algorithm can give a sufficient descent with the perturbation when $\|\bx^{(t+1)}-\bx^{(t)}\|^2\le g_{\textsf{th}}/\nu$. Since the total number of perturbation we can add is at most
\begin{equation}
n'=\frac{1}{t_{\textsf{th}}}\frac{\chi^{7}\mathcal{P}^2}{c^{7}}\frac{L^2_{\max}\Delta f}{\epsilon^2(L_{\max}\rho\epsilon)^{1/3}}
=\frac{\chi^{6}\mathcal{P}^2}{c^5}\frac{L_{\max}\Delta_f}{\epsilon^2}.
\end{equation}
Using the union bound, the probability of \leref{le.escape} being satisfied for all perturbations is
\begin{equation}\label{eq.bbdoferror}
1-n'\frac{dL_{\max}}{(L_{\max}\rho\epsilon)^{\frac{1}{3}}}e^{-\chi}=1-\frac{dL_{\max}}{(L_{\max}\rho\epsilon)^{\frac{1}{3}}}e^{-\chi}\frac{\chi^{6}\mathcal{P}^2}{c^5}\frac{L_{\max}\Delta_f}{\epsilon^2}
=1-\underbrace{\frac{dL_{\max}}{(L_{\max}\rho\epsilon)^{\frac{1}{3}}}\frac{\mathcal{P}^2}{c^5}\frac{\Delta_f}{\epsilon^2}}_{\bydef\mathcal{C}'}\chi^{6}e^{-\chi}.
\end{equation}
With chosen $\chi=6\max\{\ln(\mathcal{C}'/\delta),4\}$, we have $\chi^6e^{-\chi}\le e^{-\chi/6}$, which implies $\chi^6e^{-\chi}\mathcal{C}'\le e^{-\chi/6}\mathcal{C}'\le\delta$.
\end{enumerate}
The proof is complete.
\end{proof}

\subsection{Proof of \coref{co.eiginver}}
\begin{proof}
Recall the definitions:
\begin{equation}
\bH'_u=\left[\begin{array}{cc} 0 &  \nabla^2_{12} f(\wbx^{(t)})
\\
0 &   0\end{array}\right],\quad\bH'_l=\left[\begin{array}{cc} \nabla^2_{11} f(\wbx^{(t)}) & 0
\\
\nabla^2_{21} f(\wbx^{(t)}) & \nabla^2_{22} f(\wbx^{(t)}) \end{array}\right],\quad
\end{equation}
where $\wbx^{(t)}$ is an $\epsilon$-second order stationary point, and
\begin{equation}
\bM'\bydef\bI+\eta\bH'_l\quad\bT'\bydef\bI-\eta\bH'_u.
\end{equation}
Obviously, we also have $\bH=\bH'_l+\bH'_u$.

Note that $\det(\bT')=1$, which implies that $\det(\bT'^{-1}\bM'-\lambda\bI)=\det(\bM'-\lambda\bT')$, where $\lambda$ denotes the eigenvalue. We can analyze the determinant of $\bM'-\lambda\bT'$. We have
\begin{equation}\notag
\det[\bM'-\lambda\bT']=\left[\underbrace{\begin{array}{cc}(1-\lambda)\bI+\eta\nabla^2_{11}f(\wbx^{(t)}) & \lambda\eta\nabla^2_{12} f(\wbx^{(t)})
\\\notag
\eta\nabla^2_{21}f(\wbx^{(t)})  & (1-\lambda)\bI+\eta\nabla^2_{22}f(\wbx^{(t)})
\end{array}}_{\bydef\bQ'(\lambda)}\right].
\end{equation}

It can be observed that
\begin{equation}
\notag
\bQ'(\lambda)=\left[\begin{array}{cc}\bI &  \\ & \frac{1}{\sqrt{\lambda}}\end{array}\right]\underbrace{\left[\begin{array}{cc}(1-\lambda)\bI+\eta\nabla^2_{11}f(\wbx^{(t)})  & \eta\sqrt{\lambda}\nabla^2_{12} f(\wbx^{(t)}) \\\notag
\eta\sqrt{\lambda}\nabla^2_{21}f(\wbx^{(t)})  & (1-\lambda)\bI+\eta\nabla^2_{22}f(\wbx^{(t)})
\end{array}\right]}_{\bG'(\lambda)}\left[\begin{array}{cc}\bI &  \\ & \sqrt{\lambda} \end{array}\right],
\end{equation}
meaning that  $\bQ'(\lambda)$ is similar to $\bG'(\lambda)$. Consequently, we can conclude that $\bQ'(\delta)$ has the same eigenvalues of $\bG'(\delta)$. Furthermore, since matrix $\bG'(\lambda)$ is symmetric, we know that all eigenvalues of $\bQ'(\lambda)$ and $\bG'(\lambda)$ are real. Then, we can need to show there exists $\lambda$ such that $\det(\bQ'(\lambda))=0$.

Consider $0\le\delta\le1$. We have
\begin{equation}
\bG'(1-\delta)=\left[\begin{array}{cc}\delta\bI+\eta\nabla^2_{11}f(\wbx^{(t)})  & \eta\sqrt{1-\delta}\nabla^2_{12} f(\wbx^{(t)})
\\
\eta\sqrt{1-\delta}\nabla^2_{21}f(\wbx^{(t)})  & \delta\bI+\eta\nabla^2_{22}f(\wbx^{(t)})
\end{array}\right]\label{eq.exg}.
\end{equation}

Since we know that $\bH$ and $\bG(1-\delta)$ are diagonalizable (normal matrices), then we have the following result \cite{weyl1912asymptotische} (or \cite{holbrook1992spectral}) of quantifying the difference of the eigenvalues of the two matrices
\begin{equation}
\max_{1\le i\le d}|\lambda_i(\eta\bH)-\lambda_i(\bG'(1-\delta))|\le\|\eta\bH-\bG'(1-\delta)\|\label{eq.brehg}
\end{equation}
where $\lambda_i(\bH)$ and $\lambda_i(\bG'(1-\delta))$ denote the $i$th eigenvalue of $\bH$ and $\bG'(1-\delta)$, which are listed in a decreasing order.

With the help of \eqref{eq.brehg}, we can check
\begin{align}
\notag
&\|\bG'(1-\delta)-\eta\bH\|
\\\notag
=&\left\|\delta\bI+\left[\begin{array}{cc}0 & (\sqrt{1-\delta}-1)\eta\nabla^2_{12}f(\wbx^{(t)})
\\\notag
(\sqrt{1-\delta}-1)\eta\nabla^2_{21}f(\wbx^{(t)}) & 0
\end{array}\right]\right\|
\\\notag
\le&\delta+(\sqrt{1-\delta}-1)\eta\|\bH\|+(\sqrt{1-\delta}-1)\eta\left\|\begin{array}{cc}\nabla^2_{11}f(\wbx^{(t)}) &0 \\ 0 &\nabla^2_{22}f(\wbx^{(t)})\end{array}\right\|
\\
\mathop{\le}\limits^{(a)} & \delta+(\sqrt{1-\delta}-1)(\frac{L}{L_{\max}}+1)\label{eq.bdiffeig}
\end{align}
where $(a)$ is true since we used $\eta\le c_{\max}/L_{\max}$. Also, it can be observed that when $\delta=0$, matrix $\bG'(\delta)$ is reduced to $\eta\bH$.

We know that the minimum eigenvalue of $\eta\bH$ which is $-\eta\gamma$ and the maximum difference of the eigenvalues between  $\eta\bH$ and $\bG'(\delta)$ is upper bounded by \eqref{eq.bdiffeig}. Then, we can choose a sufficient small $\delta$ such that $\bG'(\delta)$ also has a negative eigenvalue, meaning that we need to find a $\delta\in[0,1]$ such that
\begin{align}
\delta+(\sqrt{1-\delta}-1)(\frac{L}{L_{\max}}+1)<\eta\gamma. \label{eq.bvadelta}
\end{align}

In other words, if we choose
\begin{equation}\notag
\delta^*=\frac{\eta\gamma}{2}
\end{equation}
then we can conclude that $\bG'(\delta^*)$ has a negative eigenvalue which is less than $-\eta\gamma+\delta^*=-\frac{\eta\gamma}{2}$. In the following, we will check that $\delta^*$ is a valid choice, meaning that equation \eqref{eq.bvadelta} holds when $\delta^*=\frac{\eta\gamma}{2}$.

Actually, equation \eqref{eq.bvadelta} can be rewritten as
\begin{equation}
\delta+\sqrt{1-\delta}(1+\frac{L}{L_{\max}})<\eta\gamma+(1+\frac{L}{L_{\max}}),
\end{equation}
Since $\kappa=L_{\max}/\gamma\ge1$ and $\eta\le c_{\max}/L_{\max}$ where $c_{\max}\le1/2$, we have
\begin{equation}
\sqrt{1-\delta^*}=\sqrt{1-\eta\gamma/2}<1,
\end{equation}
which implies that equation \eqref{eq.bvadelta} is true with chosen $\delta^*$
Therefore, we can conclude that $\bQ'(1+\delta^*)$ has a negative eigenvalue.

When $\delta$ is large, i.e., $\delta>1$, we have
\begin{equation}
\bQ'(1-\delta)=\left[\begin{array}{cc}\bI &  \\ & \frac{j}{\sqrt{1-\delta}}\end{array}\right]\underbrace{\left[\begin{array}{cc}\delta\bI+\eta\nabla^2_{11}f(\wbx^{(t)})  & -j\eta\sqrt{1-\delta}\nabla^2_{12} f(\wbx^{(t)}) \\
\eta\sqrt{1-\delta}\nabla^2_{21}f(\wbx^{(t)})  & \delta\bI+\eta\nabla^2_{22}f(\wbx^{(t)})
\end{array}\right]}_{\bG'(1-\delta)}\left[\begin{array}{cc}\bI &  \\ & j\sqrt{1-\delta} \end{array}\right],\label{eq.grze}
\end{equation}
where $j$ denotes the imaginary number, so $\bQ'(1-\delta)$ is similar to $\bG'(1-\delta)$ when $\delta>1$. Also, we know that $\bG'(1-\delta)$ is a Hermitian matrix. It is easy to check $\bQ'(1-\delta)$ has a positive eigenvalue, since term $\delta\bI$ dominates the spectrum of matrix $\bQ'(1-\delta)$ in \eqref{eq.grze}. Considering the eigenvalue is continuous with respect to $\delta$, we can conclude there exists a $\delta$, i.e., $\widehat{\delta}'$,  such that $\bQ'(1-\widehat{\delta}')$ has a zero eigenvalue, i.e., $\det(\bQ'(1-\widehat{\delta}'))=0$ where $1-\widehat{\delta}'$ is at least as small as
\begin{equation}
1-\delta^*=1-\frac{\eta\gamma}{2},\label{eq.bgammap}
\end{equation}
meaning that $1-\wde'\le1-\frac{\eta\gamma}{2}$.
\end{proof}

In the following, we will give the proofs of \leref{le.blayer001bl}--\leref{le.blayer21} in details.

\section{Proofs of \leref{le.blayer001grad}--\leref{le.blayer21}}

\subsection{Proof of \leref{le.blayer001grad}}
\begin{proof}
First, we have
\begin{align}\notag
\|\nabla_1 f(\bx^{(t)}_1,\bx^{(t)}_2)\|^2\le& 2\|\nabla_1f(\bx^{(t+1)}_1,\bx^{(t)}_2)-\nabla_1 f(\bx^{(t)}_1,\bx^{(t)}_2)\|^2+2\|\nabla_1 f(\bx^{(t+1)}_1,\bx^{(t)}_2)\|^2
\\\notag
\mathop{\le}\limits^{(a)}&2L^2_{\max}\|\bx^{(t+1)}_1-\bx^{(t)}_1\|^2+2\|\nabla_1 f(\bx^{(t+1)}_1,\bx^{(t)}_2)\|^2
\\\notag
\mathop{\le}\limits^{\eqref{eq.upp}}&2L^2_{\max}\eta^2 \|\nabla_1 f(\bx^{(t+1)}_1,\bx^{(t)}_2)\|^2+2\|\nabla_1 f(\bx^{(t+1)}_1,\bx^{(t)}_2)\|^2
\\
\mathop{\le}\limits^{(b)}&3\|\nabla_1 f(\bx^{(t+1)}_1,\bx^{(t)}_2)\|^2\label{eq.ur1}
\end{align}
where in $(a)$ we used block-wise Lipschitz continuity, in $(b)$ we choose $\eta\le 1/(2L_{\max})$.

\begin{align}
\notag
\|\nabla_2 f(\bx^{(t)}_1,\bx^{(t)}_2)\|^2\le& 2\|\nabla_2f(\bx^{(t+1)}_1,\bx^{(t+1)}_2)-\nabla_2 f(\bx^{(t)}_1,\bx^{(t)}_2)\|^2+2\|\nabla_2 f(\bx^{(t+1)}_1,\bx^{(t+1)}_2)\|^2
\\\notag
\le&4(\|\nabla_2f(\bx^{(t+1)}_1,\bx^{(t+1)}_2)-\nabla_2 f(\bx^{(t+1)}_1,\bx^{(t)}_2)\|^2+\|\nabla_2f(\bx^{(t+1)}_1,\bx^{(t)}_2)-\nabla_2 f(\bx^{(t)}_1,\bx^{(t)}_2)\|^2)
\\\notag
&\quad+2\|\nabla_2 f(\bx^{(t+1)}_1,\bx^{(t+1)}_2)\|^2
\\\notag
\mathop{\le}\limits^{\eqref{eq.upp}}&4(L^2_{\max}\|\bx^{(t+1)}_2-\bx^{(t)}_2\|^2+\|\bx^{(t+1)}_1-\bx^{(t)}_1\|^2)+2\|\nabla_2 f(\bx^{(t+1)}_1,\bx^{(t+1)}_2)\|^2
\\
\mathop{\le}\limits^{(a)}&\|\nabla_1 f(\bx^{(t+1)}_1,\bx^{(t)}_2)\|^2+3\|\nabla_2 f(\bx^{(t+1)}_1,\bx^{(t+1)}_2)\|^2\label{eq.ur2}
\end{align}
where $(a)$ we also choose $\eta\le 1/(2L_{\max})$.

Summing \eqref{eq.ur1} and \eqref{eq.ur2}, we have
\begin{equation}
\|\nabla f(\bx^{(t)})\|^2\le\sum^2_{k=1}\|\nabla_kf(\bx^{(t)}_k)\|^2\le 4\sum^2_{k=1}\|\nabla_{k}f(\bh^{(t)}_{-k},\bx^{(t+1)}_{k})\|^2\mathop{=}\limits^{\eqref{eq.upp}}4\nu\|\bx^{(t+1)}-\bx^{(t)}\|^2
\end{equation}
where $\bh^{(t)}_{-1}=\bx^{(t)}_2$ and $\bh^{(t)}_{-2}=\bx^{(t+1)}_1$.
\end{proof}

\subsection{Proof of \leref{le.blayer001bl}}

There proof involves two parts:
\paragraph{Upper Triangular Matrix:}
Consider three different vectors $\bx$, $\by$ and $\bz$. We can have
\begin{align}
\notag
&\left\|\left[\begin{array}{cc}\nabla^2_{11}f(\bx) &  0
\\
\nabla^2_{21}f(\by) & \nabla^2_{22}f(\by)\end{array}\right]-\left[\begin{array}{cc}\nabla^2_{11}f(\bz) & 0 \\
\nabla^2_{21}f(\bz) & \nabla^2_{22}f(\bz)\end{array}\right]\right\|
\\\notag
\le&\left\|\bI_1\left(\left[\begin{array}{cc}\nabla^2_{11}f(\bx) & \nabla^2_{12}f(\bx)
\\
\nabla^2_{21}f(\bx) & \nabla^2_{22}f(\bx)\end{array}\right]-\left[\begin{array}{cc}\nabla^2_{11}f(\bz) & \nabla^2_{12}f(\bz)\\
\nabla^2_{21}f(\bz) & \nabla^2_{22}f(\bz)\end{array}\right]\right)\bI_1\right\|
\\\notag
&+\left\|\bI_2\left(\left[\begin{array}{cc}\nabla^2_{11}f(\by) & \nabla^2_{12}f(\by)
\\
\nabla^2_{21}f(\by) & \nabla^2_{22}f(\by)\end{array}\right]-\left[\begin{array}{cc}\nabla^2_{11}f(\bz) & \nabla^2_{12}f(\bz)\\
\nabla^2_{21}f(\bz) & \nabla^2_{22}f(\bz)\end{array}\right]\right)\right\|
\\\notag
\mathop{\le}\limits^{(a)}&\left\|\left[\begin{array}{cc}\nabla^2_{11}f(\bx) & \nabla^2_{12}f(\bx)
\\\notag
\nabla^2_{21}f(\bx) & \nabla^2_{22}f(\bx)\end{array}\right]-\left[\begin{array}{cc}\nabla^2_{11}f(\bz) & \nabla^2_{12}f(\bz)\\
\nabla^2_{21}f(\bz) & \nabla^2_{22}f(\bz)\end{array}\right]\right\|
+\left\|\left[\begin{array}{cc}\nabla^2_{11}f(\by) & \nabla^2_{12}f(\by)
\\\notag
\nabla^2_{21}f(\by) & \nabla^2_{22}f(\by)\end{array}\right]-\left[\begin{array}{cc}\nabla^2_{11}f(\bz) & \nabla^2_{12}f(\bz)\\
\nabla^2_{21}f(\bz) & \nabla^2_{22}f(\bz)\end{array}\right]\right\|
\\\notag
\le&\rho\left(\|\bx-\bz\|+\|\by-\bz\|\right)
\end{align}
where in $(a)$ we use
\begin{equation}
\bI_1=\left[\begin{array}{cc}\bI & 0 \\ 0 & 0\end{array}\right]\quad\quad\bI_2=\left[\begin{array}{cc}0 & 0\\0 &\bI
\end{array}\right]
\end{equation}
and $\|\bI_1\|=\|\bI_2\|=1$.

\paragraph{Lower Triangular Matrix:}
\begin{align}
\notag
&\left\|\left[\begin{array}{cc}0 & \nabla^2_{21}f(\bx)
\\
0 & 0 \end{array}\right]-\left[\begin{array}{cc} 0 & \nabla^2_{21}f(\by)\\
0 & 0\end{array}\right]\right\|
\\
&=\left\|\bI_1\left(\left[\begin{array}{cc}\nabla^2_{11}f(\bx) & \nabla^2_{12}f(\bx)
\\\notag
\nabla^2_{21}f(\bx) & \nabla^2_{22}f(\bx)\end{array}\right]-\left[\begin{array}{cc}\nabla^2_{11}f(\by) & \nabla^2_{12}f(\by)\\
\nabla^2_{21}f(\by) & \nabla^2_{22}f(\by)\end{array}\right]\right)\bI_2\right\|
\\\notag
&\mathop{\le}\limits^{(a)}\rho\|\bx-\by\|
\end{align}
where $(a)$ is true because we know $\|\bI_1\|=\|\bI_2\|=1$.

\subsection{Proof of \leref{le.bdescent}}
\begin{proof}
Under \asref{as1}, we have (descent lemma)
\begin{align}
\notag
f(\bx^{(t+1)})\le& f(\bx^{(t)})+\sum^2_{k=1}\nabla_k f(\bh^{(t)}_{-k},\bx^{(t)}_k)^{\T}(\bx^{(t+1)}_{k}-\bx^{(t)}_k)+\sum^2_{k=1}\frac{L_k}{2}\|\bx^{(t+1)}_k-\bx^{(t)}_k\|^2
\\\notag
\le& f(\bx^{(t)})+\sum^2_{k=1}\nabla_k f(\bh^{(t)}_{-k},\bx^{(t+1)}_k)^{\T}(\bx^{(t+1)}_{k}-\bx^{(t)}_k)+\sum^2_{k=1}(\nabla_k f(\bh^{(t)}_{-k},\bx^{(t)}_k)-\nabla_k f(\bh^{(t)}_{-k},\bx^{(t+1)}_k))^{\T}(\bx^{(t+1)}_{k}-\bx^{(t)}_k)
\\\notag
&\quad+\sum^2_{k=1}\frac{L_k}{2}\|\bx^{(t+1)}_k-\bx^{(t)}_k\|^2
\\\notag
\mathop{\le}\limits^{(a)}&f(\bx^{(t)})-\sum^2_{k=1}\eta\|\nabla_k f(\bh^{(t)}_{-k},\bx^{(t+1)}_k)\|^2+\sum^2_{k=1}\frac{3\eta^2L_k}{2}\|\nabla_k f(\bh^{(t)}_{-k},\bx^{(t+1)}_k)\|^2
\\\notag
\mathop{\le}\limits^{(b)}&f(\bx^{(t+1)})-\sum^2_{k=1}\frac{\eta}{2}\|\nabla_k f(\bh^{(t)}_{-k},\bx^{(t+1)}_k)\|^2
\\
=&f(\bx^{(t+1)})-\frac{\nu}{2}\|\bx^{(t+1)}-\bx^{(t)}\|^2\label{eq.bdesbcd}
\end{align}
where (a) is true because of the update rule of APP in each block
and \asref{as1} and block-wise Lipschitz continuity, in (b) we choose $\eta\le1/(3L_{\max})$ and $\nu=1/\eta$.
\end{proof}

\subsection{Proof of \leref{le.blayer31}}
\begin{proof} Without loss of generality, let $\bu^{(0)}$ be the origin, i.e.,
$\bu^{(0)}=0$. According to the APP update rule of variables, we have
\begin{align}
\bu^{(t+1)}=&\bu^{(t)}-\eta\left[\begin{array}{c}\nabla_1 f(\bu^{(t+1)}_1,\bu^{(t)}_2)
\\ \nabla_2 f(\bu^{(t+1)}_1,\bu^{(t+1)}_2)\end{array}\right]\label{eq.sumupdateo}.
\end{align}
It can be observed that the update rule of PA-PP is very similar as the one of PA-GD. The proof of this lemma is also similar as \leref{le.layer31}. We only need to replace $\nabla_1 f(\bu^{(t)}_1,\bu^{(t)}_2)$ as $\nabla_1 f(\bu^{(t+1)}_1,\bu^{(t)}_2)$ and $\nabla_2 f(\bu^{(t+1)}_1,\bu^{(t)}_2)$ as $\nabla_2 f(\bu^{(t+1)}_1,\bu^{(t+1)}_2)$, which can give us the claimed result after following the proof of \leref{le.layer31}. Hence, we ignore the repeated part with the proof of \leref{le.layer31} for simplicity of expressions.
\end{proof}

\subsection{Proof of \leref{le.blayer32}}
\begin{proof}
Let $\bu^{(0)}=0$ and define $\bv^{(t)}\bydef\bw^{(t)}-\bu^{(t)}$. According
to the assumption of \leref{le.layer32}, we know that
$\bv^{(0)}=\upsilon[\eta L_{\max}\mathcal{S}/(\kappa\log(\frac{d\kappa}{\delta})\mathcal{P}_1)]\vec{\be}'$ when
$\upsilon\in[\delta/(2\sqrt{d}),1]$. First, we define the following auxiliary function
\begin{equation}
\notag
h(\theta)\bydef\left[\begin{array}{c}\nabla_1 f(\bu^{(t+1)}_1+\theta\bv^{(t+1)}_1,\bu^{(t)}_2+\theta\bv^{(t)}_2)
\\
\nabla_2 f(\bu^{(t+1)}_1+\theta\bv^{(t+1)}_1,\bu^{(t+1)}_2+\theta\bv^{(t+1)}_2)\end{array}\right],
\end{equation}
then have
\begin{align}
\notag
& h(0)=\left[\begin{array}{c}\nabla_1 f(\bu^{(t+1)}_1,\bu^{(t)}_2)\\ \nabla_2 f(\bu^{(t+1)}_1,\bu^{(t+1)}_2)\end{array}\right],\quad
h(1)=\left[\begin{array}{c}\nabla_1 f(\bu^{(t+1)}_1+\bv^{(t+1)}_1,\bu^{(t)}_2+\bv^{(t)}_2)\\ \nabla_2 f(\bu^{(t+1)}_1+\bv^{(t+1)}_1,\bu^{(t+1)}_2+\bv^{(t+1)}_2)\end{array}\right],
\\\notag
&g(\theta)=\frac{dh(\theta)}{d\theta}=\underbrace{\left[\begin{array}{cc}\nabla^2_{11} f(\bu^{(t+1)}_1+\theta\bv^{(t+1)}_1,\bu^{(t)}_2+\theta\bv^{(t)}_2) &0 \\
\\\notag
 \nabla^2_{21} f(\bu^{(t+1)}_1+\theta\bv^{(t+1)}_1,\bu^{(t+1)}_2+\theta\bv^{(t+1)}_2) & \nabla^2_{22} f(\bu^{(t+1)}_1+\theta\bv^{(t+1)}_1,\bu^{(t+1)}_2+\theta\bv^{(t+1)}_2)  \\
\end{array}\right]}_{\tbH'^{(t)}_l(\theta)}\bv^{(t+1)}
\\
&\quad\quad\quad+\underbrace{\left[\begin{array}{cc}0 & \nabla^2_{12} f(\bu^{(t+1)}_1+\theta\bv^{(t+1)}_1,\bu^{(t)}_2+\theta\bv^{(t)}_2) \\
\\\notag
0 & 0  \\
\end{array}\right]}_{\tbH'^{(t)}_u(\theta)}\bv^{(t)},
\\\notag
& \left[\begin{array}{c}\nabla_1 f(\bw^{(t+1)}_1,\bw^{(t)}_2)\\  \nabla_2 f(\bw^{(t+1)}_1,\bw^{(t+1)}_2)\end{array}\right]=\int^1_0g(\theta)d\theta+\left[\begin{array}{c}\nabla_1 f(\bu^{(t+1)}_1,\bu^{(t)}_2)\\ \nabla_2 f(\bu^{(t+1)}_1,\bu^{(t+1)}_2)\end{array}\right].
\end{align}

Then, we consider sequence $\bw^{(t)}$, i.e.,
\begin{align}
&\bu^{(t+1)}+\bv^{(t+1)}=\bw^{(t+1)}=\bw^{(t)}-\eta\left[\begin{array}{c}\nabla_1 f(\bw^{(t+1)}_1,\bw^{(t)}_2)\\  \nabla_2 f(\bw^{(t+1)}_1,\bw^{(t+1)}_2)\end{array}\right]\label{eq.bupdateofw}
\\\notag
=&\bu^{(t)}+\bv^{(t)}-\eta\left[\begin{array}{c}\nabla_1 f(\bu^{(t+1)}_1+\bv^{(t+1)}_1,\bu^{(t)}_1+\bv^{(t)}_1)\\  \nabla_2 f(\bu^{(t+1)}_1+\bv^{(t+1)}_1,\bu^{(t+1)}_2+\bv^{(t+1)}_2)\end{array}\right]
\\
= &\bu^{(t)}+\bv^{(t)}-\eta\left[\begin{array}{c}\nabla_1 f(\bu^{(t+1)}_1,\bu^{(t)}_2)\\  \nabla_2 f(\bu^{(t+1)}_1,\bu^{(t+1)}_2)\end{array}\right]-\int^1_0g(\theta)d\theta
\\
\mathop{=}\limits^{(a)} &\bu^{(t)}+\bv^{(t)}-\eta\left[\begin{array}{c}\nabla_1 f(\bu^{(t+1)}_1,\bu^{(t)}_2)\\  \nabla_2 f(\bu^{(t+1)}_1,\bu^{(t+1)}_2)\end{array}\right]-\eta\tDelta'^{(t)}_u \bv^{(t)} - \bH'_u\bv^{(t)}-\eta\tDelta'^{(t)}_l \bv^{(t+1)} - \eta\bH'_l\bv^{(t+1)}
\label{eq.bdynofv}
\end{align}
where in $(a)$  we used the following definitions
\begin{align}\notag
\tDelta'^{(t)}_u\bydef\int^1_0\tbH'^{(t)}_u(\theta)d\theta-\bH'_u,
\\\notag
\tDelta'^{(t)}_l\bydef\int^1_0\tbH'^{(t)}_l(\theta)d\theta-\bH'_l,
\end{align}
and
\begin{equation}
\bH'_u=\left[\begin{array}{cc} 0 &  \nabla^2_{12} f(\wbx^{(t)})
\\
0 &   0\end{array}\right]\quad\bH'_l=\left[\begin{array}{cc} \nabla^2_{11} f(\wbx^{(t)}) & 0
\\
\nabla^2_{21} f(\wbx^{(t)}) & \nabla^2_{22} f(\wbx^{(t)}) \end{array}\right].
\end{equation}
Obviously, $\bH=\bH'_l+\bH'_u$.

\paragraph{Dynamics of $\bv^{(t)}$:}

Since the first two terms at RHS of \eqref{eq.bdynofv} combined with $\bu^{(t)}$ at LHS of \eqref{eq.bdynofv} are exactly the same as \eqref{eq.sumupdateo}. It can be observed that equation \eqref{eq.bdynofv} gives the dynamic of
$\bv^{(t)}$, i.e.,
\begin{equation}\label{eq.brev}
\bv^{(t+1)}=\bv^{(t)}-\eta\tDelta'^{(t)}_u \bv^{(t)} - \eta\bH'_u\bv^{(t)}-\eta\tDelta'^{(t)}_l \bv^{(t+1)} - \eta\bH'_l\bv^{(t+1)},
\end{equation}
which can be equivalently expressed by
\begin{equation}
\underbrace{(\bI+\eta\bH'_l)}_{\bydef\bM'}\bv^{(t+1)}=\underbrace{(\bI-\eta\bH'_u)}_{\bydef\bT'}\bv^{(t)}-\eta\tDelta'^{(t)}_l \bv^{(t+1)}-\eta\tDelta'^{(t)}_u \bv^{(t)}\label{eq.biterav}.
\end{equation}

It is worth noting that matrix $\bT'$ is an upper triangular matrix where the diagonal entries are all 1s, so it is invertible. Taking the inverse of $\bT'$ on both sides of \eqref{eq.biterav}, we can obtain
\begin{equation}
\bT'^{-1}\bM'\bv^{(t+1)}\mathop{=}\limits^{\eqref{eq.bdynofv}}\bv^{(t)}-\bT'^{-1}\eta\tDelta'^{(t)}_l \bv^{(t+1)}-\bT'^{-1}\eta\tDelta'^{(t)}_u\bv^{(t)}.\label{eq.biterav2}
\end{equation}

Let $\mathbb{P}'_{\texttt{left}}$ denote the projection operator that projects the vector onto the space spanned by the eigenvector of $\bT'^{-1}\bM$ whose corresponding positive eigenvalue is minimum. Taking the projection on both sides of \eqref{eq.biterav2}, we have
\begin{equation}
\mathbb{P}'_{\texttt{left}}(\bT'^{-1}\bM')\bv^{(t+1)}+\mathbb{P}'_{\texttt{left}}\bT'^{-1}\eta\tDelta'^{(t)}_l \bv^{(t+1)}=\mathbb{P}'_{\texttt{left}}\bv^{(t)}-\mathbb{P}'_{\texttt{left}}\bT'^{-1}\eta\tDelta'^{(t)}_u\bv^{(t)}.\label{eq.bkeyre}
\end{equation}

\paragraph{Relationship of the Norm of $\bv^{(t)}$ Projected onto the Two Subspaces:}

Let $\phi^{(t)}$ denote the norm of $\bv^{(t)}$ projected onto the space spanned by the eigenvector of $\bT'^{-1}\bM'$ whose positive minimum eigenvalue of $\bM'^{-1}\bT'$ is $1-\wde'>0$ and $\theta^{(t)}$ denote the norm of $\bv^{(t)}$ projected onto the remaining space. From \eqref{eq.bkeyre}, we can have
\begin{align}
(1-\wde')\phi^{(t+1)}\mathop{\ge}\limits^{(a)}\phi^{(t)}-\eta\|\bT'^{-1}\|\|\tDelta'^{(t)}_l\|\|\bv^{(t+1)}\|-\eta\|\bT'^{-1}\|\|\tDelta'^{(t)}_u\|\|\bv^{(t)}\|,\label{eq.brecurphi}
\\
(1-\wde')\theta^{(t+1)}\le\theta^{(t)}+\eta\|\bT'^{-1}\|\|\tDelta'^{(t)}_l\|\|\bv^{(t+1)}\|+\eta\|\bT'^{-1}\|\|\tDelta'^{(t)}_u\|\|\bv^{(t)}\|.\label{eq.brecurtheta}
\end{align}
where $(a)$ is true because we applied the triangle inequality since $\eta$ is sufficiently small.

Since
$\|\bw^{(0)}-\wbx^{(t)}\|\le\|\bu^{(0)}-\wbx^{(t)}\|+\|\bv^{(0)}\|\le 2r$,
we can apply \leref{le.blayer31}. Then, we know
$\|\bw^{(t)}-\wbx^{(t)}\|\le 5\widehat{c}\mathcal{S},\forall t< T$. According to the
assumptions of \leref{le.blayer32}, we have $\|\bu^{(t)}-\wbx^{(t)}\|\le
5\widehat{c}\mathcal{S}$, and
\begin{equation}\label{eq.bsizeofv}
\|\bv^{(t)}\|=\|\bw^{(t)}-\bu^{(t)}\|\le\|\bu^{(t)}-\wbx^{(t)}\|+\|\bw^{(t)}-\wbx^{(t)}\|\le 10\widehat{c}\mathcal{S}.
\end{equation}

From \eqref{eq.bdofdiffu}, we know that
\begin{equation}\notag
\|\bw^{(t+1)}-\bw^{(t)}\|\le\frac{4.3\eta\mathcal{G}}{\kappa}=
\frac{4.3\eta^3L^3_{\max}\frac{\gamma}{\rho}}{\kappa^2\log^3\frac{d\kappa}{\delta}\mathcal{P}}\le\mathcal{S},
\end{equation}
where we choose $\eta\le c_{\max}/L_{\max}$ and $c_{\max}=1/10$. Similarly, we also have $\|\bu^{(t+1)}-\bu^{(t)}\|\le\mathcal{S}$.

Then, we need to quantify the upper bounds of $\|\bM'^{-1}\|$, $\|\bv^{(t+1)}\|$, $\|\tDelta'^{(t)}_u\|$ and  $\|\tDelta'^{(t)}_l\|$.
\begin{enumerate}

\item Upper bound of $\|\bM'^{-1}\|$: applying the steps of deriving \eqref{eq.minv}, we can quantify the inverse of matrix $\bT'$ as follows
\begin{align}
\notag
\|\bT'^{-1}\|\le&1+\eta\|\bH'_u\|=1+\eta\|\bH'^{\T}_u\|
\\\notag
=&1+\|\eta\bH\odot\bD-\eta\bH_d\|
\\\notag
<&2(1+\frac{L\log(2d)}{L_{\max}}).
\end{align}

\item Relation between $\|\bv^{(t)}\|$ and $\|\bv^{(t+1)}\|$:
We also know that
\begin{align}
\notag
\|\bv^{(t+1)}\|^2=&\|\bw^{(t+1)}-\bu^{(t+1)}\|^2=\left\|\bw^{(t)}-\eta\left[\begin{array}{c}\nabla_1 f(\bw^{(t+1)}_1,\bw^{(t)}_2)\\  \nabla_2 f(\bw^{(t+1)}_1,\bw^{(t+1)}_2)\end{array}\right] -\left(\bu^{(t)}-\eta\left[\begin{array}{c}\nabla_1 f(\bu^{(t+1)}_1,\bu^{(t)}_2)\\ \nabla_2 f(\bu^{(t+1)}_1,\bu^{(t+1)}_2)\end{array}\right]\right)\right\|^2
\\\notag
\le&2\|\bv^{(t)}\|^2+4\eta^2\left\|\left[\begin{array}{c}\nabla_1 f(\bw^{(t+1)}_1,\bw^{(t)}_2)\\ \nabla_2 f(\bw^{(t+1)}_1,\bw^{(t+1)}_2)\end{array}\right]-\left[\begin{array}{c}\nabla_1 f(\bu^{(t+1)}_1,\bw^{(t)}_2)\\  \nabla_2 f(\bu^{(t+1)}_1,\bw^{(t+1)}_2)\end{array}\right]\right\|^2
\\\notag
\quad&+4\eta^2\left\|\left[\begin{array}{c}\nabla_1 f(\bu^{(t+1)}_1,\bw^{(t)}_2)\\  \nabla_2 f(\bu^{(t+1)}_1,\bw^{(t+1)}_2)\end{array}\right]-\left[\begin{array}{c}\nabla_1 f(\bu^{(t+1)}_1,\bu^{(t)}_2)\\ \nabla_2 f(\bu^{(t+1)}_1,\bu^{(t+1)}_2)\end{array}\right]\right\|^2
\\
\mathop{\le}\limits^{(a)}& 2\|\bv^{(t)}\|^2+8\eta^2L^2_{\max}\|\bv^{(t)}_1\|^2+4\eta^2L^2_{\max}(\|\bv^{(t+1)}_2\|^2+\|\bv^{(t)}_2\|^2)\label{eq.brevv}
\end{align}
where $(a)$ is true due to Lipschitz continuity.

We can express \eqref{eq.brevv} as
\begin{equation}\notag
(1-4\eta^2L^2_{\max})\|\bv^{(t+1)}\|\le(2+8\eta^2L^2_{\max})\|\bv^{(t)}\|^2,
\end{equation}
which implies
\begin{equation}
\|\bv^{(t+1)}\|\le\sqrt{\frac{2+\frac{8}{100}}{1-\frac{4}{100}}}\|\bv^{(t)}\|<\sqrt{2.2}\|\bv^{(t)}\|<1.5\|\bv^{(t)}\|\label{eq.brelavv}
\end{equation}
where we choose $\eta\le c_{\max}/L_{\max}$ and $c_{\max}=1/10$.

\item Upper bound of $\|\tDelta'^{(t)}_l\|$: applying \leref{le.blayer001bl}, we can also get the upper bound of $\|\tDelta'^{(t)}_l\|$, i.e.,
\begin{align}
\notag
\|(\tDelta'^{(t)}_l)\|\le&\int^1_0\|\tbH'^{(t)}_l(\theta)-\bH'_l\|d\theta
\\\notag
\mathop{\le}\limits^{\eqref{eq.bdeltatbd}}&\int^1_0\rho\left(\|\bu^{(t+1)}+\theta\bv^{(t+1)}-\wbx^{(t)}\|+\left\|\left[\begin{array}{c}\bu^{(t+1)}_1+\theta\bv^{(t+1)}_1 \\ \bu^{(t)}_2+\theta\bv^{(t)}_2\end{array}\right]-\wbx^{(t)}\right\|\right)d\theta
\\\notag
\le&\int^1_0\rho\left(2\|\bu^{(t+1)}+\theta\bv^{(t+1)}-\wbx^{(t)}\|+\|\bu^{(t)}+\theta\bv^{(t)}-\wbx^{(t)}\|\right)d\theta
\\\notag
\le&\rho(2\|\bu^{(t+1)}-\wbx^{(t)}\|+\|\bu^{(t)}-\wbx^{(t)}\|)+\rho\int^1_0\theta(2\|\bv^{(t+1)}\|+\|\bv^{(t)}\|)d\theta
\\\notag
\le&\rho\left(2\|\bu^{(t+1)}-\bu^{(t)}\|+2\|\bu^{(t)}-\wbx^{(t)}\|+\|\bu^{(t)}-\wbx^{(t)}\|)+0.5\|\bv^{(t+1)}\|+0.5\|\bv^{(t)}\|\right)
\\\notag
\mathop{\le}\limits^{\eqref{eq.brelavv}} &\rho\left(2\|\bu^{(t+1)}-\bu^{(t)}\|+3\|\bu^{(t)}-\wbx^{(t)}\|+1.25\|\bv^{(t)}\|\right)
\\\notag
\le & \rho (2+27.5\widehat{c})\mathcal{S}.
\end{align}

\item Upper bound of $\|\tDelta'^{(t)}_u\|$: according to $\rho$-Hessian Lipschitz continuity and \leref{le.blayer001bl}, we have the size of $\tDelta'^{(t)}_u$ as the following.
\begin{align}
\notag
\|(\tDelta'^{(t)}_u)\|\le&\int^1_0\|\tbH'^{(t)}_u(\theta)-\bH'_u\|d\theta
\\
\mathop{\le}\limits^{\eqref{eq.bdeltatbd2}}&\int^1_0\rho\left\|\left[\begin{array}{c}\bu^{(t+1)}_1+\theta\bv^{(t+1)}_1 \\ \bu^{(t)}_2+\theta\bv^{(t)}_2\end{array}\right]-\wbx^{(t)}\right\|d\theta
\\\notag
\le&\int^1_0\rho(\|\bu^{(t)}+\theta\bv^{(t)}-\wbx^{(t)}\|+\|\bu^{(t+1)}+\theta\bv^{(t+1)}-\wbx^{(t)})\|d\theta
\\\notag
\le&\rho(\|\bu^{(t+1)}-\wbx^{(t)}\|+\|\bu^{(t)}-\wbx^{(t)}\|)+\rho\int^1_0\theta(\|\bv^{(t+1)}\|+\|\bv^{(t)}\|)d\theta
\\\notag
\mathop{\le}\limits^{\eqref{eq.brelavv}} &\rho\left(\|\bu^{(t+1)}-\bu^{(t)}\|+2\|\bu^{(t)}-\wbx^{(t)}\|+1.25\|\bv^{(t)}\|\right)
\\\notag
\le&\rho(1+22.5\widehat{c})\mathcal{S}.
\end{align}

\end{enumerate}

With the bounds of $\|\bv^{(t+1)}\|$, $\|\tDelta'^{(t)}_u\|$, $\|\tDelta'^{(t)}_l\|$ and relation between $\|\bv^{(t+1)}\|$ and $\|\bv^{(t)}\|$, we can further simply \eqref{eq.brecurphi} and \eqref{eq.brecurtheta} as follows,
\begin{align}
\notag
(1-\wde')\phi^{(t+1)}\mathop{\ge}\limits^{\eqref{eq.brecurphi}}\phi^{(t)}-\eta(1.5\|\tDelta'^{(t)}_l\|+\|\tDelta'^{(t)}_u\|)\|\bT'^{-1}\|\sqrt{(\phi^{(t)})^2+(\theta^{(t)})^2},
\\\notag
(1-\wde')\theta^{(t+1)}\mathop{\le}\limits^{\eqref{eq.brecurtheta}}\theta^{(t)}+\eta(1.5\|\tDelta'^{(t)}_l\|+\|\tDelta'^{(t)}_u\|)\|\bT'^{-1}\|\sqrt{(\phi^{(t)})^2+(\theta^{(t)})^2},
\end{align}
since $\|\bv^{(t)}\|=\sqrt{(\phi^{(t)})^2+(\theta^{(t)})^2}$.

Consequently, we can arrive at
\begin{align}
(1-\wde')\phi^{(t+1)}\ge\phi^{(t)}-\mu\sqrt{(\phi^{(t)})^2+(\theta^{(t)})^2},\label{eq.bphii}
\\
(1-\wde')\theta^{(t+1)}\le\theta^{(t)}+\mu\sqrt{(\phi^{(t)})^2+(\theta^{(t)})^2},\label{eq.bthei}
\end{align}
where $\mu$ is the upper bound of term $\eta(1.5\|\tDelta'^{(t)}_l\|+\|\tDelta'^{(t)}_u\|)\|\bT'^{-1}\|$ and can be obtained by
\begin{equation}
\mu\bydef\eta\rho\mathcal{S}\mathcal{P}(4+62\widehat{c})\label{eq.bdefofmu}.
\end{equation}

\paragraph{Quantifying the Norm of $\bv^{(t)}$ Projected at Different Subspaces:}
Then, we will use mathematical induction to prove
\begin{equation}\label{eq.bphinorm}
\theta^{(t)}\le 4\mu t\phi^{(t)}.
\end{equation}
It is true when $t=0$ since $\|\theta^{(0)}\|\mathop{=}\limits^{\eqref{eq.inicond2}}0$.

Assuming that equation $\eqref{eq.bphinorm}$ is
true at the $t$th iteration, we need to prove
\begin{equation}
\theta^{(t+1)}\le 4\mu (t+1)\phi^{(t+1)}.\label{eq.bgoa}
\end{equation}

Applying \eqref{eq.bphii} into RHS of \eqref{eq.bgoa}, we have
\begin{equation}\label{eq.brephi}
4\mu(t+1)\phi^{(t+1)}\ge\frac{4\mu(t+1)}{1-\wde'}\left(\phi^{(t)}-\mu\sqrt{(\phi^{(t)})^2+(\theta^{(t)})^2}\right),
\end{equation}
and substituting \eqref{eq.bthei} into LHS of \eqref{eq.bgoa}, we have
\begin{equation}\label{eq.bretheta}
\theta^{(t+1)}\le \frac{(4\mu t\phi^{(t)})+\mu\sqrt{(\phi^{(t)})^2+(\theta^{(t)})^2}}{1-\wde'}.
\end{equation}

Then, our goal is to prove RHS of \eqref{eq.brephi} is greater
than RHS of \eqref{eq.bretheta}. After some manipulations, it is sufficient to show
\begin{equation}
\left(1+4\mu(t+1)\right)\left(\sqrt{(\phi^{(t)})^2+(\theta^{(t)})^2}\right)\le4\phi^{(t)}.
\end{equation}

In the following, we will show that the above relation is true.
\paragraph{First step}:
We know that
\begin{equation}\label{eq.bbdofmu}
4\mu(t+1)\le4\mu T\mathop{\le}\limits^{\eqref{eq.bdefofmu}}4\eta\rho\mathcal{S}\mathcal{P}(4+62\widehat{c})\widehat{c}\mathcal{T}\mathop{\le}\limits^{\eqref{eq.bsandt}\eqref{eq.bdefofmu}}\frac{4\widehat{c}\eta^2 L^2_{\max}(4+62\widehat{c})}{\kappa\log(\frac{d\kappa}{\delta})}\mathop{\le}\limits^{(a)}1
\end{equation}
where $(a)$ is true because we choose $c'_{\max}=1/(2\widehat{c}(4+62\widehat{c}))$ and
$\eta\le c'_{\max}/L_{\max}$.

\paragraph{Second step}:
Also, we know that
\begin{equation}\notag
4\phi^{(t)}\ge2\sqrt{2(\phi^{(t)})^2}\mathop{\ge}\limits^{\eqref{eq.bphinorm},\eqref{eq.bbdofmu}}(1+4\mu(t+1))\sqrt{(\phi^{(t)})^2+(\theta^{(t)})^2}.
\end{equation}

With the above two steps, we have $\theta^{(t+1)}\le 4\mu(t+1)\phi^{(t+1)}$, which completes the induction.

\paragraph{Recursion of $\phi^{(t)}$}:Using \eqref{eq.bphinorm}, we have $\theta^{(t)}\mathop{\le}\limits^{\eqref{eq.bphinorm}} 4\mu t\phi^{(t)}\mathop{\le}\limits^{\eqref{eq.bbdofmu}} \phi^{(t)}$, and have
\begin{equation}
\notag
(1-\wde')\phi^{(t+1)}\mathop{\ge}\limits^{\eqref{eq.bphii}}\phi^{(t)}-\mu\sqrt{(\phi^{(t)})^2+(\theta^{(t)})^2},
\end{equation}
which implies
\begin{align}\notag
\phi^{(t+1)}\mathop{\ge}\limits^{(a)}&\frac{1}{1-\wde'}\left(\phi^{(t)}-\mu\sqrt{(\phi^{(t)})^2+(\theta^{(t)})^2}\right)
\\\notag
\mathop{\ge}\limits^{(b)}&\frac{1}{1-\frac{\eta\gamma}{2}}\left(\phi^{(t)}-\mu\sqrt{(\phi^{(t)})^2+(\theta^{(t)})^2}\right)
\\
\mathop{\ge}\limits^{(c)}&\frac{1-\frac{\gamma^2\eta^2}{4}}{1-\frac{\eta\gamma}{2}}\phi^{(t)}=(1+\frac{\eta\gamma}{2})\phi^{(t)}\label{eq.brelation}
\end{align}
where $(a)$ is true because $1-\wde'>0$,  in $(b)$ we used \coref{co.eiginver}, i.e., $0<1-\wde'\le1-\frac{\eta\gamma}{2}$, and $(c)$ is true because $\theta^{(t)}\le\phi^{(t)}$ and
\begin{equation}
\notag
\mu=\eta\rho\mathcal{S}\mathcal{P}(4+62\widehat{c})
\mathop{\le}\limits^{\eqref{eq.seg}}\gamma^2\eta^2\frac{\eta L_{\max}(4+62\widehat{c})}{\log^{2}(\frac{d\kappa}{\delta})}
\mathop{\le}\limits^{(a)}\frac{\gamma^2\eta^2}{4\sqrt{2}}
\end{equation}
where in $(a)$ we choose $c''_{\max}= 1/(4\sqrt{2}(4+62\widehat{c}))$ and $\eta\le c''_{\max}/L_{\max}$.

\paragraph{Quantifying Escaping Time:}

From  \eqref{eq.bsizeofv}, we have
\begin{align}
\notag
10\mathcal{S}\widehat{c}\ge&\|\bv^{(t)}\|\ge\phi^{(t)}
\mathop{\ge}\limits^{\eqref{eq.brelation}}(1+\frac{\gamma\eta}{2})^{t}\phi^{(0)}
\mathop{\ge}\limits^{(a)}(1+\frac{\gamma\eta}{2})^{t}\frac{\delta}{2\sqrt{d}}\frac{\eta L_{\max}\mathcal{S}}{\kappa}\log^{-1}(\frac{d\kappa}{\delta})
\\
\mathop{\ge}\limits^{(b)}&(1+\frac{\gamma\eta}{2})^{t}\frac{\delta}{2\sqrt{d}}\frac{c\mathcal{S}}{\kappa}\log^{-1}(\frac{d\kappa}{\delta}) \quad\forall t<T \label{eq.bupoft}
\end{align}
where in $(a)$ we use condition $\upsilon\in[\delta/(2\sqrt{d}),1]$, in $(b)$ we used $\eta =c/L_{\max}$.

Since \eqref{eq.bupoft} is true for all $t<T$,  we can have
\begin{align}
\notag
T-1\le&\frac{\log(20\frac{\widehat{c}}{c}(\frac{\kappa\sqrt{d}}{\delta})\log(\frac{d\kappa}{\delta}))}{\log(1+\frac{\eta\gamma}{2})}
\mathop{<}\limits^{(a)}\frac{4\log(20(\frac{\sqrt{d}\kappa}{\delta})\frac{\widehat{c}}{c}\log(\frac{d\kappa}{\delta}))}{\eta\gamma}
\\
\mathop{<}\limits^{(b)}&\frac{4\log(20(\frac{d\kappa }{\delta})^2\frac{\widehat{c}}{c})}{\eta\gamma}
\mathop{<}\limits^{(c)}4(2+\log(20\frac{\widehat{c}}{c}))\mathcal{T}\label{eq.bbdt}
\end{align}
where $(a)$ comes from inequality $\log(1+x)>x/2$ when $x<1$,  in $(b)$ we used relation $\log(x)<x, x>0$,   and $(c)$ is true because $\delta\in(0,\frac{d\kappa}{e}]$ and $\log(d\kappa/\delta)>1$.

From \eqref{eq.bbdt}, we know that
\begin{equation}
T<4(2+\log(20\frac{\widehat{c}}{c}))\mathcal{T}+1\mathop{<}\limits^{(a)}4(2\frac{1}{4}+\log(20\frac{\widehat{c}}{c})\mathcal{T}\label{eq.bdtt}
\end{equation}
where $(a)$ is true due to the fact that $\eta L_{\max}\ge1$ and $\log(d\kappa/\delta)>1$ so we know $\mathcal{T}\ge1$.

Applying the proof from \eqref{eq.chat} to \eqref{eq.bdhatc}, we can also conclude that there exists a universal $\widehat{c}$ such that \eqref{eq.bdtt} holds. The proof is complete.

\end{proof}

\subsection{Proof of \leref{le.blayer21}}

First, after the random perturbation, the objective function value in the worst case is increased at most by
\begin{align}
\notag
f(\bu^{(0)})-f(\wbx^{(t)})\le& \sum^2_{k=1}\nabla_kf(\tbh^{(t)}_{-k},\wbx^{(t)}_k)^{\T}\xi_k+\frac{L_k}{2}\|\xi_k\|^2
\\\notag
\le&\sum^2_{k=1}\left(\nabla_kf(\tbh^{(t)}_{-k},\wbx^{(t)}_k)-\nabla_kf(\tbh^{(t)}_{-k},\wbx^{(t+1)}_k)\right)^{\T}\xi_k+\sum^2_{k=1}\nabla_kf(\tbh^{(t)}_{-k},\wbx^{(t+1)}_k)^{\T}\xi_k+\frac{L_k}{2}\|\xi_k\|^2
\\\notag
\le&\sum^2_{k=1}L_{\max}\left\|\bx^{(t+1)}_k-\bx^{(t)}_k\right\|\|\xi_k\|+\sum^2_{k=1}\|\nabla_kf(\tbh^{(t)}_{-k},\wbx^{(t+1)}_k)\|\|\xi_k\| +\frac{L_{\max}}{2}\|\xi\|^2
\\\notag
\mathop{\le}\limits^{(a)}&1.25\sum^2_{k=1}\|\nabla_kf(\tbh^{(t)}_{-k},\wbx^{(t+1)}_k)\|\|\xi_k\| +\frac{L_{\max}}{2}\|\xi\|^2
\\\notag
\mathop{\le}\limits^{(b)}&1.25\|\xi\|\sqrt{\sum^2_{k=1}2\|\nabla_kf(\tbh^{(t)}_{-k},\wbx^{(t+1)}_k)\|^2}+\frac{L_{\max}}{2}\|\xi\|^2
\\
\mathop{\le}\limits^{(c)}& 1.25\frac{\mathcal{G}}{\kappa}\frac{\eta L_{\max}\mathcal{S}}{\kappa\log(\frac{d\kappa}{\delta})\mathcal{P}}+\frac{L_{\max}}{2}(\frac{\eta L_{\max}\mathcal{S}}{\kappa\log(\frac{d\kappa}{\delta})\mathcal{P}})^2\le\frac{3}{2}\mathcal{F}\label{eq.bsufin}
\end{align}
where $\bu^{(0)}$ is a vector that follows uniform distribution within the ball $\mathbb{B}^{(d)}_{\wbx^{(t)}}(r)$, $\mathbb{B}^{(d)}_{\wbx^{(t)}}$ denotes the $d$-dimensional ball centered at $\wbx^{(t)}$ with radius $r$, $\xi_k$ represents the $k$th block of the vector which is the difference between random generated vector $\bu^{(0)}$ and saddle point $\wbx^{(t)}$, and in $(a)$ we choose $\eta\le1/(4L_{\max})$ and $(b)$ is true because $\xi\bydef[\xi_1,\ldots,\xi_K]$, $\|\xi_k\|\le\|\xi\|,\forall k$, and in $(c)$ we used $\kappa>1$, $\log(d\kappa/\delta)>1$, $\mathcal{P}\ge2$ and \conref{cond2} where $g_{\textsf{th}}$ is defined in \eqref{eq.defgth}.

Then, the rest of proof of \leref{le.blayer21} is the same as the rest of \leref{le.layer21}, therefore ignored for simplicity.

\newpage
\section{Numerical Results}

\subsection{Proof of \leref{le.simup}}
\begin{proof}
Consider function
\begin{equation}
f(\bx)=\bx^{\T}\bA\bx+\frac{1}{4}\|\bx\|^4_4\label{eq.siobj2}
\end{equation}
where $\bx\in\mathcal{S}$, $\mathcal{S}=\{\bx|\|\bx\|^2\le\tau\}$ and $\tau\ge\lambda_{\max}(\bA)$.

\paragraph{To prove L-smooth Lipschitz continuity}:
\begin{align}
\notag
\|\nabla f(\bx)-\nabla f(\by)\|=&\left\|2(\bA\bx-\bA\by)+\left[\begin{array}{c}x^3_1-y^3_1\\\vdots\\x^3_d-y^3_d\end{array}\right]\right\|,\quad\forall\bx,\by\in\mathcal{S}
\\\notag
\le&2\lambda_{\max}(\bA)\|\bx-\by\|+\left\|\left[\begin{array}{c}(x_1-y_1)(x^2_1+x_1y_1+y^2_1)\\\vdots\\(x_d-y_d)(x^2_d+x_dy_d+y^2_d)\end{array}\right]\right\|
\\\notag
\mathop{\le}\limits^{(a)}&2\tau\|\bx-\by\|+3\tau\|\bx-\by\|\le5\tau\|\bx-\by\|
\end{align}
where $x_i$ denotes the $i$th entry of vector $\bx$, and $(a)$ is true because
\begin{equation}
x^2_i\le\tau,\quad y^2_i\le\tau,\quad x_iy_i\le (x^2_i+y^2_i)/2\le\tau,\forall i.\label{eq.scain}
\end{equation}

\paragraph{To prove block-wise Lipschitz continuity}:
Without loss of generality, consider first block $\bx_1\in\mathcal{S}'$ where $\mathcal{S}'=\{\bx_1|\|\bx_1\|^2\le\tau',\bx_1\in\mathbb{R}^{d'\times 1}\}$ and  $d'$ denotes the dimension of $\bx_1$. Consider $\tau'\ge\lambda_{\max}(\bA')$ where $\bA'\in\mathbb{R}^{d'\times d'}$ is the leading principal minor of matrix $\bA$ of order $d'$. Obviously, we have $\tau'\le\tau$.
\begin{align}
\notag
\|\nabla_1 f(\bx_{-1},\bx_1)-\nabla_1 f(\bx_{-1},\bx'_1)\|=&\left\|2\bI'_1\left(\bA\left[\begin{array}{c}\bx_1\\\bx_{-1}\end{array}\right]-\bA\left[\begin{array}{c}\bx'_1\\\bx'_{-1}\end{array}\right]\right)+\left[\begin{array}{c}x^3_1-x'^3_1\\\vdots\\x^3_{d'}-x'^3_{d'}\end{array}\right]\right\|,\quad\forall\bx,\bx'\in\mathcal{S}'
\\\notag
\le&2\|\bI'_1\left(\bA\left[\begin{array}{c}\bx_1\\\bx_{-1}\end{array}\right]-\bA\left[\begin{array}{c}\bx'_1\\\bx_{-1}\end{array}\right]\right)\|+\left\|\left[\begin{array}{c}(x_1-x'_1)(x^2_1+x_1x'_1+x'^2_1)\\\vdots\\(x_{d'}-x'_{d'})(x^2_{d'}+x_{d'}x'_{d'}+x'^2_{d'})\end{array}\right]\right\|
\\\notag
\mathop{\le}\limits^{(a)}&2\lambda_{\max}(\bA')\|\bx_1-\bx'_1\|+3\tau'\|\bx_1-\by_1\|
\\\notag
\le&5\tau'\|\bx_1-\bx'_1\|,\quad\forall \bx,\bx'
\end{align}
where $(a)$ is true because we used $\bI'_1\bydef\left[\begin{array}{cc}\bI_{d'} & 0 \\ 0 & 0\end{array}\right]$ which selects the first $d'$ rows of $\bA\left(\left[\begin{array}{c}\bx_1\\\bx_{-1}\end{array}\right]-\left[\begin{array}{c}\bx'_1\\\bx_{-1}\end{array}\right]\right)$.

\paragraph{To prove Hessian Lipschitz continuity}:
\begin{align}\notag
\|\nabla^2 f(\bx)-\nabla^2 f(\by)\|=&3\left\|\begin{array}{ccc}x^2_1-y^2_1 & \cdots & 0 \\ \vdots & \ddots & \vdots \\ 0 & \cdots & x^2_d-y^2_d\end{array}\right\|
\\\notag
\le&6\sqrt{\tau}\left\|\begin{array}{ccc}x_1-y_1 & \cdots & 0 \\ \vdots & \ddots & \vdots \\ 0 & \cdots & x_d-y_d\end{array}\right\|=6\sqrt{\tau}\|\bx-\by\|
\end{align}
where $(a)$ is true because $x_i+y_i\le\sqrt{(x_i+y_i)^2}=\sqrt{x^2_1+2x_iy_i+y^2_i}\mathop{\le}\limits^{\eqref{eq.scain}}2\sqrt{\tau},\forall i$.

\end{proof}
%

\subsection{Additional Simulation}
\begin{figure}[ht]
\vskip 0.2in
\begin{center}
\centerline{\includegraphics[width=0.45\columnwidth]{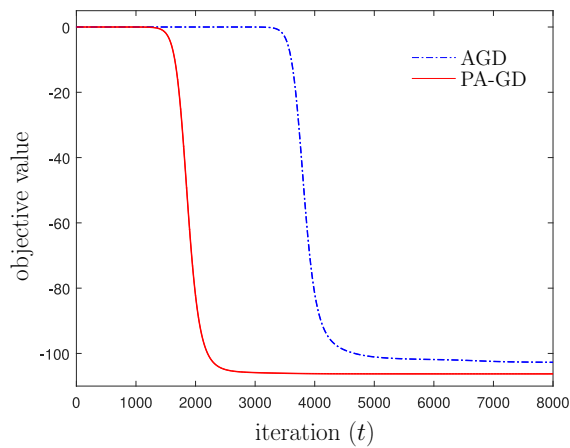}}
\caption{Convergence comparison between AGD and PA-GD, where $d=100$, $\epsilon=10^{-4}$, $g_{\textsf{th}}=\epsilon/10$, $\eta=1\times 10^{-3}$, $t_{\textsf{th}}=10/\epsilon^{1/3}$, $r=\epsilon/10$.}
\label{fig:4}
\end{center}
\vskip -0.2in
\end{figure}

\paragraph{Random matrix $\bA$}: we also test the algorithms with a randomly generated symmetric matrix $\bA$ by the following steps: 1) randomly generate a diagonal matrix $\bD$ whose entries follow \emph{i.i.d.} Gaussian distribution with zero mean and variance two; 2) generate an orthogonal matrix $\bU\in\mathbb{R}^{d\times d}$; 3) obtain matrix $\bA=\bU\bD\bU^{\T}$. We initialize the PA-GD/AGD algorithms around the saddle point which is at the origin. The results are shown in \figref{fig:4} where $d=100$. It can be observed that PA-GD can still escape from the strict saddle point faster than ordinary AGD, illustrating the benefit of adding the random perturbation when the gradient size is small.
 %

\end{document}